\documentclass[12pt,reqno]{amsart}
\usepackage{hyperref}
\usepackage{enumerate,amssymb}
\usepackage{amsopn,amsfonts,amsmath}
\usepackage{graphicx}
\usepackage{subfigure}
\usepackage{lipsum}
\usepackage{lscape}
\usepackage{multicol}
\usepackage{setspace}
\usepackage{placeins}
\usepackage{changepage}
\usepackage[foot]{amsaddr}
\allowdisplaybreaks
\theoremstyle{plain}
\newtheorem{thrm}{Theorem}[section]
\newtheorem{lmm}{Lemma}[section]

\newtheorem{prpstn}{Proposition}[section]
\newtheorem{crllr}{Corollary}[section]
\newtheorem{dfntn}{Definition}[section]
\newtheorem{rmrk}{Remark}[section]
\newtheorem{ass}{Assumption}[section]
\newtheorem{algm}{Algorithm}[section]

\def\E{{\mathbb E}}
\def\indic{{\rm {\large 1}\hspace{-2.3pt}{\large l}}}

\def\R{{\mathbb R}}

\def\N{{\mathbb N}}

\def\argmin{{\rm argmin}}

\def\bb1{1}

\def\bb{{\bf b}}

\usepackage{multicol,enumerate}
\usepackage{tabularx,booktabs}
\usepackage{changepage}

\begin{document}


\title{High-dimensional Instrumental Variables Regression and Confidence Sets}

%
%
%
\author[Gautier]{Eric Gautier$^{{(1)}}$}
\address{$ ^{(1)}$ Toulouse School of Economics, Toulouse Capitole University, 1 esplanade de l'universit\'e, 31000 Toulouse, France.}
\email{\href{mailto:eric.gautier@tse-fr.eu}
{eric.gautier@tse-fr.eu}}
\author[Rose]{Christiern Rose$^{{(2)}}$}
\address{$ ^{(2)}$ School of Economics, University of Queensland, St Lucia, Brisbane, Australia, 4072.}
\email{\href{mailto:christiern.rose@uq.edu.au}{christiern.rose@uq.edu.au}}

\thanks{We thank James Stock, Elie Tamer, Ulrich M\"uller, and the referees for their useful comments. We warmly thank Alexandre Tsybakov whose reflections on this paper have been fundamental. We acknowledge financial support from ERC POEMH and ANR-17-EURE-0010. This is the text of our most recent revision}

\begin{abstract}
This article considers inference in linear instrumental variables models with many regressors, all of which could be endogenous. We propose the STIV estimator. Identification robust confidence sets are derived by solving linear programs. We present results on rates of convergence, variable selection, confidence sets which adapt to the sparsity, and analyze confidence bands for vectors of linear functions using bias correction. We also provide solutions to some instruments being endogenous. 
The application is to the EASI demand system.
\end{abstract}

\maketitle
\section{Introduction}\label{s1}
The high-dimensional paradigm concerns models in which the number of regressors $d_X$ is large relative to the number of observations $n$ but there is an unknown small set of relevant regressors. This can happen for various reasons. Researchers increasingly have access to large datasets and theory is often silent on the correct regressors. The number of observations can be limited because data is costly to obtain, because there simply exist few units (e.g., countries), or because the researcher is interested in a stratified analysis. The usual fixed $d_X$ large $n$ asymptotic framework does not necessarily provide a good approximation when there is high-dimensionality. A challenging situation is when $d_X$ is much larger than $n$  ($d_X\gg n$). Comparing models for all subsets of regressors is impossible when ${d_X}$ is even moderately large. The main focus of the high-dimensional literature is therefore the analysis of computationally feasible methods. For high-dimensional regression, the Lasso \cite{Tib} involves an $\ell^1$-penalty. The Dantzig Selector \cite{CT} is a linear program (henceforth LP).\\ 
\indent We study the high-dimensional linear instrumental variables (henceforth IVs) model where all regressors can be endogenous but the parameter $\beta$ is {\it sparse}, meaning it has few nonzero entries, or {\em approximately sparse}, meaning it is well approximated by a sparse vector. Sparsity can arise naturally when $\beta$ has an economic interpretation. Examples include social effects models with unobserved networks, models with uncertain exclusion restrictions, and treatment models with group heterogeneity in the treatment effect and many groups. 
Approximate sparsity is more appropriate when a linear model is used to approximate a function and the regressors and IVs comprise functions (e.g., splines) of baseline regressors and IVs. This can arise due to linearization or the use of series. The latter is relevant if the structural model has a nonparametric component or such components arise by including controls to justify IV exogeneity.\\ 
\indent With many endogenous regressors the number of IVs, $d_Z$, can be large. We allow for $d_Z\gg n$ but of order less than $\exp(n)$. 
Strong IVs are often scarce, particularly when there are many endogenous regressors. 
For these reasons, we pay 
attention to finite sample validity and robustness to identification.  Indeed, if there are weak and/or many IVs, inference based on standard asymptotic approximations can fail 
even if ${d_X}$ is small. To achieve this, we use a $\ell^{\infty}$-norm statistic derived from the moment condition. This is close in spirit to identification robust test inversion in which the exogeneity is the null hypothesis and a confidence set is formed by parameters which are not rejected. The Anderson-Rubin test is an example. In practice, such tests are conducted over a grid, which is only feasible for small ${d_X}$. To allow for large $d_X$, we use 
convex relaxations (linear or conic) instead. 
Our approach does not estimate reduced form equations 
and imposes no structure on them (e.g., sparsity).\\ 
\indent 
We propose the {\em Self Tuning Instrumental Variables} (henceforth STIV) estimator and establish its error bounds, which are used to obtain confidence sets for a vector of functions of $\beta$ and rates of convergence. 
Some confidence sets 
are uniform over identifiable parameters 
and 
distributions of the data among classes which leave the dependence between the regressors and IVs unrestricted, implying robustness to identification. 
Under stronger assumptions, including on the joint distribution of the IVs and regressors, the STIV estimator can be a pilot estimator to perform variable selection, obtain confidence sets which adapt to the sparsity, 
and conduct joint inference on linear functions of $\beta$ based on a data-driven bias correction. 
We also propose solutions to the problem that, in this data rich environment, 
a few IVs can be endogenous. All of our methods are pivotal because they jointly estimate standard deviations of structural errors or moments, making the tuning parameters data-driven. \\
\indent The application is to 
Engel curves in the EASI system of \cite{LP}. We show that the first-order in prices approximation error can be large and propose a second order approximation leading to a linear system 
with thousands of endogenous regressors. 
To achieve this, our theory allows for empirically relevant specificities, including knowledge of the relevance of certain regressors (e.g., price and quadratic expenditure), parameter restrictions (e.g., symmetry of the Slutsky matrix), approximation error, and systems of equations.

\section{Preliminaries}\label{s2}
{\bf Notations.} 
To simplify the exposition we consider an i.i.d. sample of size $n$. The population model comprises an outcome $Y$, regressors $X\in\R^{d_X}$, and IVs $Z\in\R^{d_Z}$ of joint distribution $\mathbb{P}$. $\E$ is the expectation under $\mathbb{P}$. For a mean zero random variable $A$, $\sigma_{A}\triangleq\E[A^2]^{1/2}$. We denote stacked matrices in bold, e.g.,  
$\bold{X}\in\mathcal{M}_{n,d_X}$, where $\mathcal{M}_{{d},{d'}}$ is the set of ${d}\times {d'}$ matrices. For 
$d\in\mathbb{N}$ and a random vector $W\in\R^{d}$, $\E_n[W]$ is the sample mean, 
$\bold{D}_{\bold{W}}$ is the diagonal matrix with entries 
$\mathbb{E}_n\left[W_k^2\right]^{-1/2}$ for $k\in[d]\triangleq \{1,2,...,d\}$, and $D_W$ its population counterpart. For 
$b\in\R^{d_X}$, 
$U(b)\triangleq Y-X^{\top}b$,  
 $\mathbb{P}(b)$ is the distribution of $\left(X,Z,U(b)\right)$ implied by $\mathbb{P}$, and $\widehat{\sigma}(b)^2\triangleq\E_n[U(b)^2]$.
We write $\widehat{\Psi}\triangleq \bold{D}_{\bold{Z}}\E_n[ZX^{\top}]\bold{D}_{\bold{X}}$ and $\Psi\triangleq D_Z\E[ZX^{\top}]D_X$. 
The set $S_I\subseteq [{d_X}]$ collects the indices of the regressors which are also IVs 
and $S_Q\subseteq [{d_X}]$ of size $d_Q$ collects those of the regressors of questionable relevance. 
When we make inference on a vector of functions,  
its dimension is $d_{\varPhi}$. 
Some results are asymptotic in $n\to\infty$ in which case 
${d_Z}$, ${d_X}$, $d_{\varPhi}$ and $d_Q$ can increase with $n$ and triangular arrays are permitted.  
Inequality between vectors is entrywise.   $M_{k,\cdot}$ (resp. $M_{\cdot,k}$) is the $k^{\rm th}$ row (resp. column) of $M$. 
$\indic$ is the indicator function. 
For $S\subseteq [d]$, $|S|$ is its cardinality and $S^c$ its complement. For $\Delta\in\mathbb{R}^d$, $S(\Delta)\triangleq\{k\in[d]:\ \Delta_{k}\ne0\}$ and  
$\Delta_S\triangleq\left(\Delta_{k}\indic{\{k\in S\}}\right)_{k\in[d]}$. 
$|\Delta|_p$ is the $\ell^p$-norm of 
$\Delta$ or a vectorization if $\Delta$ is a matrix.\vspace{0.2cm}

\indent {\bf Baseline moments model.} The linear IV model is
\begin{align}
&\E\left[ZU(\beta)\right]=0,\label{einstr}\\
&\beta\in\mathcal{B},\ \mathbb{P}(\beta)\in\mathcal{P},\label{econstraints}
\end{align}
where $\mathcal{B}\subseteq \mathbb{R}^{d_X}$ accounts for restrictions on $\beta$
and $\mathcal{P}$ is a nonparametric class, e.g.,\\ 
{\bf Class 1:} 
\emph{$Z_lU(\beta)$ is 
symmetric for all $l\in[d_Z]$ and ${d_Z}<9\alpha/\left(4e^3\Phi\left(-\sqrt{n}\right)\right)$}\\ 
where $\alpha\in(0,1)$ is a confidence level and $\Phi$ the normal CDF. Other classes allow for non i.d. and dependent data, and asymmetry (see Section \ref{cs3}). Their basic versions do not restrict the joint distribution of $Z$ and $X$. 
All but Class 4 allow for conditional heteroscedasticity.  
The set $\mathcal{I}$ collects the vectors which satisfy \eqref{einstr}-\eqref{econstraints}. Our results are for all $\beta\in\mathcal{I}$, hence for the true $\beta^*$. 
\vspace{0.2cm}

\indent {\bf The $\ell^{\infty}$-norm statistic.} 
Our confidence sets and estimators use slack versions of \eqref{einstr} 
based on the statistic $\widehat{t}(b)\triangleq\left|\bold{D}_{\bold{Z}}\E_n[ZU(b)]\right|_\infty/\widehat{\sigma}(b)$ for  $b\in\R^{d_X}$. We use, for $\beta\in\mathcal{I}$, the event 
$\mathcal{G}\triangleq\left\{\widehat{t}(\beta)\le \widehat{r}\right\}$. Taking $\widehat{r}=\underline{r}_n\left|\bold{D}_{\bold{Z}}\bold{Z}^{\top}\right|_{\infty}$, the base choice in the main text, and 
$\underline{r}_n=-\Phi^{-1}\left(9\alpha/(4{d_Z}e^3)\right)/\sqrt{n}$ for Class 1,\footnote{$\underline{r}_n\le2\log\left(4{d_Z}e^3/(9\alpha)\right)/\sqrt{n},$ $\forall\alpha\in[0,1],d_Z\ge1$(because $\Phi^{-1}(a)\ge2\log(a)$ if $0<a\le\exp(-1/(4\pi))$).}
yields $\mathbb{P}(\mathcal{G})\ge1-\alpha$ for all $n$ and $(\beta,\mathbb{P})$ such that $\beta\in\mathcal{I}$. Such a simple bound is possible due to the division by $\widehat{\sigma}(b)$. The set $\{b\in\mathcal{B}:\ \widehat{t}(b)\le \widehat{r}\}$ is a confidence set but it is infeasible because it is nonconvex and $b$ is high-dimensional and (approximately) sparse. 
Class 4 determines $\widehat{r}$ by bootstrap under conditional homoscedasticity.\vspace{0.2cm}

\indent {\bf Sparsity certificate.} 
A \emph{sparsity certificate} is a bound $s\in[d_Q]$ on the sparsity and $\mathcal{I}_s\triangleq\mathcal{I}\cap\left\{b\in\R^{d_X}:\ \left|S(b)\cap S_Q\right|\le s\right\}$ is the set of {\em $s$-sparse identifiable parameters}. For asymptotic results and triangular arrays, $s$ can depend on $n$. $\mathcal{I}_s$ can be a singleton when $s+d_X-d_Q<d_Z<d_X$ and sparsity implies 
exogenous regressors have a zero coefficient (i.e., they are excluded). This 
occurs when 
exclusion restrictions are uncertain (see \cite{KCFGI,KZCS}). 
When $d_Q=d_X$, 
 $\mathcal{I}_s$ is a singleton if there is a solution for only one of the 
${d_X \choose s}$ overdetermined systems based on  \eqref{einstr}-\eqref{econstraints} 
 and it is unique. Another condition (see \cite{CT}) 
is that all matrices formed from $2s$ columns of $\E[ZX^\top]$ have full rank.\vspace{0.2cm}

\indent {\bf Example SE.} The outcome $Y_j$ of individual $j\in[m]$ depends on peer outcomes. When the peers are unknown and there are endogenous peer effects, a linear model is $Y_j=\sum_{k\ne j}^m\rho_{j,k} Y_k+Z_j^\top\pi_j+U_j$, where $\rho_{j,k}$ is the effect of $k$ on $j$ and $Z_j$ are low-dimensional exogenous characteristics. 
If we set $\beta=(\rho_{j,1},...,\rho_{j,j-1},\rho_{j,j+1},...,\rho_{j,m},\pi^\top_j)^\top$ and all peers are unknown then $S_Q=[m-1]$, $S_I=S_Q^c$ and $P_j=S(\beta)\cap S_Q$ is the set of $j$'s peers.  If the network is sparse (e.g., due to costly link formation) then $|P_j|\ll d_X$. A sparsity certificate is an upper bound on the number of peers. 
The IVs are $Z=(Z_1^\top,Z_2^\top,...,Z_m^\top)^\top$, so $d_Z\geq d_X$.  If there are also exogenous peer effects then $Z_j$ is replaced by $Z$ in the structural equation and $d_Z<d_X$. 
\vspace{0.2cm}

\indent {\bf Example NP.} The model is $Y=f(\widetilde{X})+U$ with a nonparametric $f$ and an IV $\widetilde{Z}$ such that $\mathbb{E}[U|\widetilde{Z}]=0$. Assuming no approximation error, which we cover in Section \ref{approx}, \eqref{einstr} holds with $X_k=g_k^X(\widetilde{X})$ and $Z_l=g_l^Z(\widetilde{Z})$ for approximating functions $(g_k^X)_{k\in[d_X]}$ and $(g_l^Z)_{l\in[d_Z]}$. 
\vspace{0.2cm}

\indent {\bf Roadmap.} The paper is organized to progressively strengthen the assumptions. Here, we summarize our methods with the simplifications $d_Q=d_X$, and, for all $\beta\in\mathcal{I}$, $U(\beta)|Z$ is normally distributed with mean 0 and known variance $\sigma^2$. 
The simplifications permit to replace $\widehat\sigma(b)$ by $\sigma$ and $\widehat r$ by $r_n=-\Phi^{-1}(\alpha/(2d_Z))/\sqrt{n}$ of order $\log(d_Z)/\sqrt{n}$.\\
\indent A starting point is to find $b\in \R^{d_X}$ by minimizing $|\bold{D}_{\bold{X}}^{-1}b|_1$ subject to $\widehat{t}(b)\leq r_n$. A solution 
$\widehat\beta$ is obtained by solving a LP and is called the nonpivotal STIV estimator. 
STIV does not require a known $\sigma$ and is introduced in Section \ref{sSTIVintro}. To analyze the estimation error and construct confidence sets, we introduce sensitivity characteristics in Section \ref{s4}. To explain their role, we now take a $\beta\in\mathcal{I}$. 
Since $\widehat\beta$ is a minimizer, on the event $\mathcal{G}$ we can use $|\bold{D}_{\bold{X}}^{-1}\widehat\beta|_1\leq |\bold{D}_{\bold{X}}^{-1}\beta|_1$, $\widehat{t}(\widehat\beta)\leq r_n$ and $\widehat{t}(\beta)\leq r_n$. Letting $\widehat\Delta\triangleq\bold{D}_{\bold{X}}^{-1}(\widehat\beta-\beta)$, the first inequality implies 
$|\bold{D}_{\bold{X}}^{-1}\widehat\beta_{S(\beta)}|_1+|\bold{D}_{\bold{X}}^{-1}\widehat\beta_{S(\beta)^c}|_1\leq |\bold{D}_{\bold{X}}^{-1}\beta_{S(\beta)}|_1$, hence $|\bold{D}_{\bold{X}}^{-1}(\widehat\beta-\beta)_{S(\beta)^c}|_1\leq |\bold{D}_{\bold{X}}^{-1}\beta_{S(\beta)}|_1-|\bold{D}_{\bold{X}}^{-1}\widehat\beta_{S(\beta)}|_1$, and by the triangle inequality 
$|\widehat\Delta_{S(\beta)^c}|_1\leq |\widehat\Delta_{S(\beta)}|_1$.  The last two imply $|\widehat\Psi\widehat\Delta|_\infty\le \sigma(|\widehat{t}(\widehat{\beta})|_\infty+|\widehat{t}(\beta)|_\infty)\leq 2\sigma r_n$. For $k\in[d_X]$, we introduce a sensitivity 
\begin{equation}\widehat\kappa_{\ell_{k},S(\beta)}=\min_{\Delta\in\R^{d_X}:|\Delta_k|=1,|\Delta_{S(\beta)^c}|_1\leq |\Delta_{S(\beta)}|_1}|\widehat\Psi\Delta|_\infty,\label{kkk}\end{equation}
which gives $|\widehat\Delta_k| \leq |\widehat\Psi\widehat\Delta|_\infty/\widehat\kappa_{\ell_{k},S(\beta)}\leq 2\sigma r_n/\widehat\kappa_{\ell_{k},S(\beta)}$. The first inequality holds if $\widehat{\Delta}_k=0$ and otherwise follows by $|\widehat\Delta_{S(\beta)^c}|_1\leq |\widehat\Delta_{S(\beta)}|_1$, homogeneity, and
$$|\widehat{\Psi}\widehat{\Delta}|_{\infty}/|\widehat{\Delta}_k|
\ge\min_{\Delta:\ \Delta\ne0,
|\Delta_{S(\beta)^c}|_1\leq |\Delta_{S(\beta)}|_1}|\widehat{\Psi}
\Delta|_{\infty}/|\Delta_k|.$$ 
Omitting the constraint $|\Delta_{S(\beta)^c}|_1\leq |\Delta_{S(\beta)}|_1$ from \eqref{kkk} leads to a smaller $\widehat\kappa_{\ell_{k},S(\beta)}$, hence wider confidence sets below.
 We do not know $S(\beta)$ but if we know $|S(\beta)|\le s$, we replace $ \widehat\kappa_{\ell_{k},S(\beta)}$ by a lower bound $\widehat\kappa_{\ell_{k}}(s)$ (see Section \ref{sCS}). Adding  $|\widehat\Delta_{S(\beta)}|_1$ to $|\widehat\Delta_{S(\beta)^c}|_1\leq |\widehat\Delta_{S(\beta)}|_1$ yields $|\widehat\Delta|_1\leq 2|\widehat\Delta_{S(\beta)}|_1$ and using $|\widehat\Delta_{S(\beta)}|_1\leq s|\widehat\Delta_{S(\beta)}|_\infty\leq s|\widehat\Delta|_\infty$, we obtain 
 \begin{align}\label{LBCWS}
 \widehat\kappa_{\ell_{k},S(\beta)}&\geq
 \min_{j\in[d_X]}
 \min_{\Delta,\mu\in\R^{d_X}:-\mu\leq\Delta\leq\mu\leq\Delta_j,
 \mu_k=1,\mu^\top1\leq 2s\mu_j}|\widehat\Psi\Delta|_\infty
 \triangleq\widehat\kappa_{\ell_{k}}(s).
 \end{align}
This yields the bounds $\widehat\beta_k\pm2\E_n[X_k^2]^{-1/2}\sigma r_n/\widehat\kappa_{\ell_{k}}(s)$. $\widehat\kappa_{\ell_{k}}(s)$ is obtained by solving $d_X$ LPs so this is a computationally feasible confidence set. 
The coverage guarantee is uniform over $(\beta,\mathbb{P})$ such that $\beta\in\mathcal{I}_s$ for $\mathcal{P}$ from Class 1. The set is robust to arbitrarily weak IVs because $\mathcal{P}$ does not restrict the dependence between $X$ and $Z$. Also,  
 $\mathcal{I}_s$ need not be a singleton.  
In Section \ref{sec:main} 
we obtain rates of convergence for $\widehat\beta$ (possibly to a set) by replacing sensitivities with 
population analogues. 
A simple condition to analyze the rates is that, for every relevant regressor $X_k$  (i.e., $k\in S(\beta)$), there is a linear combination $\lambda_k^{\top}Z$ of IVs, where $|\lambda_k|_1\le 1$, which has low correlation (if $\E[Z]=0$) with the other regressors relative to $X_k$. This yields the rate
\begin{equation}\label{ratesE1}
\left|D_X^{-1}(\widehat\beta-\beta)\right|_q \lesssim r_n|S(\beta)|^{1/q}\left(\max_{k\in S(\beta)}\lambda_k^{\top}\Psi_{\cdot,k}\right)^{1/q}\left(\max_{k\in [d_X]}\lambda_k^{\top}\Psi_{\cdot,k}\right)^{1-1/q}
\end{equation}
for all $q\in[1,\infty]$. Assuming the nonzero entries of $\beta$ are sufficiently large relative to the $\ell^{\infty}$-rate, we obtain $S(\beta)\subseteq S(\widehat\beta)$, and if they are larger still, $S(\beta)=S(\widehat\beta^{\widehat{\omega}})$ for a thresholded estimator $\widehat\beta^{\widehat{\omega}}$. We then build confidence sets based on 
$\widehat\kappa_{\ell_{k}}(\widehat S)$ for such an estimator $\widehat{S}$.\\ 
\indent 
Our confidence sets need few assumptions, can be robust to identification, and are useful for inference on a function of the entire parameter $\beta$. However they may be conservative when stronger assumptions can be maintained and the object of interest has dimension much smaller than $d_X$. In Section \ref{CI} we present confidence bands for a vector of linear functions $\varPhi\beta$ based on a bias correction of $\widehat\beta$. A special case is a confidence interval.  
These are obtained by applying a variant of STIV 
to estimate 
$\Lambda$ satisfying $\Lambda\mathbb{E}[ZX^{\top}]=\varPhi$, and then combining $\widehat\Lambda$ with $\widehat\beta$. \\
\indent In Section \ref{sec:endiv}, we present a method to detect endogenous IVs. The basic idea is to use a variant of STIV to estimate the correlation of the IVs with the residuals (i.e., $\bold{U}(\widehat\beta)$) from a first-stage STIV estimator which uses only the IVs known to be exogenous.\\
\indent In Section \ref{s8}, we conduct simulations and apply STIV to build confidence bands around the EASI Engel curves. This is similar in spirit to Example NP, but the approximation error arises due to linearization rather than use of approximating functions. To make full use of economic theory, we require $S_Q\subset[d_X]$, 
a set $\mathcal{B}$ 
based on theory, and a minor modification of STIV to permit approximation errors and systems of structural equations. Proofs are in the appendix.\vspace{0.2cm}

\indent {\bf References.} High-dimensional estimation and inference has become an active field. To name a few; 
\cite{BCCH} uses Lasso type methods to estimate the optimal IV and make inference on a low-dimensional structural equation, \cite{FL} consider a nonconvex approach to IV estimation,  
\cite{Caner,CZ,CL} consider GMM with large dimensions but do not handle the high-dimensional regime. 
Inference for subvectors in high-dimension is an active topic related to Section \ref{CI} (see \cite{BCH,Jd_L,vdG,ZZ}, but also \cite{GLT,BCHN,Simoni,Kock} in the case of IVs. \cite{handbook} reviews results based on the nonpivotal STIV and others. Our results are applied to social effects models with unknown networks in \cite{Rose,GR} and \cite {BGMPR}.

\section{Self-Tuning IV Estimator and Confidence Sets}\label{sSTIVintro}
\begin{dfntn}
For $c>0$, a 
STIV estimator is any solution
$(\widehat{\beta},\widehat{\sigma})$ of
\begin{equation}\label{IVS}
\min_{b\in \widehat{\mathcal{I}}(\widehat{r},\sigma),\sigma\ge0} \left(\left|\bold{D}_{\bold{X}}^{-1}b_{S_Q}\right|_1+c\sigma \right),
\end{equation}
where, for $r,\sigma>0$, 
\begin{equation}\label{IVC}
\widehat{\mathcal{I}}(r,\sigma)\triangleq\left\{b\in\mathcal{B}, \left| \bold{D}_{\bold{Z}}\E_n[ZU(b)]\right|_{\infty}\le r\sigma,\widehat{\sigma}(b)\le \sigma\right\}.
\end{equation}
\end{dfntn} 
The $\ell^1$-norm 
is a convex relaxation of $|S(b)\cap S_Q|$. The term $c\sigma$ favors small $\sigma$, hence increasing $c$ tightens the set $\widehat{\mathcal{I}}(r,\sigma)$. 
$\bold{D}_{\bold{X}}^{-1}$ and $\bold{D}_{\bold{Z}}$ guarantee invariance to scale of the regressors and IVs. If $\mathcal{B}$ comprises linear (in)equality restrictions, a STIV  estimator is computed by solving a convex (second-order) conic program, similarly to the Square-root Lasso of \cite{BCW}. Linearity of \eqref{einstr} in $\beta$ is key to obtain such a simple program. 
$\widehat{\sigma}(\widehat{\beta})$ and $\widehat{\sigma}$ are estimators of the standard deviation of the structural error which need not be known.  
Taking $Z=X$ in the nonpivotal STIV estimator gives the Dantzig Selector. 

\indent Minimizing $\mathcal{O}(b)\triangleq\max\left(\widehat{\sigma}(b),\left|\bold{D}_{\bold{Z}}\E_n[ZU(b)]\right|_{\infty}/r\right)$ trades-off least-squares and exogeneity of the IVs, which is desirable in the presence of weak IVs (see \cite{AS}). STIV implements this in high-dimension because
\begin{equation}\label{esq}
\widehat{\beta}\in\argmin_{b\in\mathcal{B}}\left(\frac{1}{c}\left|\bold{D}_{\bold{X}}^{-1}b_{S_Q}\right|_1+\mathcal{O}(b)\right),\ \widehat{\sigma}=\mathcal{O}(\widehat{\beta}).
\end{equation}
If $\mathcal{O}$ were a differentiable and strictly convex function of $\bold{W}b$ 
and the entries of $\bold{W}$ drawn from a continuous distribution, minimizers of \eqref{esq} 
would be unique and 
one could obtain regularization paths (see \cite{Tib2}) for ad hoc determination of the penalty level.  
Our analysis is valid for all minimizers and determination of the penalty level is not an issue because STIV is pivotal.
Non uniqueness also occurs for LIML, which minimizes the Anderson-Rubin statistic. 

\subsection{Sensitivity Characteristics}\label{s4}
If $Z=X$, the minimal eigenvalue of $\E_n[XX^{\top}]$ can be used to obtain error bounds for quantities such as the mean squared error. It is the minimum of $b^{\top}\E_n[XX^{\top}]b/|b|_2^2$ over $b\in\R^{d_X}$, and is equal to zero if $d_X>n$. Under sparsity, 
$\R^{d_X}$ can be replaced by a subset in the case of the Dantzig Selector and Lasso. 
This is typically expressed via the 
restricted eigenvalue 
condition of \cite{BRT}. 
The sensitivity characteristics 
introduced in this paper are core elements to analyze 
STIV and provide sharper results for the Dantzig Selector and Lasso (see Section \ref{s:LBS}). They are 
related to the action of $\widehat{\Psi}$ on 
a 
subset $\widehat K_S$ for $S\subseteq[d_X]$.\\ 
\indent As in the Roadmap, we 
bound $|\widehat\Psi\widehat{\Delta}|_\infty$ on $\mathcal{G}$. 
To bound 
$\ell(\widehat{\Delta})$ for a loss $\ell$, 
we use 
a \emph{sensitivity} 
\begin{equation}\label{sensdef}\widehat\kappa_{\ell,S}\triangleq\min_{\Delta\in \widehat K_{S}:\
\ell(\Delta)=1}|\widehat{\Psi}\Delta |_{\infty}.\end{equation}
When $S=[d_X]$ we use the shorthand $\widehat\kappa_{\ell}$. 
We require that 
$\ell\in\mathcal{L}$, where $\mathcal{L}$ are the continuous functions from $\mathbb{R}^{d_X}$ to $[0,\infty)$ which are homogeneous of degree 1. An important loss is $\ell^{p}_{S_0}(\Delta)\triangleq|\Delta_{S_0}|_p$ for $p\in[1,\infty]$ and $S_0\subseteq [d_X]$. For $S_0=[d_X]$ and $S_0=\{k\}$ for $k\in[d_X]$, we use the shorthand notations   
$\ell^{p}$ and $\ell_k$. 
The sensitivities for these losses can be related to one another as expressed in Proposition \ref{p4}. Due to the $\ell^\infty$-norm in \eqref{sensdef}, 
additional IVs can only increase $|\widehat{\Psi}\Delta|_{\infty}$.  
Their cost is mild because it appears only through the $\log(d_Z)$ factor in $\underline{r}_n$.  

The cone $\{\Delta\in\R^{d_X}:\ |\Delta_{S(\beta)^c}|_1\leq |\Delta_{S(\beta)}|_1\}$ for the nonpivotal STIV is modified to be
\begin{equation}\label{eq:cone}
\widehat K_{S}\triangleq\left\{\Delta\in\R^{d_X}:\ \Delta_{S^c\cap S(\widehat{\beta})^c}=0,\
\left|\Delta_{S^c\cap S_Q}\right|_1\le \left|\Delta_{S\cap S_Q}\right|_1+c \widehat{g}(\Delta)\right\}.
\end{equation}
$\widehat{g}(\Delta)\triangleq \min(\widehat{r},1) |\Delta_{S_I}|_1+ \left|\Delta_{S_I^c}\right|_1$ is used because, by convexity and since the regressors of index in $S_I$ are used as IVs, $\widehat{\sigma}(\beta)-\widehat\sigma(\widehat\beta)\le\widehat{g}(\widehat\Delta)$. $\widehat{g}(\Delta)=|\widehat\Delta|_1$ when all regressors are endogenous. 
Similarly to the Roadmap, for every $\beta\in\mathcal{I}$, 
on the event $\mathcal{G}$, we have $\widehat\Delta\in \widehat K_{S(\beta)}$.
The error bounds for STIV in Proposition \ref{t1} are decreasing in the sensitivities, hence it is important that $\widehat K_{S}$ be small so that the sensitivities can be bounded away from zero. The researcher's knowledge components $\mathcal{B}$, $S_I$, and $S_Q$ serve this purpose. When, e.g., $\mathcal{B}$ comprises linear equalities $\{b:Mb=m\}$, we add 
$M\bold{D}_{\bold{X}}\Delta=0$ 
to $\widehat K_{S}$. Because $\widehat{r}<1$ is typical,  
accounting for $S_I$ yields a smaller set. 
If we omit $\Delta_{S^c\cap S(\widehat{\beta})^c}=0$, take $d_Q=d_X$, $\mathcal{B}=\R^{d_X}$, and replace $\widehat{g}(\Delta)$ by $|\Delta|_1$  
we obtain the simple cone of dominant coordinates
$\left\{\Delta\in\mathbb{R}^{d_X}:  (1-c)|\Delta_{S^c}|_1\le (1+c)|\Delta_S|_1\right\}$,  
due to which $\Delta$ has most of its $\ell^1$-norm concentrated on the indices in $S$. 
This cone is $\R^{d_X}$ if $c\ge 1$. Using the smaller $\widehat{K}_S$ is empirically relevant because 
in practice 
we find that STIV performs better for $c>1$. The constraint $\Delta_{S^c\cap S(\widehat{\beta})^c}=0$ can be removed 
to obtain rates of convergence, but is useful 
to construct confidence sets which are as small as possible. 
\\ 
\indent If $\beta$ is not sparse, $\widehat K_{S(\beta)}$ can be large (e.g., $\R^{d_X}$ when $S(\beta)=[d_X]$), so   
the sensitivities, denoted by $\widehat{\overline{\kappa}}$ 
instead of $\widehat\kappa$, are defined by 
replacing $\widehat{K}_{S}$ by 
$$\widehat{\overline{K}}_{S}\triangleq\left\{\Delta\in\mathbb{R}^{d_X}: \left|\Delta_{S^c\cap S_Q}\right|_1
\le 2\left(\left|\Delta_{S\cap S_Q}\right|_1+c\widehat{g}(\Delta)\right)+|\Delta_{S_Q^c}|_1\right\}.$$
Due to the additional terms on the right-hand side, 
$\widehat{\overline{K}}_{S}$ is larger than $\widehat{K}_{S}$. However, in our analysis these sensitivities need not be computed at $S=S(\beta)$. The slackness allows  
$\widehat\Delta\in \widehat{\overline{K}}_{S}$ on $\mathcal{G}$ provided that $\left|\bold{D}_{\bold{X}}^{-1}\beta_{S^c\cap S_Q}\right|_1$ is sufficiently small. The form of the additional terms is related to the factor 6  
in the second inequality in Proposition \ref{t1}. 

Our results involve the weakly increasing function 
$$\gamma(x)\triangleq1/\max(1-x,0),$$ where  
by convention $1/0= \infty$. It is close to 1 for small $x$ and is $\infty$ when $x\ge1$. We also let $\overline{\sigma}\triangleq(\widehat\sigma+\widehat{\sigma}(\widehat{\beta}))/2$ and $\widehat{h}(\Delta)\triangleq \min(|\Delta_{S_Q}|_1,(3|\Delta_{S\cap S_Q}|_1+c\widehat{r}|\Delta_{S_I}|_1+c|\Delta_{S_I^c}|_1+|\Delta_{S_Q^c}|_1)/2)$. 
\begin{prpstn}\label{t1} 
For all $(\beta,\mathbb{P})$ such that $\beta\in\mathcal{I}$, any STIV estimator, 
$\ell\in\mathcal{L}$, $q\in[1,\infty]$, $S_0,S\subseteq[{d_X}]$, and $c>0$, we have, on $\mathcal{G}$, 
$$
\ell\left(\bold{D}_{\bold{X}}^{-1}\left(\widehat\beta-\beta\right)\right) \le 
\frac{2\widehat{r}}{\widehat\kappa_{\ell,S(\beta)}}\min\left(
\overline{\sigma}\gamma\left(\frac{\widehat{r}}{\widehat\kappa_{\widehat{g},S(\beta)}}\right),
\widehat{\sigma}(\beta)\gamma\left(\frac{\widehat{r}}{c\widehat\kappa_{\ell^1_{S(\beta)\cap S_Q},S(\beta)}}\right)\right),
$$
\vspace{-.5cm}
\begin{align}
\left|\bold{D}_{\bold{X}}^{-1} \left(\widehat\beta-\beta\right)_{S_0}\right|_q
\le \max\left(
\frac{2\widehat{r}}{\widehat{\overline{\kappa}}_{\ell^q_{S_0},S}}
\min\left(\overline{\sigma}\gamma\left(\frac{\widehat{r}}{\widehat{\overline{\kappa}}_{\widehat{g},S}}\right),\widehat{\sigma}(\beta)\gamma\left(\frac{\widehat{r}}{c\widehat{\overline{\kappa}}_{\widehat{h},S}}\right)\right),
6 \left|\bold{D}_{\bold{X}}^{-1}\beta_{S^c\cap S_Q}\right|_1\right).\notag
\end{align}
\end{prpstn}
The first term in the minimum in the first inequality 
is used for confidence sets, 
and the second 
for rates of convergence. 
The second inequality is used when $\mathcal{I}$ contains nonsparse vectors, and is the basis of the sparsity oracle inequality in Theorem \ref{t1b} \eqref{t1biii}. 
 To obtain a confidence set one needs to circumvent the dependence of the sensitivities on the unknown $S(\beta)$ in a computationally feasible way. 
This is the focus of Section \ref{sCS}. For rates of convergence one requires population analogues of the upper bounds. 
This is the focus of Section \ref{sec:main}. 

\subsection{Computable Bounds on the Sensitivities and Confidence Sets}\label{sCS}
Confidence sets can be obtained by using lower bounds on the sensitivities. 
To obtain \eqref{LBCWS} in the Roadmap,   
we use a sparsity certificate $s$. Alternatively, one replaces $S(\beta)$ by $\widehat{S}$ such that $\widehat{S}\supseteq S(\beta)$ with probability converging to 1. We explain how STIV can be used to obtain such $\widehat{S}$ in Section \ref{sec:main}. 
We now present 
our base result, which provides bounds through LPs. 
\begin{prpstn}\label{prop:LBsensitivities}
For all $S\subseteq\widehat{S}\subseteq[{d_X}]$, $\ell\in\mathcal{L}$, $|S\cap S_Q|\le s$, and $c>0$, 
\begin{equation*}
\widehat\kappa_{\ell,S} \ge \max\left(\widehat\kappa_{\ell}\left(\widehat{S}\right),\widehat\kappa_{\ell}(s)\right),\quad \widehat{\overline{\kappa}}_{\ell,S} \ge \max\left(\widehat{\overline{\kappa}}_{\ell}\left(\widehat{S}\right),\widehat{\overline{\kappa}}_{\ell}(s)\right),
\end{equation*}
where the quantities in the bounds and losses $\ell$ are in Table \ref{TabConst}.
\end{prpstn}
Using a sparsity certificate $s$ and 
\begin{align}\label{CC}
\widehat C(s)\triangleq\left\{b\in\mathcal{B}: \forall \ell\in\mathcal{L},\forall c>0,\ \ell\left(\bold{D}_{\bold{X}}^{-1}\left(\widehat\beta-b\right)\right)\le 
\frac{2\widehat{r}\overline{\sigma} \gamma\left(\widehat{r}/\widehat\kappa_{\widehat{g}}(s)\right)}{\widehat\kappa_{\ell}(s)}\right\},
\end{align}
the confidence set for $\varPhi\beta$, where $\varPhi\in\mathcal{M}_{d_{\varPhi},d_X}$, denoted $\widehat{C}_{\varPhi}(s)\triangleq\{\varPhi b: b\in \widehat C(s) \}$
verifies 
\begin{align}\label{cover}
\min_{s\in[d_Q]}\inf_{(\beta,\mathbb{P}):\beta\in\mathcal{I}_s}
\mathbb{P}\left(\varPhi\beta\in\widehat{C}_{\varPhi}(s)\right)\ge 1-\alpha. 
\end{align}
The set is robust to identification if, as for Class 1, $\mathcal{P}$ does not restrict the dependence between $X$ and $Z$. Though we do not make it explicit, the bound in \eqref{CC}  depends on $c$. Increasing $c$ decreases $\overline{\sigma}$ (by increasing the penalty on $\sigma$ in the STIV  objective function) but increases 
$\gamma\left(\widehat{r}/\widehat\kappa_{\widehat{g}}(s)\right)/\widehat\kappa_{\ell}(s)$
(by enlarging $\widehat K_{S}$). 
The set \eqref{CC} can be made computationally feasible with correct coverage by replacing $\forall \ell\in\mathcal{L},\forall c>0$ with a finite intersection. 
If $\varPhi=I$ and  $\Xi$ is a grid for $c$, we can define the confidence set $[\widehat{\underline{C}}_{k}(s),\widehat{\overline{C}}_{k}(s)]$ for $k\in[d_X]$, with
\begin{align}\label{kbnds}
\hspace{-.7cm}\widehat{\underline{C}}_{k}(s)=\max_{c\in\Xi}\left(\widehat\beta_k-\frac{2\widehat{r}\overline{\sigma} \gamma\left(\widehat{r}/\widehat\kappa_{\widehat{g}}(s)\right)}{\widehat\kappa_{\ell_{k}}(s)\E_n[X_k^2]^{1/2}}\right),&\quad \widehat{\overline{C}}_{k}(s)=\min_{c\in\Xi}\left(\widehat\beta_k+\frac{2\widehat{r}\overline{\sigma} \gamma\left(\widehat{r}/\widehat\kappa_{\widehat{g}}(s)\right)}{\widehat\kappa_{\ell_{k}}(s)\E_n[X_k^2]^{1/2}}\right).
\end{align}
\begin{table}
{\small
  \caption{Lower bounds on sensitivities}
\label{TabConst}  
\begin{tabular}{rlrl}
\hline
$\widehat\kappa_{\ell^\infty_{S_0}}(\widehat{S})\triangleq$ &\hspace{-.5cm} $\displaystyle\min_{j\in S_0}
\min_{\substack{(\Delta,\mu,\nu)\in \widehat{B}(\widehat{S})\\ \Delta_j=1,\mu_{S_0}\le1}}
\nu$&$\widehat\kappa_{\ell^\infty_{S_0}}(s)\triangleq$ &\hspace{-.5cm} $\displaystyle\min_{j\in S_0} \min_{\substack{(\Delta,\mu,\nu)\in \widehat{B}(j)\\ \Delta_j=1,\mu_{S_0}\le1}}\nu$\\
$\widehat\kappa_{\ell_{k}}(\widehat S)\triangleq$ &\hspace{-.5cm} $\displaystyle\min_{\substack{j\in[{d_X}]\\
\eta=\pm1
}}\min_{
\substack{(\Delta,\mu,\nu)\in \widehat{B}(\widehat{S})\\ \mu_k= 1,\Delta_k=\eta,\mu\le \Delta_j}}\nu$&$\widehat\kappa_{\ell_{k}}(s)\triangleq$ &\hspace{-.5cm} $\displaystyle\min_{\substack{j\in[{d_X}]\\
\eta=\pm1}}\min_{
\substack{(\Delta,\mu,\nu)\in \widehat{B}(j)\\ \mu_k= 1,\Delta_k=\eta,\mu\le \Delta_j
}}\nu
$\\
$\widehat\kappa_{\widehat{g}}(\widehat{S})\triangleq$ &\hspace{-.5cm} $\displaystyle\min_{j\in[{d_X}]}\hspace{-.2cm}
\min_{
\substack{(\Delta,\mu,\nu)\in \widehat{B}(\widehat{S})\\ 
\sum_{k\in S_I}\widehat{r}\mu_k+\sum_{k\in S_{I}^c}\mu_k=1,\mu\le \Delta_j}}\hspace{-1.2cm}\nu$&
$\widehat\kappa_{\widehat{g}}(s)\triangleq$ &\hspace{-.5cm} $\displaystyle\min_{j\in[{d_X}]}\hspace{-.2cm}
\min_{\substack{(\Delta,\mu,\nu)\in \widehat{B}(j)\\ 
\sum_{k\in S_I}\widehat{r}\mu_k+\sum_{k\in S_{I}^c}\mu_k=1,\mu\le \Delta_j}}\hspace{-1.2cm}\nu$\\
\hline
$\widehat{B}(\widehat{S})\triangleq$&\multicolumn{3}{l}{\hspace{-.3cm}$\left\{
\begin{array}{l}
-\mu\le\Delta\le \mu,\mu_{\widehat{S}^c\cap S(\widehat{\beta})^c}=0,-\nu1\le\widehat{\Psi}\Delta\le \nu1\\
(1-c\widehat{r})\sum_{k\in S_I}\mu_k+(1-c)\sum_{k\in S_{I}^c}\mu_k\le 2
\sum_{k\in \widehat{S}\cap S_Q}\mu_k+\sum_{k\in S_Q^c}\mu_k
\end{array}
\right\}$}\\
$\widehat{B}(j)\triangleq$&\multicolumn{3}{l}{\hspace{-.3cm}$\displaystyle\left\{
\begin{array}{l}
-\mu\le\Delta\le \mu,-\nu1\le\widehat{\Psi}\Delta\le \nu1\\ 
(1-c\widehat{r})\sum_{k\in S_I}\mu_k+(1-c)\sum_{k\in S_{I}^c}\mu_k\le 2
s\mu_j+\sum_{k\in S_Q^c}\mu_k
\end{array}
\right\}$}\\
\hline
\multicolumn{4}{p{120mm}}{\scriptsize \textbf{Notes:} $\widehat{\kappa}_{\ell^{\infty}_{S_0},S}$ is also bounded by \eqref{p44i} in Proposition \ref{p4}. Bounds for $\widehat\kappa_{\ell^1,S}$ (resp. $\widehat\kappa_{\ell_{\varphi},S}$) 
replace $\sum_{k\in S_I}\widehat{r}\mu_k+\sum_{k\in S_{I}^c}\mu_k=1$ (resp. $\mu_k=1,\Delta_k=\eta$) by $1^{\top}\mu= 1$ (resp. $\varphi^{\top}\bold{D}_{\bold{X}}^{-1}\Delta= \eta$) in the bounds for $\widehat\kappa_{\widehat{g},S}$ (resp. $\widehat\kappa_{\ell_l,S}$). Section \ref{s:LBS} gives sharper but more computationally demanding bounds. }
\end{tabular}

}
\vspace{-5mm}
\end{table}
We can replace $\widehat\kappa_{\ell_{k}}(s)$ by $\widehat\kappa_{\ell^\infty}(s)$ to obtain a larger but less computationally demanding set (i.e., with less LPs).  
If $\varPhi\ne I$, we use the loss $\ell_{\varphi}\triangleq|\varphi^{\top}\bold{D}_{\bold{X}}^{-1}\Delta|$, where $\varphi^{\top}=\varPhi_{f,\cdot}$ for $f\in[d_{\varPhi}]$.\\ 
\indent The above confidence sets are nonempty hyperrectangles, and are infinite if $\widehat\kappa_{\widehat{g}}(s)\le \widehat{r}$. This is unavoidable for sets which are robust to weak IVs (see \cite{Dufour}). 
 Section \ref{sec:sim} provides a rule of thumb to determine a single value of $c$. Even if $c$ is determined from the data, the set has coverage at least $1-\alpha$ 
 due to \eqref{CC}. 
Because the researcher may be unsure about an appropriate value of $s$, the minimum over $s$ in 
\eqref{cover} allows to construct nested sets over different values. This can be used to assess the information content of progressively stronger sparsity assumptions.

{\bf Example SE continued.} A sparsity certificate (upper bound on the number of peers) yields a confidence set for the peer effects. By \eqref{cover}, the estimator $\widehat{P}_j=\{k\in S_Q:0\notin [\widehat{\underline{C}}_{k}(s),\widehat{\overline{C}}_{k}(s)]\}$ of the peers
satisfies 
$\min_{s\in[d_Q]}\inf_{(\beta,\mathbb{P}):\beta\in\mathcal{I}_s}
\mathbb{P}(\widehat{P}_j\subseteq P_j)\geq1-\alpha$. 
A subset is unavoidable because the peer effects can be arbitrarily close to zero, an issue to which we return in Section \ref{sec:main}. A confidence interval for the average peer effect uses $\varPhi\beta=(\sum_{k\neq j}\rho_{j,k})/(m-1)$.
\subsection{Deterministic Error Bounds, Model Selection, and Refined Confidence Sets}\label{sec:main}
We give deterministic counterparts of the 
bounds in Proposition \ref{t1} based on an 
event $\mathcal{G}_{A1}$, on which $(\tau_n)_{n\in\N}\in(0,1)^\N$ controls the deviation of the sample from the population, where 
\begin{equation}
\log(\max(d_Z,d_X,d_{\varPhi}))/(n\tau_n^2)\to0.\label{etau}
\end{equation} 
On $\mathcal{G}_{A1}$, $\widehat{r}$ can be replaced by a deterministic upper bound $r_n$ (see \eqref{GA1}, which also defines $\mathcal{G}_{A1}$)  and $\widehat{\sigma}(\beta)$ and 
the sensitivities by 
population analogues $\sigma_{U(\beta)}$ and 
$\kappa$ and $\overline{\kappa}$, obtained by replacing $\widehat{\Psi}$ by $\Psi$ and $\widehat{K}_S$, $\widehat{\overline{K}}_{S}$ by 
$K_S$, 
$\overline{K}_{S}$ (see Lemma \ref{thrm:DLBsensitivities}). 
We restrict $\mathcal{P}$ using Assumption \ref{ass:Nemirovski}, which places mild restrictions on second moments and the tails of the IVs so that 
$\mathbb{P}(\mathcal{G}_{A1})\ge 1-\alpha^{A1}_n$, where $\alpha^{A1}_n\to 0$ is defined in \eqref{alphapsi}. 
Asymptotic statements allow $c$ to depend on $n$.  
In the discussion of orders below Theorem \ref{t1b} and in Section \ref{CI}, $Z$ and $X$ are assumed uniformly bounded, so $\alpha_n^{A1}\to 0$ under 
\eqref{etau} and, 
for the choice of $\widehat{r}$ using classes 1-3, 
$r_n$ has same order as $\underline{r}_n$ (i.e., $\log(d_Z/\alpha)/\sqrt{n}$). Section \ref{cs3} presents the general case. 

\subsubsection{Deterministic Error Bounds and Rates of Convergence}\label{s54}
The deterministic bounds below hold without additional assumptions. We 
leave the dependence between $Z$ and $X$ unrestricted, allow for 
partial identification, and  $c\gg1$ when $S_I\ne\emptyset$, which works well in practice. For the bounds to be orders in probability, 
we replace the confidence level $\alpha$ used to set $\widehat{r}$ by $(\alpha_n)_{n\in\N}\in(0,1)^{\N}$ converging to 0, so 
$$\mathbb{P}(\mathcal{G}\cap\mathcal{G}_{A1})\geq1-\alpha_n^{S},\ \text{where}\ \alpha_n^{S}\triangleq\alpha_n+\alpha^{A1}_n\to0.$$ 
For a function $\omega$ from $\R^{d_X}$ to $\R^{d_X}$ and given $\mathbb{P}$, we set 
$$\mathcal{I}(\omega)\triangleq\{b\in \mathcal{I}: \forall k\in S(b),1_{n}\E[X_k^2]^{1/2}|b_k|>\omega_k(b)\}\ \text{and}\ 1_{n}\triangleq\sqrt{(1-\tau_n)/(1+\tau_n)}.$$

\begin{thrm}\label{t1b} Let $c>0$ and $\mathcal{P}$ be such that Assumption \ref{ass:Nemirovski} holds. 
\begin{enumerate}[\textup{(}i\textup{)}]  
\item\label{t1bi} 
For all  $(\beta,\mathbb{P})$ such that $\beta\in\mathcal{I}$ and any STIV estimator, 
we have, on $\mathcal{G}\cap\mathcal{G}_{A1}$, for all $\ell\in\mathcal{L}$,
\begin{align}
&\ell\left( D_X^{-1}\left(\widehat\beta-\beta\right)\right)\le 
\frac{2r_n\sigma_{U(\beta)}}{1_{n}\kappa_{\ell,S(\beta)}}\Gamma_{\kappa}(S(\beta)),\notag
\end{align}
where for $S\subseteq[d_X]$, 
$\Gamma_{\kappa}(S)\triangleq(1+\tau_n)\gamma
(\tau_n/\kappa_{\ell^1,S}+r_n(1+\tau_n)/(c \kappa_{\ell^1_{S\cap S_Q}, S}))$.
\item
\label{t1bii} For all  $(\beta,\mathbb{P})$ such that $\beta\in\mathcal{I}(\underline{\omega})$, where $\underline{\omega}_k:\ b\in\R^{d_X}\to 2r_n\sigma_{U(b)}\Gamma_{\kappa}(S(b))/\kappa_{\ell_{k},S(b)}$,  and any STIV estimator, 
we have, on $\mathcal{G}\cap\mathcal{G}_{A1}$, 
$S(\beta) \subseteq S(\widehat\beta)$.
\item
\label{t1biii}
For all  $(\beta,\mathbb{P})$ such that $\beta\in\mathcal{I}$ and any STIV estimator, 
we have, on $\mathcal{G}\cap\mathcal{G}_{A1}$, 
for all $q\in[1,\infty]$ and $S_0\subseteq[{d_X}]$,
$$
\left| D_X^{-1} \left(\widehat\beta-\beta\right)_{S_0}\right|_q
\le \min_{S\subseteq[{d_X}]}\max\left(
\frac{2r_n\sigma_{U(\beta)}}{1_{n}\overline{\kappa}_{\ell^q_{S_0},S}}
\Gamma_{\overline{\kappa}}(S),\frac{6}{1_{n}}\left| D_X^{-1}\beta_{S^c\cap S_Q}\right|_1\right),$$
where $\Gamma_{\overline{\kappa}}$ 
(resp. $h$) replaces $\kappa_{\ell^1\hspace{-.1cm},S},\kappa_{\ell^1_{S\cap S_Q}, S}$ by $\overline{\kappa}_{\ell^1\hspace{-.1cm},S},\overline{\kappa}_{h,S}$ in $\Gamma_{\kappa}$ 
(resp. $\widehat{r}$ by $r_n$ in $\widehat{h}$).  
\end{enumerate}
\end{thrm} 
Theorem \ref{t1b} \eqref{t1bi}-\eqref{t1bii} provide a bound and a model selection result suited to the sparse case. 
Result \eqref{t1biii} gives an alternative bound 
better suited to approximate sparsity. 
The dependence of the bounds on the function $\gamma$ is unavoidable. It means that they can be infinite, and so hold with high probability regardless of the dependence between $Z$ and $X$. Bounds for $\ell^1$-loss are used in the next section, in which STIV is used as a pilot estimator. The loss $\ell^{\infty}_{S_0}$ for $S_0\subseteq [d_X]$ is 
used in \eqref{t1bii} with $S_0=\{k\}$. It can be used to obtain  uniform rates of convergence for the coefficients of index in $S_0$ (e.g., $S(\beta)$). 
The \emph{bona fide} loss for model selection is $\ell^{\infty}$ (see \cite{Lounici}).\\ 
\indent Result \eqref{t1bii} means that, for all $\beta\in\mathcal{I}(\underline\omega)$, STIV finds a superset of the regressors. The term $\underline{\omega}_k$ corresponds to the upper bound on $1_{n}\E[X_k^2]^{1/2}|\widehat\beta_k-\beta_k|$.  Using $\ell_k$ in its definition allows a larger $\mathcal{I}(\underline\omega)$ than using $\ell^{\infty}$. 
 Due to \eqref{t1bii}, the confidence set $\widehat{C}_{\varPhi}\triangleq\{\varPhi b: b\in \widehat C\}$, where 
\begin{align}\label{ESSB}
\widehat C\triangleq\left\{b: \forall \ell\in\mathcal{L},\ell\left(\bold{D}_{\bold{X}}^{-1}\left(\widehat\beta-b\right)\right)\le 
\frac{2\widehat{r}\overline{\sigma}\gamma\left(\widehat{r}/\widehat\kappa_{\widehat{g}}(S(\widehat \beta))\right)}{\widehat\kappa_{\ell}(S(\widehat \beta))}\right\},
\end{align}
is such that
\begin{equation}\label{ecoverb}\inf_{(\beta,\mathbb{P}):\ \beta\in\mathcal{I}(\underline\omega)}
\mathbb{P}\left(\varPhi \beta\in\widehat{C}_{\varPhi}\right)\ge 1-\alpha-\alpha^{A1}_n.  
\end{equation}
It is not robust to identification because 
 $\underline\omega$ depends on the population sensitivities (which depend on $\Psi$), hence on the joint distribution of $Z$ and $X$ (recalling $\Psi=D_Z\E[ZX^\top]D_X$). The condition $\beta\in\mathcal{I}(\underline\omega)$ in Theorem \ref{t1b} \eqref{t1bii} is a beta-min condition. It requires that the nonzero entries of $\beta$ be large enough, and is interpretable if $\beta$ is a structural parameter. It is not intended to be used if the regressors are used to approximate a function as in Example NP. 

{\bf Example SE continued.} The beta-min condition means that the peer effects are sufficiently large so as to be distinguishable from zero. It is reasonable because a typical parameterization when the social effect is via the mean (see \cite{bram}) is $\rho_{j,k}=\overline{\rho} \indic_{\{k\in P_j\}}/|P_j|$ where $\overline{\rho}$ is a scalar, so when the network is sparse the relevant effects are bounded away from zero. STIV finds a superset of the peers with asymptotic (uniform) probability at least $1-\alpha$.

 \indent If $\beta\in\mathcal{I}$ is sparse, the upper bound in Theorem \ref{t1b} \eqref{t1biii}, which holds for all $S\subseteq [d_X]$, also applies to $S=S(\beta)$, for which the second term in the maximum is zero. We are then left with a bound similar to the right-hand side of \eqref{t1bi}. 
When $q=1$ and $S_0=S_Q=[d_X]$, it is 6 times the error made when $\widehat{\beta}=\beta_S$ (estimating perfectly the components in $S$). In this sense the second term is an approximation error.
STIV performs a data-driven trade-off for nonsparse parameter vectors. 
Result \eqref{t1biii} implies that, for an optimal set $S_*\subseteq[d_X]$ (not necessarily $S(\beta)$), 
\begin{equation}\label{eq:UBASdecay}
\left| D_X^{-1} \left(\widehat\beta-\beta\right)\right|_1 \le 
\frac{2r_n\sigma_{U(\beta)}}{1_{n}\overline{\kappa}_{\ell^1,S_*}}\Gamma_{\overline{\kappa}}(S_*).
\end{equation}
This allows us to define formally approximately sparse parameter vectors as vectors which are sufficiently well approximated by a sparse vector so that the right-hand side of \eqref{eq:UBASdecay} is small.
\begin{rmrk}\label{rset}
Theorem \ref{t1b} applies if $\mathcal{I}$ is not a singleton, in which case, for a given $\widehat\beta$, one can take the infimum over $\beta\in\mathcal{I}$ on both sides of the inequality in \eqref{t1bi} and \eqref{t1biii}. The left-hand side can be viewed as the distance to a set, and the right-hand side defines the elements of $\mathcal{I}$ to which $\widehat\beta$ is closest.  
The discussion below uses such $\beta$. For a model which is not indexed by $n$, if there is $\beta\in\mathcal{I}$ such that, for a constant $C<\infty$, for $n$ large enough, $\sigma_{U(\beta)}\Gamma_{\kappa}(S(\beta))\le C\kappa_{\ell,S(\beta)}$, then $\widehat{\beta}$ converges to such $\beta\in\mathcal{I}$. It will become apparent from the lower bounds on $\kappa_{\ell,S(\beta)}$ that these are usually sparse vectors in $\mathcal{I}$. 
When $\beta$ are coefficients of a function on a collection of simple functions and $d_X$ increases with $n$, due to \eqref{t1biii}, $\widehat{\beta}$ can converge to the 
coefficients of a smooth function and the population sensitivities vary with $n$. 
It is typically the case in nonparametric IV that the coefficients decay rapidly to zero (see Example NP continued below).
\end{rmrk}

\begin{rmrk}
Deterministic bounds for $\widehat{\sigma}(\widehat{\beta})$ and $\widehat{\sigma}$ are given in Lemma \ref{sigmabounds}. These can be used to justify applying nonpivotal STIV in two-stages, or to obtain confidence bands such as those in Section \ref{CI} under conditional homoskedasticity.  
\end{rmrk}
We now discuss rates of convergence based on the bounds in Theorem \ref{t1b}, which 
depend on the population sensitivities. For ease of exposition we focus on explaining \eqref{t1bi}. Proposition \ref{p5} relates the population sensitivities to one another, so we start by considering the following alternative expression for $\kappa_{\ell^{\infty}_{S_0},S}$ 
for all $S_0,S\subseteq [d_X]$, 
\begin{equation}\label{epcc}
\kappa_{\ell^{\infty}_{S_0},S}
=\min_{k\in S_0}\min_{
\Delta\in K_{S}:
\Delta_k=1,|\Delta_{S_0}|_{\infty}\le1
}\max_{l\in[d_Z]}\left|\Psi_{l,k}-\sum_{k'\ne k}\Psi_{l,k'}\Delta_{k'}\right|,
\end{equation}
which has a natural interpretation
as a measure of the strength of the IVs for the regressors in $S_0$. If the 
IVs are centered, the second minimum in \eqref{epcc} is 
a maximum absolute normalized partial 
covariance between regressor $k$ and the IVs, where the partialling-out of the other regressors is restricted (i.e., $\Delta$ is constrained). Ignoring for the moment the constraints on $\Delta$, the second minimum in \eqref{epcc} is zero if 
$D_Z\E[Z(XD_X)_{k}]\in\R^{d_Z}$ lies in a vector space of dimension at most $\min(d_X-1, d_Z )$ of 
Lebesgue measure zero 
if $d_Z\geq d_X$. The vector space has a smaller maximum dimension $\min(n,d_X-1, d_Z )$ for $\widehat{\kappa}_{\ell^{\infty}_{S_0},S}$. This contrasts 
with the restricted and sparse eigen and singular values (see \cite{BRT,BCHN}) which can be zero even if $|S|< \min(n,d_X-1, d_Z )$ (else are always zero) and depend on $S$ only via its size. 

For simplicity of exposition we now use Condition IC  (see the appendix for the general case), under which we provide interpretable 
conditions to derive explicit rates. 
\vspace{.2cm}

\noindent{\bf Condition IC.} $d_Q=d_X$ 
and $c$ is a constant such that $c<1_{n}$.\vspace{.2cm}

\noindent Under Condition IC, by Proposition \ref{p5}, we have, for all $S\subseteq S_0\subseteq[d_X]$ and $q\in[1,\infty]$, 
\begin{align}
&u_{\kappa}|S|^{1-1/q}\kappa_{\ell^1,S}\ge \kappa_{\ell^q_S,S}\ge\kappa_{\ell^q_{S_0},S}\ge\kappa_{\ell^\infty_{S_0},S}/(\min(u_\kappa|S|,|S_0|))^{1/q},\label{kbndIC}\\
&K_{S}=\{\Delta\in\R^{d_X}: (1_{n}-c)|\Delta_{S^c}|_1\le (1+c) |\Delta_{S}|_1\}=
\{\Delta\in\R^{d_X}: |\Delta|_1\le u_{\kappa} |\Delta_{S}|_1\}\label{inclCone}, 
\end{align}
where $u_\kappa\triangleq (1+1_{n})/(1_{n}-c)$. 
Due to the form of $\Gamma_{\kappa}(S(\beta))$ and \eqref{kbndIC}, the upper bound in Theorem \ref{t1b} \eqref{t1bi} and \eqref{t1bii} is finite if 
\begin{equation}\label{cnontrivial}
\kappa_{\ell^1_{S(\beta)},S(\beta)}>\tau_n u_{\kappa}+r_n(1+\tau_n)/c\quad(\sim r_n/c).
\end{equation}
Sufficient conditions for consistency can be obtained based on the easier to interpret $\kappa_{\ell^{\infty}_{S(\beta)},S(\beta)}$. If $S\subseteq S_0$ (as when both are $S(\beta)$), 
we can further interpret the expression of  $\kappa_{\ell^{\infty}_{S_0},S}$ in \eqref{epcc} due to the constraints on $\Delta$. The constraints \eqref{inclCone}, $\Delta_k=1$, and $|\Delta_{S_0}|_\infty\leq 1$ imply that the subvector of $\Delta$ appearing in $\sum_{k'\ne k}\Psi_{l,k'}\Delta_{k'}$ has $\ell^{1}$-norm $m\le u_{\kappa}|S|-1$. Also, if $k\in S$, the subvector with indices in $S\setminus\{k\}$ has $\ell^1$-norm larger than $(m+1)/u_{\kappa}+1$ and $\ell^{\infty}$-norm smaller than 1. This restricts the set of vectors used to perform the partialling-out to have mass at most $u_{\kappa}|S|-1$, predominantly concentrated on $S$ (e.g., $S(\beta)$). 
Moreover, by Proposition \ref{p5} \eqref{p5i}
\begin{align}
&\forall S\subseteq S_0,\ 
\kappa_{\ell^{\infty}_{S_0},S}\ge \sup_{\eta\in(0,1)} \eta
\min_{k\in S_0}w_k(S,\eta),\label{slb2}\\
\hspace{-.7cm}\text{where }\quad&S_k(S,\eta)\triangleq\left\{
\lambda\in\R^{d_Z}:|\lambda|_1\le1,
(u_\kappa |S|-1)\max_{k'\ne k}|\lambda^{\top}\Psi_{\cdot,k'}|\le (1-\eta)\lambda^{\top}\Psi_{\cdot,k}
 \right\},\notag\\ 
 &w_k(S,\eta)\triangleq
\max_{\lambda\in S_k(S,\eta)}\lambda^{\top}\Psi_{\cdot,k}= \max_{\lambda\in S_k(S,\eta)}\E[(\lambda^{\top}D_ZZ)(XD_X)_{k}].\notag
 \end{align}

\noindent{\bf Assumption C$(S_0,q,\eta,\eta_0)$.} Condition IC holds, 
$\beta\in\mathcal{I}$, $\mathbb{P}\in\mathcal{P}$ satisfies Assumption \ref{ass:Nemirovski}, $\alpha=\alpha_n$, and 
\begin{enumerate}[\textup{(}i\textup{)}]  
\item\label{CSi} for all $k\in S_0$ if $q>1$ (resp. for all $k\in S(\beta)$ if $q=1$),  $S_k(S(\beta),\eta)\neq\emptyset$, 
\item\label{CSii} 
$c\eta_0\eta\min_{k\in S(\beta)}w_k(S(\beta),\eta)\ge r_n|S(\beta)|$ for $n$ large enough,
\item\label{CSiii} $\min_{k\in S(\beta)}w_k(S(\beta),\eta)/(r_n|S(\beta)|)\to \infty$  if $q=1$, else $\min_{k\in S_0}w_k(S(\beta),\eta)/r_n
\to \infty$. 
\end{enumerate}
\vspace{.2cm}

By Theorem \ref{t1b} \eqref{t1bi}, \eqref{kbndIC}, and H\"older's inequality, we obtain the following corollary.
\begin{crllr}
Let $S_0\subseteq [d_X]$, $q\in[1,\infty]$, and $(\eta,\eta_0)\in(0,1)^2$. For all $(\beta,\mathbb{P})$ satisfying Assumption C$(S_0,q,\eta,\eta_0)$, $S_0\supseteq S(\beta)$, 
and any STIV estimator, 
we have, 
for $n$ large enough, 
\begin{equation}
\hspace{-.2cm}\mathbb{P}\left(\left|D_X^{-1}(\widehat\beta-\beta)_{S_0}\right|_q \ge
\frac{
2r_n\sigma_{U(\beta)}|S(\beta)|^{1/q}
(1+\tau_n)\gamma(\eta_0)/(1_{n}\eta)}{
\min_{k\in S(\beta)}w_k(S(\beta),\eta)^{1/q}
\min_{k\in S_0}w_k(S(\beta),\eta)^{1-1/q}}\right)\le \alpha_n^{S}
.\label{lqcon}
\end{equation}
\end{crllr}
Taking $S_0=S(\beta)\cup\{k\}$ and $q=\infty$ yields rates for the $\ell_k$-loss and upper bounds on $\underline\omega_k$ in Theorem \ref{t1b} \eqref{t1bii} 
(see also \eqref{eellk}). For simplicity, the discussion now uses the word 
`correlation' as if the IVs and/or regressors were mean zero. Assumption C($S_0,q,\eta,\eta_0$) \eqref{CSi} means that for regressor $k$ there exists a nonempty set of 
linear combinations of the IVs of small enough relative absolute correlation with the other regressors. It is similar to the coherence condition for symmetric matrices of \cite{DET}, but more general because it allows for rectangular matrices and linear combinations of the IVs (i.e., $\lambda^{\top}\Psi_{\cdot,k}$ rather than $|\Psi_{l,k}|$ for $l\in[d_Z]$). 
The coherence condition is used to study $\ell^{\infty}$-norm convergence rates and model selection in \cite{Lounici}.  
 Item \eqref{CSii} 
is introduced to guarantee \eqref{cnontrivial} and \eqref{CSiii} for consistency.  
They require that for each regressor of index $k\in S(\beta)$ (but not for the other regressors), there is a linear combination of the IVs, which does not need to be known, of large enough absolute correlation with $X_k$. 
If $q=1$ then, by \eqref{CSiii}, \eqref{CSii} holds for all $\eta_0$ for $n$ large enough. Consistency can hold with $d_Z<d_X$.
\begin{rmrk}\label{stiv2sls}
Assume $d_X$ is fixed, 
we add to Assumption C$([d_X],1,\eta,\eta_0)$ that, for all $k\in S(\beta)$, $S_k(S(\beta),\eta)$ contains the vectors from the canonical basis of $\mathbb{R}^{d_Z}$ with a 1 at the  indices of the $d_R$ largest entries in absolute value of $\Psi_{\cdot,k}$, and $|\Psi_{\cdot,k}|_2=\psi$ where $\psi$ does not vary with $n$. STIV is consistent when $\log(d_Z)^2/(n\rho_n^2)\to0$, where $\rho_n=\min_{k\in S(\beta)}|\Psi_{\cdot,k}|_\infty$.  
Assume $\rho_n=\psi/\sqrt{d_R}$, so the $d_R$ IVs are equally relevant. 
If $d_R=d_Z/d_X$ then STIV is consistent if $d_Z\log(d_Z)^2/n\to0$. Like 2SLS, it may not be consistent if $d_Z/n$ converges to a nonzero constant. 
When $d_Z/n\to 0$ but $d_Z\log(d_Z)^2/n\not\to0$, 2SLS is consistent (see \cite{chaoswanson}) but STIV may not be. If $d_R=1$ then, for each relevant regressor, all but one of the IVs can be arbitrarily irrelevant and 
STIV is consistent when $d_Z \lesssim \exp(\sqrt{n}\epsilon_n)$ with $\epsilon_n\to0$. 
\end{rmrk}
Remarkably, for $q=1$, \eqref{lqcon} is not affected if all IVs are irrelevant for an irrelevant regressor. This is important to handle ill-posed inverse problems such as the following. 

{\bf Example NP continued.} Assume the baseline endogenous regressor and IV are related via $\widetilde{X}=\pi\widetilde{Z}+\sigma V$ and the approximating functions are $X_k=h_{k-1}(\widetilde{X}/\sqrt{\pi^2+\sigma^2})$ and $Z_k=h_{k-1}(\widetilde{Z})$ for $k\in[d_X]$, where $h_k$ is the $k^{\text{th}}$ Hermite polynomial. If, for simplicity, $(\widetilde{Z},V)$ follows a standard normal distribution, $\Psi$ is diagonal with $\Psi_{k,k}=(1+(\sigma/\pi)^2)^{(1-k)/2}$ 
(see Section \ref{mainonline}), so the $\ell^1$-rate depends on the exponentially small $\Psi_{k,k}$ for the largest $k\in S(\beta)$. 
Due to Theorem \ref{t1b} \eqref{t1biii}, a bound on the $\ell^1$-rate without sparsity is 
\begin{equation}\label{rNP}
\min_{S\subseteq[d_X]}\max\left(\frac{\log(d_Z)}{\sqrt{n}}\max_{k\in S}\left(\exp\left(\frac{(k-1)}{2}\left(\frac{\sigma}{\pi}\right)^2\right)\right),|\beta_{S^c}|_1\right).
\end{equation}

Assumption C$(S_0,q,\eta,\eta_0)$ \eqref{CSi} 
is in line with the common empirical practice of, for each endogenous regressor, finding an IV which is correlated more specifically with that regressor, and arises naturally in our application. To obtain adaptive nonparametric estimators in statistical inverse problems using series, it is common to use basis functions adapted to the operator so that $\Psi$ is (nearly) diagonal 
(see, e.g., 
\cite{HR} and \cite{GLP} in conjunction with wavelet/needlet thresholding and Galerkin approximation), and so Assumption C$(S_0,q,\eta,\eta_0)$ \eqref{CSi} is not restrictive.

For the sake of comparison, we present an assumption similar to that in \cite{BCHN}.\vspace{.2cm}

\noindent{\bf Assumption SV$(q,(\underline{\delta}_n)_{n\in\N},(\overline{\delta}_n)_{n\in\N},(l_n)_{n\in\N},\eta_0)$.} Condition IC holds, 
$\beta\in\mathcal{I}$, $\mathbb{P}\in\mathcal{P}$ satisfies Assumption \ref{ass:Nemirovski}, $\alpha=\alpha_n$, and 
\begin{enumerate}[\textup{(}i\textup{)}]  
\item\label{SVi} $4|S(\beta)|\overline{\delta}_n^2u_{\kappa}^2/\underline{\delta}_n^2$ is an integer smaller than $|S(\beta)| l_n$ and 
\begin{equation}\underline{\delta}_n\le 
\min_{\substack{K\subseteq[d_X]\\|K|\leq |S(\beta)| l_n}}\max_{\substack{L\subseteq[d_Z]\\|L|\leq |S(\beta)| l_n}}\hspace{-.3cm} 
\sigma_{\min}(\Psi_{L,K})\le 
\max_{\substack{K\subseteq [d_X]\\
|K|\leq |S(\beta)| l_n}}\max_{\substack{L\subseteq[d_Z]\\|L|\leq|S(\beta)| l_n}} \hspace{-.3cm}
\sigma_{\max}(\Psi_{L,K})\le \overline{\delta}_n,\label{SV}
\end{equation}
 where $\sigma_{\min}(M)$ and $\sigma_{\max}(M)$ are respectively the smallest and largest singular values of $M$ and $M_{L,K}$ is the submatrix obtained by extracting the rows in $L$ and columns in $K$, 
\item\label{SVii} 
$c\eta_0 \underline{\delta}_n^2/(4(1+u_{\kappa})u_{\kappa}^2\overline{\delta}_n)\ge r_n |S(\beta)|$,
\item\label{SViii} $\underline{\delta}_n^2/(4(1+u_{\kappa})u_{\kappa}^2\overline{\delta}_nr_n|S(\beta)|^{1/q})\to \infty$. 
\end{enumerate}
\vspace{.2cm}  

By adapting the proofs to apply to population sensitivities and all $c\in(0,1_n)$, we obtain
\begin{crllr}\label{crllrrate}
Let $q\in\{1,2\}$, $(\underline{\delta}_n)_{n\in\N}$, $(\overline{\delta}_n)_{n\in\N}$, $(l_n)_{n\in\N}$, and $\eta_0\in(0,1)$. For all $(\beta,\mathbb{P})$ satisfying Assumption SV$(q,(\underline{\delta}_n)_{n\in\N},(\overline{\delta}_n)_{n\in\N},(l_n)_{n\in\N},\eta_0)$, 
and any STIV estimator, 
we have
\begin{equation}
\hspace{-.2cm}\mathbb{P}\left(\left|D_X^{-1}(\widehat\beta-\beta)\right|_q \ge
8r_n\sigma_{U(\beta)}|S(\beta)|^{1/q}(1+\tau_n)\gamma(\eta_0)(1+u_{\kappa})u_{\kappa}^2
\overline{\delta}_n^2/(\underline{\delta}_n1_{n})\right)\le\alpha_n^{S}
.\label{svdl1}
\end{equation}
\end{crllr}
%

For $q=1$, $\underline{\delta}_n^2/\overline{\delta}_n$ plays the same role as $\min_{k\in S(\beta)}w_k(S(\beta),\eta)$ 
 under C$([d_X],1,\eta,\eta_0)$.  Item \eqref{SVii} 
guarantees \eqref{cnontrivial} and \eqref{SViii} gives consistency. Assumption SV$(1,(\underline{\delta}_n)_{n\in\N},(\overline{\delta}_n)_{n\in\N},(l_n)_{n\in\N},\eta_0)$ can be more appealing than C$([d_X],1,\eta,\eta_0)$ \eqref{CSi}. However, 
 $\underline{\delta}_n$ can be small (even 0) due to one irrelevant regressor. For example, suppose there is $k\in S(\beta)^c$ such that $\Psi_{\cdot,k}=(r_n,0,\hdots,0)^\top$. Taking $K=\{k\}$ in the first inequality of \eqref{SV}, $\underline\delta_n\leq r_n$, hence (unlike \eqref{lqcon}) the upper bound in \eqref{svdl1} does not converge to 0. The fundamental issue is that Assumption SV provides a rate based on the worst-case subset of regressors, regardless of their relevance. This is less costly for regression (i.e., $Z=X$) than for IV because exogenous regressors have higher correlation with the IVs than do endogenous regressors, and given many endogenous regressors it is likely that one is weakly correlated with the IVs. 
 
{\bf Example NP continued.} The bound on the rate in \eqref{svdl1} is  
$r_n|S(\beta)|^{1/q}(1+(\sigma/\pi)^2)^{(d_X-1)/2}$ while $d_X$ is replaced by $k=\max\{l:l\in S(\beta)\}$ under Assumption C$([d_X],1,\eta,\eta_0)$. 

Assumption SV does not apply to $\ell^{\infty}_{S_0}$-losses, hence cannot be used for model selection. We provide a more technical comparison with Assumption SV and the results of \cite{BCHN}, and a sharper condition in the same spirit in Section \ref{mainonline}. 

\subsubsection{Selection of Variables and Confidence Sets with Estimated Support}\label{sthreshold}
Theorem \ref{t1b} \eqref{t1bii} provides a superset of the relevant regressors. Under a stronger beta-min condition exact selection can be performed. For this purpose, we use a purely data-driven 
thresholded STIV estimator $\widehat{\beta}^{\widehat\omega}$ which uses a sparsity certificate. It is defined by
\begin{equation}
\widehat{\beta}_{k}^{\widehat\omega}\triangleq\widehat{\beta}_k\indic{\left\{\mathbb{E}_n[X_k^2]^{1/2}|\widehat{\beta}_k|>\widehat\omega_k(s)\right\}},\quad
\widehat\omega_k(s)\triangleq \frac{2\widehat{r}\overline{\sigma} \gamma\left(\widehat{r}/\widehat\kappa_{\widehat{g}}(s)\right)}{\widehat\kappa_{\ell_{k}}(s)}.\label{eq:thresh}
\end{equation}
for $k\in[{d_X}]$. 
The following theorem shows that this estimator
achieves sign consistency and hence, $S(\widehat{\beta}^{\widehat\omega})=S(\beta)$ for all $\beta\in\mathcal{I}_s\cap\mathcal{I}(2 \omega(s))$. It uses ${\rm sign}(b)\triangleq\left({\rm sign}(b_k)\right)_{k\in[d_X]}$, where ${\rm sign}(t)\triangleq \indic{\{t>0\}}-\indic{\{t<0\}}$, and makes use of the population counterparts, for all $k\in[d_X]$, 
\begin{align*}
\omega_k(s):\ b\in\R^{d_X}\to \frac{2r_n\sigma_{U(b)}
\gamma(r/\kappa_{g}(s))}{\kappa_{\ell_k}(s)}
\sqrt{1+\tau_n}
\left(1+\frac{2r_n\Gamma_\kappa(S(b))}{c\kappa_{\ell^1_{S(b)\cap S_Q},S(b)}}\right),
\end{align*}
\vspace{-.3cm}

\noindent where $g$ and $\kappa_{\ell}(s)$ 
are defined 
before \eqref{cone1} and in \eqref{thrm:DLBsensitivities1bomega}. 
Under Condition IC, $\kappa_{g}(s)=\kappa_{\ell^1}(s)$ (and $\widehat\kappa_{\widehat{g}}(s)=\widehat{\kappa}_{\ell^1}(s)$).  
\begin{thrm}\label{th:threshold}
Let $s\in[d_Q]$ and $\mathcal{P}$ be such that Assumption \ref{ass:Nemirovski} holds. For all $(\beta,\mathbb{P})$ such that $\beta\in\mathcal{I}_s\cap\mathcal{I}(2 \omega(s))$ and any STIV estimator,  
we have, 
 on $\mathcal{G}\cap\mathcal{G}_{A1}$,  
 ${\rm sign} (\widehat{\beta}^{\widehat\omega})={\rm sign}(\beta)$. 
\end{thrm}
{\bf Example SE continued.} By Theorem \ref{th:threshold}, the peers are exactly recovered with asymptotic probability at least $1-\alpha$ if the endogenous effects are sufficiently large. 

Theorem \ref{th:threshold} yields a confidence set by replacing $S(\widehat{\beta})$ by $S(\widehat\beta^{\widehat{\omega}})$ in \eqref{ESSB}, which satisfies \eqref{ecoverb} with $\mathcal{I}_s\cap\mathcal{I}(2 \omega(s))$ in place of $\mathcal{I}(\underline\omega)$.
The value of $s$ can be large (possibly ${d_X}$). 
The set's width 
matches the error bound in Proposition \ref{t1} with respect to $S(\beta)$, hence it adapts to the sparsity. To achieve this we remove a small set from $\mathcal{I}$ (vectors too close to $|S(\beta)|$-sparse vectors). 

\section{Confidence Bands using Bias Correction}\label{CI} 
Confidence sets for $\varPhi \beta$, where $\varPhi\in\mathcal{M}_{d_{\varPhi},d_X}$,  can be robust to identification and are  particularly useful when one is interested in a feature of the whole parameter vector such as the network in Example SE. But they can be conservative when stronger assumptions on the data generating process can be maintained and $d_{\varPhi}$ is small relative to $d_X$. The confidence bands below address this. 
Using $d_\varPhi=1$ yields a confidence interval (e.g., for the average peer effect in Example SE). Using $d_\varPhi>1$ one can build a confidence band for a structural function such as $f$ in Example NP or the Engel curves in Section \ref{sec:appl}. A first estimator is the plug-in $\varPhi\widehat{\beta}$. 
Another uses 
\begin{equation}\label{lInv}
\exists\Lambda\in\mathcal{M}_{d_{\varPhi},{d_Z}}:\ 
\Lambda\mathbb{E}[ZX^{\top}]=\varPhi.
\end{equation}
Indeed, by \eqref{einstr}, for all $\beta\in\mathcal{I}$ and $\Lambda$ which solves \eqref{lInv}, $\varPhi\beta=\Lambda\E[ZY]$. \eqref{lInv} is a system of equations of the same form as the following equation derived from \eqref{einstr}
\begin{equation}\label{par}
\exists \beta\in\R^{d_X}:\ \beta^{\top}\mathbb{E}[ZX^{\top}]^{\top}=\E[ZY]^{\top}.
\end{equation} 
A STIV estimator $\widehat{\Lambda}$ is obtained 
by solving \eqref{Mhat}. For simplicity, we assume \eqref{lInv} holds exactly but, as in Section \ref{approx}, one can handle an approximation error going to zero with $n$. 

\indent Using either plug-in strategy poses problems 
because STIV is 
"biased" towards zero and can converge slowly. 
To deal with this we combine the two to form the bias corrected estimator 
\begin{align}\label{BCE}
\widehat{\varPhi\beta}\triangleq \varPhi\widehat{\beta}+\frac1{n}\widehat{\Lambda}\bold{Z}^{\top}\bold{U}(\widehat{\beta})
\end{align}
and build a confidence band around $\widehat{\varPhi\beta}$. 
\begin{rmrk}\label{rmrkDR} 
\eqref{BCE} is close in spirit to the bias correction in \cite{Jd_L}. Another motivation for it 
is that $O:\ (b,L)\to \varPhi b+L (\E[Z^{\top}Y]-\E[Z^{\top}X]b)$ has zero partial derivatives 
at, respectively, identified $\Lambda$ and $\beta$ (due to \eqref{lInv} and \eqref{par}) and $O(\beta,\Lambda)=\varPhi\beta$. This 
is a type of double-robustness (see \cite{CCDDHNR}). In this paper $\beta$ appears in a structural equation and  
our analysis does not involve machine learning for regressions.
\end{rmrk}
\begin{dfntn} For $\lambda\in(0,1)$, a BC-STIV estimator is any solution $(\widehat{\Lambda},\widehat{\nu})$ of 
\begin{equation}\label{Mhat}
\min_{L\in\widehat{\mathcal{I}}_{\varPhi}(\underline{r}'_n,\nu),\nu>0} \left|L\bold{D}_{\bold{Z}}^{-1}\right|_{1}+\frac{\lambda\nu}{\widehat{\rho}^{ZX}},
\end{equation}
where $\widehat{\rho}^{ZX}\triangleq\widehat{\rho}^{ZX}_{[d_Z]}$ (see \eqref{rhoZX}) and for $r,\nu>0$,  
\begin{align*}
&\widehat{\mathcal{I}}_{\varPhi}(r,\nu)\triangleq\Big\{L\in\mathcal{M}_{d_{\varPhi},d_Z}:\ \left|\left(\varPhi-L\E_n[ZX^\top]\right)\bold{D}_{\bold{X}}\right|_{\infty}\le r\nu,\widehat{\Sigma}\left(L\right)\le\nu\Big\},\\
&\widehat{\Sigma}\left(L\right)\triangleq\max_{(f,k)\in [d_{\varPhi}]\times[{d_X}]}\widehat{\sigma}_{f,k}(L),\ \widehat{\sigma}_{f,k}(L)^2\triangleq
\E_n\left[
(\varPhi-LZX^\top)_{f,k}^2\right]\left(\bold{D}_{\bold{X}}\right)_{k,k}^2.
\end{align*}
\end{dfntn}
To choose $\underline{r}'_n$, one uses one of classes 1-3, replacing $\alpha$ by $\alpha_n$ 
and $d_Z$ by $d_{\varPhi}d_X$. 
If $\varPhi=I$ and $Z=X$, $\widehat{\Lambda}$ is an approximate inverse of 
$\mathbb{E}_n[XX^{\top}]$, which improves on the CLIME estimator of \cite{CLL} by estimating standard errors. BC-STIV 
differs from STIV in that it is for a system of $d_{\varPhi}$ (rather than 1) equations, each with $d_Z$ (rather than 1) second-order cones, making it more computationally intensive. We provide a computational solution and its analysis in Section \ref{Lambdaest}. The counterpart of $\mathcal{I}$ is $\mathcal{I}_{\varPhi}\triangleq\left\{\beta\in\mathcal{I},\Lambda:
\Lambda\E\left[ZX^\top\right]=\varPhi,\mathbb{P}\left(\beta,\Lambda\right)\in\mathcal{P}_{\varPhi}
\right\}$,  
where $\mathcal{P}_{\varPhi}$ is a 
class for the distribution of 
$(X,Z,U(\beta),\varPhi-\Lambda ZX^\top\hspace{-.1cm}, 
\Lambda ZU(\beta))$.
Asymptotically uniformly valid 
confidence bands are obtained as
\begin{align}
\hspace{-0.5cm}\widehat{C}_{\varPhi}\triangleq\left[\widehat{\varPhi\beta}-\widehat{q},\widehat{\varPhi\beta}+\widehat{q}\right],\quad \widehat{q}\triangleq\frac{
q_{G_{\varPhi}|\bold{F}(\widehat{\beta})\widehat{\Lambda}^{\top}}(1-\alpha)+3\zeta_n
}{\sqrt{n}}\bold{D}_{\bold{F}(\widehat{\beta})\widehat{\Lambda}^{\top}}^{-1}1,\label{bndOmega}
\end{align}
where $q_{G_{\varPhi}|\bold{F}(\widehat{\beta})\widehat{\Lambda}^{\top}}(1-\alpha)$ is the $1-\alpha$ quantile of $G_{\varPhi}=|\bold{D}_{\bold{F}(\widehat{\beta})\widehat{\Lambda}^{\top}}\widehat{\Lambda}\bold{F}(\widehat{\beta})^\top\bold{E}|_{\infty}/\sqrt{n}$ given $\bold{F}(\widehat{\beta})\widehat{\Lambda}^\top$ (obtained 
by simulation),  $F(b)\triangleq ZU(b)$, 
$\bold{E}\in\R^n$ is a standard Gaussian vector independent of $\bold{F}(\widehat{\beta})\widehat{\Lambda}^{\top}$, and $(\zeta_n)_{n\in\N}$ is a positive sequence. 
\begin{thrm}\label{tCI} Let $\mathcal{P}_{\varPhi}$ and $(\zeta_n)_{n\in\N}$ be such that Assumption 
\ref{assCI1a} holds. Then, for $\alpha^{A2}_n\to0$ defined in \eqref{alphaPhi}, for all $(\beta,\Lambda,\mathbb{P})$ such that $(\beta,\Lambda)\in\mathcal{I}_{\varPhi}$, 
$$\forall n\in\N,\ \mathbb{P}\left(\varPhi\beta\in \widehat{C}_{\varPhi}\right)\ge 1-\alpha-\alpha^{A2}_n.$$
\end{thrm}
The sequence 
$(\zeta_n)_{n\in\N}$ restricts $\mathcal{P}_{\varPhi}$ and $\mathcal{I}_{\varPhi}$ on which uniformity over distributions and parameters holds. We denote by $(v^\beta_n)_{n\in\mathbb{N}}$ (resp. $(v^{\Lambda,\beta}_n)_{n\in\mathbb{N}}$) the deterministic upper bound on $|D_X^{-1}(\widehat\beta-\beta)|_1$ (resp. on $|D_{\Lambda F(\beta)}(\widehat\Lambda-\Lambda)D_Z^{-1}|_{\infty,\infty}\sigma_{U(\beta)}$, see  Proposition \ref{tCLambda}),  
where $|\cdot|_{\infty,\infty}$ is the maximum row-wise $\ell^1$-norm. These bounds hold on an event of probability at least $1-\alpha^{S}_n-\alpha^{BC}_n$, where $\alpha^{BC}_n$ is defined in Assumption \ref{assCI0}, converging to one, and $v^{\Lambda,\beta}_n$ and $v^{\beta}_n$ can depend on $(\beta,\Lambda)\in\mathcal{I}_{\varPhi}$. 
Proposition \ref{tCLambda} is the analogue 
of Theorem \ref{t1b}. It provides deterministic upper bounds for useful losses and characterizes the limit of $\widehat{\Lambda}$ when there are multiple solutions to \eqref{lInv}, as discussed for $\widehat\beta$ in Remark \ref{rset}. 

For the discussion, we now take $\zeta_n\to0$ (e.g., $\log(n)^{-1}$), $(\beta,\Lambda)\in\mathcal{I}_{\varPhi}$ to which $(\widehat{\beta},\widehat{\Lambda})$ converges, 
and assume $X$ and $Z$ are uniformly bounded. Assumption \ref{assCI1a} holds if
\begin{align}
&\max(v^{\Lambda,\beta}_n,\tau_n)\max(\log(d_Z/\alpha_n)\log(n/\alpha_n), \log(d_\varPhi/\alpha_n)
=o(\zeta_n),\notag\\
&v^\beta_n|D_{\Lambda F(\beta)}|_\infty\max(\log(d_Xd_\varPhi),\left|\Lambda D_Z^{-1}
\right|_{\infty,\infty}\log(d_\varPhi/\alpha_n))=o(\zeta_n)\label{evb}.
\end{align}
The requirement on $v^{\Lambda,\beta}_n$ is mild. It can be logarithmic if $\max(d_Z,d_\varPhi)$ is of polynomial order in $n$. When $d_\varPhi=1$ (i.e., a confidence interval) and we use $|D_{\Lambda F(\beta)}(\widehat\Lambda-\Lambda)D_Z^{-1}|_{\infty,\infty}\sigma_{U(\beta)}\le
|(\widehat\Lambda-\Lambda)D_Z^{-1}|_{\infty,\infty}/(\sigma_{\min}(D_Z\mathbb{E}[ZZ^{\top}U(\beta)^2/\sigma_{U(\beta)}^2]D_Z)\min_{f}|\Lambda_{f,\cdot}D_Z^{-1}|_2))$, we can obtain similar upper bounds on $v_n^{\Lambda,\beta}$ as on 
$v_n^\beta$ under Condition IC in Section \ref{sec:main} because the set $K'_S$ used to analyze $\widehat\Lambda$ is equal to the set $K_S$ in \eqref{inclCone}, replacing $c$ by $\lambda$ and $d_X$ by $d_Z$. The key difference 
is that $Z$ and $X$ switch roles 
and $\Psi$ is replaced by $\Psi^\top$ in assumptions C and SV for $q=1$. 
So $S(\Lambda)$ plays the role of $S(\beta)$
under sparsity, and otherwise bounds in the spirit of \eqref{rNP} can be derived under approximate sparsity.   
The requirement on $v^\beta_n$ is also mild if the rows of $\Lambda$ are (approximately) sparse. For example, if $d_\varPhi=1$, 
$v^\beta_n|D_{\Lambda F(\beta)}|_\infty\left|\Lambda D_Z^{-1}
\right|_{\infty,\infty}\log(d_\varPhi/\alpha_n)$ 
is of order at most $(v_n^\beta/\sigma_{\min}(D_Z\mathbb{E}[ZZ^{\top}U(\beta)^2]D_Z))\sqrt{|S(\Lambda)|}\log(d_\varPhi/\alpha_n)$. Hence, if $|S(\Lambda)|$ grows slowly the rate of estimation of $\beta$ could also be logarithmic.  
 We provide alternative confidence bands under conditional homoscedasticity and their analysis in Section \ref{conhomban}. 
 
\section{Endogenous IVs}\label{sec:endiv}
With many endogenous regressors and IVs, the exogeneity of some IVs could fail. We now consider a high-dimensional framework for the problem of  
IV exogeneity (see, e.g., 
\cite{Sargan}). 
Introducing $\theta\in\R^{d_Z}$ to account for the possible failure of \eqref{einstr},  
we replace \eqref{einstr}-\eqref{econstraints} by
\begin{align}
&\E[ZU(\beta)-\theta]=0,\label{einstrNV}\\
&\left(\beta,\theta\right)\in\mathcal{B}\times\Theta,\ \mathbb{P}\left(\beta,\theta\right)\in\mathcal{P}_{\not\perp},\label{econstraintsNV}
\end{align}
 where $\theta_l\ne0$ means that $Z_l$ is endogenous, $\mathbb{P}\left(b,t\right)$ is the distribution of $\left(X,Z,ZU(b)-t\right)$ implied by $\mathbb{P}$ and $\Theta\subseteq\R^{d_Z}$ encodes restrictions on $\theta$. 
For example, the sign of the correlation of a regressor and the structural error could be known. 
Another restriction is $\theta_{ S_{\perp}}=0$ for $S_{\perp}\subseteq[d_Z]$ of cardinality $d_{\perp}$ which indexes the IVs known to be exogenous. 
The counterpart of $\mathcal{I}_s$, denoted by $\mathcal{I}_{s,\widetilde{s}}$, collects the vectors $(b,t)\in\mathcal{B}\times\Theta$ which satisfy \eqref{einstrNV}-\eqref{econstraintsNV}, $|S(b)\cap S_Q|\leq s$, and $|S(t)|\le \widetilde{s}\}$, where $\widetilde{s}\in[d_Z-d_{\perp}]$ is a sparsity certificate for the possibly endogenous IVs. 

To detect endogenous IVs, we use a variant of STIV to estimate $\theta$ by replacing $\bold{U}(\beta)$ by the residuals $\bold{U}(\widehat\beta)$ from a pilot 
STIV estimator which uses only the IVs in $S_\perp$ and 
$\widehat{r}^{\perp}$ (in place of $\widehat{r}$) based on $\mathcal{G}_{\perp}$, which differs from $\mathcal{G}$ by using only $d_\perp$ moments and $\alpha^{\perp}$ (in place of $\alpha$). Based on STIV, one then computes $\widehat\delta$ and $\widehat\delta^{\Sigma}$ such that for all $(\beta,\mathbb{P})$ such that $\beta\in\mathcal{I}$, on $\mathcal{G}_{\perp}$ 
\begin{align}\label{UB}
&\left|\left(\widehat{\Psi} \bold{D}_{\bold{X}}^{-1}\left(\widehat{\beta}-\beta\right)\right)_{ S_{\perp}^c}\right|_{\infty}
\le\widehat{\delta},\quad
\widehat{\rho}_{S_\perp^c}^{ZX}\left|\bold{D}_{\bold{X}}^{-1}\left(\widehat{\beta}-\beta\right)\right|_{1}
\le\widehat{\delta}^{\Sigma},\\
\label{rhoZX}
&\widehat{\rho}_{S}^{ZX}\triangleq \max_{l\in S,k\in [d_X]}\left(\bold{D}_{\bold{Z}}\right)_{l,l}\left(\bold{D}_{\bold{X}}\right)_{k,k}\mathbb{E}_n\left[Z_{l}^2X_k^2\right]^{1/2}.
\end{align}
Though the analysis does not require the sparsity certificate approach, if $\beta\in\mathcal{I}_s$ we can use 
\begin{equation}\label{bhatnvstiv}
\widehat{\delta}=2\widehat{r}^{\perp}\overline{\sigma}\gamma\left(\widehat{r}^{\perp}/\widehat\kappa_{\widehat{g}}(s)\right)/\widehat\kappa^{\Psi}(s)\quad \text{and}\quad
\ \widehat{\delta}^{\Sigma}=2\widehat{r}^{\perp}\widehat{\rho}_{S_\perp^c}^{ZX}\overline{\sigma}\gamma\left(\widehat{r}^{\perp}/\widehat\kappa_{\widehat{g}}(s)\right)/\widehat\kappa_{\ell^1}(s),\end{equation}
where the lower bound $\widehat\kappa^\Psi(s)$ on $\widehat\kappa_{S}^{\Psi}$, the sensitivity for the loss in the first inequality of \eqref{UB},\footnote{The sensitivities of the pilot STIV replace $|\widehat{\Psi}\Delta|_\infty$ (resp. $\widehat{r}$) by $|(\widehat{\Psi}\Delta)_{S_\perp}|_\infty$ (resp. $\widehat{r}^\perp$) in \eqref{sensdef} (resp. in $\widehat{g}$ in \eqref{eq:cone}).} is obtained by linear programming. 
Unlike the Hansen-Sargan test, we can use a pilot STIV estimator when $s< d_\perp<d_Z$ and with an approximately sparse reduced form. 
 \begin{dfntn}\label{dfntn51} For $\widetilde{c}>0$, a NV-STIV estimator is any solution $(\widehat{\theta}, \widehat{\widetilde{\sigma}})$ of 
\begin{equation}\label{def:STIV_est_nonvalid}
\min_{t\in
\widehat{\mathcal{I}}_{\not\perp}(\underline{r}^{\not\perp}_n,\widetilde{\sigma}),\widetilde{\sigma}\ge0}\left(\left|\bold{D}_{\bold{Z}}t_{ S_{\perp}^c}\right|_1+\widetilde{c}\widetilde{\sigma}\right),
\end{equation}
\vspace{-.4cm}

\noindent where, for $r,\sigma>0$, 
\vspace{-.6cm}

\begin{align*}
&\widehat{\mathcal{I}}_{\not\perp}(r,\sigma)\triangleq\left\{t\in\Theta:
\left|\bold{D}_{\bold{Z}}\left(\frac1n\bold{Z}^{\top}\bold{U}(\widehat\beta)-t\right)_{S_{\perp}^c}\right|_{\infty}\hspace{-0.3cm}\le r\sigma+ \widehat{\delta},\ \widehat{\Sigma}_{\not\perp}\left(\widehat{\beta},t\right) \le
\sigma + \widehat{\delta}^{\Sigma}\right\},\\
&\widehat{\Sigma}_{\not\perp}\left(b,t\right)\triangleq\max_{l\in  S_{\perp}^c}\ \widehat{\sigma}_l(b,t),\quad \widehat{\sigma}_l(b,t)^2\triangleq (\bold{D}_{\bold{Z}})_{l,l}^2\mathbb{E}_n[(Z_{l}U(b)-t_l)^2].
\end{align*}
\vspace{-.4cm}

\end{dfntn}
\noindent To set $\underline{r}^{\not\perp}_n$ to control $\mathbb{P}(\underline{\mathcal{G}}_{\not\perp})$ (see \eqref{G2}), 
we use one of classes 1-3, replacing $d_Z$ by $d_Z-d_\perp$ and $\alpha$ by $\alpha^{\not\perp}$. $\mathcal{P}_{\not\perp}$ in \eqref{econstraintsNV} combines the classes for $\widehat{r}^{\perp}$ and $\underline{r}^{\not\perp}_n$. For deterministic bounds, 
it is further restricted using a minor modification of Assumption \ref{ass:Nemirovski} and we modify the event $\mathcal{G}_{A1}$ of probability $1-\alpha_n^{A1}$ accordingly (see Section \ref{nvstiv}). We still refer to them as Assumption \ref{ass:Nemirovski} and $\mathcal{G}_{A1}$ in Theorem \ref{th:nonvalid}. 
Also, for a function $\omega$ from $\R^{d_X}\times\R^{d_Z}$ to $\R^{d_Z}$ and given $\mathbb{P}$, we set $\mathcal{I}_{s,\widetilde{s}}(\omega)\triangleq\{(b,t)\in \mathcal{I}_{s,\widetilde{s}}: \forall l\in S(t),|t_l|>\omega_l(b,t)((1+\tau_n)\E[Z_l^2])^{1/2}\}$. The definitions of $\widetilde{\underline{\omega}}$ 
and $\widetilde{\omega}$ 
used below 
are in Section \ref{nvstiv} and $\overline{\widetilde{\sigma}}\triangleq(\widehat{\widetilde{\sigma}} +\widehat{\Sigma}_{\not\perp}(\widehat\beta,\widehat\theta))/2$.

\begin{thrm}\label{th:nonvalid}
For all $s\in[d_Q]$ and $\widetilde{s}\in\left[d_Z-d_{\perp}\right]$, $(\beta,\theta,\mathbb{P})$ such that 
$\left(\beta,\theta\right)\in \mathcal{I}_{s,\widetilde{s}}$, $\widetilde{c}\in\left(0,1/\underline{r}^{\not\perp}_n\right)$ and $c>0$ and any NV-STIV estimator,
we have, on $\mathcal{G}_{\perp}\cap\underline{\mathcal{G}}_{\not\perp}$, 
\begin{align}
\hspace{-.6cm}\left|\bold{D}_{\bold{Z}}\left(\widehat{\theta}-\theta\right)_{S_{\perp}^c}\right|_\infty &\le\gamma\left(2(\underline{r}^{\not\perp}_n)^2\widetilde{s}\gamma(\widetilde{c}\underline{r}^{\not\perp}_n)\right)\left(\underline{r}^{\not\perp}_n\left(2\overline{\widetilde{\sigma}}+(1+\widetilde{c}\underline{r}^{\not\perp}_n)\gamma(\widetilde{c}\underline{r}^{\not\perp}_n)\widehat{\delta}^{\Sigma}\right)+2\widehat{\delta} \right)\triangleq \widehat{\widetilde{\omega}}\left(\widetilde{s}\right).\label{eq:th:nonvalidCSinfty}
\end{align}
For fixed $\widetilde{c}\in\left(0,1/\underline{r}^{\not\perp}_n\right)$, $c$, $s\in[d_Q]$ and $\widetilde{s}\in\left[d_Z-d_{\perp}\right]$, 
if we restrict $\mathcal{I}_{s,\widetilde{s}}$ so that Assumption \ref{ass:Nemirovski} holds,  
then, for all $(\beta,\theta,\mathbb{P})$, 
if $\left(\beta,\theta\right)\in\mathcal{I}_{s,\widetilde{s}}(\widetilde{\underline{\omega}}(s,\widetilde{s}))$, then $S(\theta)\subseteq S(\widehat{\theta})$. If $\left(\beta,\theta\right)\in\mathcal{I}_{s,\widetilde{s}}(2\widetilde{\omega}(s,\widetilde{s}))$, 
then ${\rm sign} (\widehat{\theta}^{\widehat{\widetilde\omega}})={\rm sign}(\theta)$, where $\widehat{\theta}^{\widehat{\widetilde\omega}}\triangleq(\widehat{\theta}_l\indic\{|\widehat{\theta}_l|>\mathbb{E}_n[Z_l^2]^{1/2}$ $\widehat{\widetilde\omega}(\widetilde{s})\})_{l=1}^{d_Z}$  on $\mathcal{G}_{\perp}\cap\underline{\mathcal{G}}_{\not\perp}\cap\mathcal{G}_{A1}$. 
Setting $\alpha^{\perp}=\alpha^{\not\perp}$ to $\alpha_n$, we have $\mathbb{P}(\mathcal{G}_{\perp}\cap\underline{\mathcal{G}}_{\not\perp}\cap\mathcal{G}_{A1})\ge 1-2\alpha_n-\alpha_n^{A1}\to1$. 
\end{thrm}
The first statement of Theorem \ref{th:nonvalid} provides a $1-\alpha^{\perp}-\alpha^{\not\perp}$ confidence band based on \eqref{eq:th:nonvalidCSinfty}. As for the STIV set in \eqref{CC}, uniformity in $\widetilde{c}$ and $c$ permits intersection over a grid. 
As in Example SE in Section \ref{sCS}, it provides 
$\widehat{S(\theta)}$ such that 
$\min_{s\in[d_Q],\widetilde{s}\in\left[d_Z-d_{\perp}\right]}\inf_{(\beta,\theta,\mathbb{P}):(\beta,\theta)\in\mathcal{I}_{s,\widetilde{s}}}
\mathbb{P}(\widehat{S(\theta)}\subseteq S(\theta))\geq1-\alpha^{\perp}-\alpha^{\not\perp}$.
The second statement of Theorem \ref{th:nonvalid} concerns model selection.  If the endogenous IVs induce a large enough violation of \eqref{einstr} then either superset or exact recovery of $S(\theta)$ is achieved. 
We provide NV-STIV rates of convergence in Section \ref{nvstiv}. 
The C-STIV in Section \ref{scstiv} 
is an extension estimating simultaneously $(\beta,\theta)$ and allowing for unknown $S_{\perp}$.

\section{Inference In Practice}\label{s8}
\subsection{Monte-Carlo}\label{sec:sim}
We study model \eqref{einstr}-\eqref{econstraints} with $\mathcal{B}=\mathbb{R}^{d_X}$ and $d_Q=d_X$, set $\beta^*=(1,-2,-0.5,0.25,0,\dots,0)^\top$ and let $Z$ be a standard Gaussian vector in $\mathbb{R}^{d_Z}$.
The exogenous regressors are the first $|S_I|$ IVs. For an endogenous $X_k$ we set $X_k=Z^\top\Pi_{\cdot,k}+\widetilde{U}_k$ where $\Pi\in\mathcal{M}_{d_Z,|S_I^c|}$ and $\widetilde{U}\in\R^{|S_I^c|}$. We let $(U(\beta^*),\widetilde{U}^\top)^\top$ be a mean zero Gaussian vector in $\mathbb{R}^{1+|S_I^c|}$ with variance having entries .05 but in the diagonal where the first entry is 1 and the others are $1-\pi$.
\setlength{\tabcolsep}{4pt}
\begin{table}[t]
\caption{STIV estimator}\label{mc_est_tab}
{\footnotesize
\begin{tabular}{lcccccccccccc}
\hline
\multicolumn{13}{c}{$d_Z=1500,d_X=1750,n=750,\pi=0.8$}\\
\hline
&\multicolumn{3}{c}{$c\widehat r=0.95$}&\multicolumn{3}{c}{$c\widehat r=0.75$}&\multicolumn{3}{c}{$c\widehat r=0.5$}&\multicolumn{3}{c}{$c\widehat r=0.25$}\\
\cmidrule(lr){2-4}\cmidrule(lr){5-7} \cmidrule(lr){8-10} \cmidrule(lr){11-13}
&p2.5& p50&p97.5&p2.5& p50&p97.5&p2.5& p50&p97.5&p2.5& p50&p97.5\\
\cmidrule(lr){2-4}\cmidrule(lr){5-7} \cmidrule(lr){8-10} \cmidrule(lr){11-13}
$\beta^*_1(=0.5)$&0.8	&	0.88	&	0.95	&	0.74	&	0.82	&	0.91	&	0.67	&	0.78	&	0.88	&	0.27	&	0.55	&	0.78
\\
$\beta^*_2(=-2)$&-1.9	&	-1.83	&	-1.75	&	-1.9	&	-1.83	&	-1.75	&	-1.9	&	-1.82	&	-1.74	&	-1.89	&	-1.81	&	-1.73
\\
$\beta^*_3(=-0.5)$&-0.41	&	-0.33	&	-0.26	&	-0.41	&	-0.33	&	-0.26	&	-0.41	&	-0.33	&	-0.26	&	-0.41	&	-0.32	&	-0.25
\\
$\beta^*_4(=0.25)$&0.01	&	0.08	&	0.16	&	0.01	&	0.08	&	0.16	&	0	&	0.08	&	0.16	&	0	&	0.08	&	0.16
\\
$\beta^*_5(=0)$&0	&	0	&	0	&	0	&	0	&	0	&	0	&	0	&	0	&	0	&	0	&	0
\\
$\beta^*_6(=0)$&0	&	0	&	0	&	0	&	0	&	0	&	0	&	0	&	0	&	0	&	0	&	0
\\
$\sigma^*(=1)$&1	&	1.05	&	1.1	&	1.01	&	1.06	&	1.12	&	1.02	&	1.07	&	1.13	&	1.05	&	1.12	&	1.2
\\
$\left|(\widehat\beta-\beta^*)_{S(\beta^*)}\right|_\infty$&0.15	&	0.2	&	0.25	&	0.16	&	0.21	&	0.27	&	0.17	&	0.23	&	0.33	&	0.23	&	0.45	&	0.73
\\
$\left|(\widehat\beta-\beta^*)_{S(\beta^*)^c}\right|_\infty$	&0	&	0	&	0.03	&	0	&	0	&	0	&	0	&	0	&	0.03	&	0	&	0.17	&	0.39
\\
\hline
$S(\widehat\beta)\supseteq S(\beta^*)$&\multicolumn{3}{c}{.98}&\multicolumn{3}{c}{.98}&\multicolumn{3}{c}{.98}&\multicolumn{3}{c}{.96} \\
$S(\widehat\beta)= S(\beta^*)$&\multicolumn{3}{c}{.62}&\multicolumn{3}{c}{.95}&\multicolumn{3}{c}{.91}&\multicolumn{3}{c}{.06}\\
\hline
\multicolumn{13}{c}{$d_Z=1500,d_X=1750,n=750,\pi=0.5$}\\
\hline
$\beta^*_1(=0.5)$&0	&	0.79	&	0.96	&	0	&	0.78	&	0.96	&	0	&	0.79	&	0.98	&	0	&	0.8	&	0.98

\\
$\beta^*_2(=-2)$&-1.93	&	-1.83	&	-1.5	&	-1.91	&	-1.83	&	-1.48	&	-1.9	&	-1.83	&	-1.47	&	-1.93	&	-1.83	&	-1.49

\\
$\beta^*_3(=-0.5)$&-0.39	&	-0.32	&	0	&	-0.4	&	-0.34	&	0	&	-0.4	&	-0.34	&	0	&	-0.4	&	-0.33	&	0

\\
$\beta^*_4(=0.25)$&0	&	0.08	&	0.18	&	0	&	0.07	&	0.18	&	0	&	0.08	&	0.16	&	0	&	0.07	&	0.15

\\
$\beta^*_5(=0)$&0	&	0	&	0	&	0	&	0	&	0	&	0	&	0	&	0	&	0	&	0	&	0

\\
$\beta^*_6(=0)$&0	&	0	&	0	&	0	&	0	&	0	&	0	&	0	&	0	&	0	&	0	&	0

\\
$\sigma^*(=1)$&1	&	1.09	&	3.3	&	0.99	&	1.09	&	3.42	&	0.99	&	1.08	&	3.4	&	0.98	&	1.08	&	3.38

\\
$\left|(\widehat\beta-\beta^*)_{S(\beta^*)}\right|_\infty$ & 0.15	&	0.24	&	1	&	0.15	&	0.23	&	1	&	0.15	&	0.23	&	1	&	0.15	&	0.24	&	1

\\
$\left|(\widehat\beta-\beta^*)_{S(\beta^*)^c}\right|_\infty$	&0	&	0.02	&	0.42	&	0	&	0.02	&	0.43	&	0	&	0.02	&	0.45	&	0	&	0.02	&	0.43

\\
\hline
$S(\widehat\beta)\supseteq S(\beta^*)$&\multicolumn{3}{c}{.82}&\multicolumn{3}{c}{.82}&\multicolumn{3}{c}{.80}&\multicolumn{3}{c}{.78} \\
$S(\widehat\beta)= S(\beta^*)$&\multicolumn{3}{c}{.49}&\multicolumn{3}{c}{.45}&\multicolumn{3}{c}{.47}&\multicolumn{3}{c}{.46}\\
\hline

\multicolumn{13}{p{107mm}}{\scriptsize \textbf{Notes:} 1000 replications. $\underline{r}_n=0.16$.}
\end{tabular}
}
\vspace{-5mm}
\end{table}
We set $\Pi$ so that $|\Pi_{\cdot,k}|_2^2=\pi$ for $k\in[|S_I^c|]$ and $\pi\in\{0.5,0.8\}$, hence all regressors have unit variance. Since the IVs are uncorrelated with one another, the concentration matrix $\Pi^\top\bold{Z}^\top\bold{Z}\Pi/(1-\pi)$ has diagonal elements close to $n\pi/(1-\pi)$ and the degree of endogeneity is $0.05/\sqrt{1-\pi}$ (see \cite{AS}). In low dimensions the IVs could be considered strong. However, the first-stage is not approximately sparse so even a first-stage Lasso would not be consistent (so it is impossible to estimate the concentration matrix and apply a method akin to 2SLS) and most of the IVs are weakly correlated with the endogenous regressors. We take $\Pi_{d_Z,1}=\Pi_{d_Z-1,2}=\cdots=\Pi_{d_Z-|S_I^c|+1,|S_I^c|} =\sqrt{3\pi/4}$. For the remaining entries we set 
$$\Pi_{i,j}=\begin{cases}-\sqrt{(\pi/4)/(d_Z-1)}&\text{i is odd and j is odd, or i is even and }j\geq d_Z/2\\ +\sqrt{(\pi/4)/(d_Z-1)}&\text{otherwise}\end{cases}$$
This means that there is one stronger IV and $d_Z-1$ weaker IVs for each endogenous regressor. Though each weaker IV accounts for a small fraction of their variance, collectively the weaker IVs account for fraction $\pi/4$. If $d_Z\leq d_X$ each IV has a stronger correlation with one regressor.
\setlength{\tabcolsep}{3pt}
\begin{table}[t!]
\caption{0.95 confidence sets and bands}\label{mc_des1_res}
{\footnotesize
\begin{tabular}{lcccccccccccccc}
\hline
\multicolumn{13}{c}{$d_Z=2050,d_X=50,n=2000,\pi=0.8$}\\
\hline
&\multicolumn{3}{c}{STIV}& SC 4 &SC 5&SC 6&SC 10  & ES  &\multicolumn{3}{c}{Bias-corrected STIV}& CB  \\
\cmidrule(lr){2-4} \cmidrule(lr){5-9} \cmidrule(lr){10-12} \cmidrule(l){13-13}
				&p2.5 & p50 & p97.5&\multicolumn{5}{c}{Median width/2}&p2.5 & p50 & p97.5&Width/2\\
\cmidrule(lr){2-4} \cmidrule(lr){5-9} \cmidrule(lr){10-12} \cmidrule(l){13-13}
$\beta^*_1(=1)$	    &0.9	&	0.95	&	0.99	&	0.8	&	1.02	&	1.32	&	6.07	&	0.33	&	0.94	&	1	&	1.06	&	0.1\\
$\beta^*_2(=-2)$    &-1.95	&	-1.9	&	-1.85	&	0.58	&	0.73	&	0.94	&	4.55	&	0.26	&	-2.04	&	-1.99	&	-1.95	&	0.07\\
$\beta^*_3(=-0.5)$ &-0.45	&	-0.4	&	-0.36	&	0.57	&	0.73	&	0.94	&	4.64	&	0.26	&-0.54	&	-0.49	&	-0.45	&	0.07\\
$\beta^*_4(=0.25)$    &0.11	&	0.15	&	0.19	&	0.57	&	0.73	&	0.95	&	4.65	&	0.26	&	0.2	&	0.24	&	0.29	&	0.07\\
$\beta^*_5(=0)$	   &0	&	0	&	0	&	0.8	&	1.02	&	1.31	&	6.03	&	0	&	-0.05	&	0	&	0.06	&	0.1\\
$\beta^*_6(=0)$    &0	&	0	&	0	&	0.57	&	0.73	&	0.95	&	4.62	&	0	&	-0.04	&	0	&	0.04	&	0.07\\
\hline
$S(\widehat\beta)\supseteq S(\beta^*)$&&1            &  Cover     &1             &1            & 1       &1  &.98 &&&&.94    \\
$S(\widehat\beta)=S(\beta^*)$          &&.98               &         &(.996,1)                &                &         &     & (.97,.98)    &&&&(.92,.95) \\
\hline
\multicolumn{13}{c}{$d_Z=49,d_X=50,n=2000,\pi=0.8$}\\
\hline
$\beta^*_1(=1)$	    &0.84	&	0.9	&	0.96	&	$\infty$	&	$\infty$	&	$\infty$	&	$\infty$	&	0.24	&	0.94	&	0.99	&	1.05	&	0.07\\
$\beta^*_2(=-2)$    &-1.96	&	-1.91	&	-1.87	&	$\infty$	&	$\infty$	&	$\infty$	&	$\infty$	&	0.2	&	-2.04	&	-2	&	-1.95	&	0.07\\
$\beta^*_3(=-0.5)$ & -0.47	&	-0.43	&	-0.39	&	$\infty$	&	$\infty$	&	$\infty$	&	$\infty$	&	0.2	&	-0.54	&	-0.5	&	-0.46	&	0.07\\
$\beta^*_4(=0.25)$    &0.13	&	0.18	&	0.23	&	$\infty$	&	$\infty$	&	$\infty$	&	$\infty$	&	0.2	&	0.2	&	0.25	&	0.3	&	0.07\\
$\beta^*_5(=0)$	   &0	&	0	&	0	&	$\infty$	&	$\infty$	&	$\infty$	&	$\infty$	&	0	&	-0.03	&	0.02	&	0.08	&	0.1\\
$\beta^*_6(=0)$    &0	&	0	&	0	&	$\infty$	&	$\infty$	&	$\infty$	&	$\infty$	&	0	&	-0.04	&	0	&	0.04	&	0.07\\
\hline
$S(\widehat\beta)\supseteq S(\beta^*)$&&1            &  Cover     &1             &1            & 1       &1  &1 &&&&.93    \\
$S(\widehat\beta)= S(\beta^*)$&&.96      &      &  (.996,1)                &                 &            &      &     &&&&(.91,.95) \\
\hline
\multicolumn{13}{c}{$d_Z=2050,d_X=50,n=2000,\pi=0.5$}\\
\hline
$\beta^*_1(=1)$	    &0.98	&	1.02	&	1.07	&	1.55	&	2.31	&	3.78	&	$\infty$	&	0.43&	0.92	&	1	&	1.07	&	0.12
\\
$\beta^*_2(=-2)$    &-1.95	&	-1.9	&	-1.86	&	0.85	&	1.26	&	2.08	&	$\infty$	&	0.27
&	-2.04	&	-1.99	&	-1.95	&	0.07
\\
$\beta^*_3(=-0.5)$ &-0.45	&	-0.4	&	-0.36	&	0.86	&	1.27	&	2.1	&	$\infty$	&	0.27
&	-0.54	&	-0.49	&	-0.45	&	0.07
\\
$\beta^*_4(=0.25)$    &0.11	&	0.16	&	0.2	&	0.85	&	1.26	&	2.08	&	$\infty$	&	0.26
	&	0.2	&	0.24	&	0.29	&	0.07
\\
$\beta^*_5(=0)$	   &0	&	0.03	&	0.07	&	1.56	&	2.32	&	3.78	&	$\infty$	&	0
&	-0.07	&	0	&	0.07	&	0.12
\\
$\beta^*_6(=0)$    &0	&	0	&	0	&	0.85	&	1.26	&	2.1	&	$\infty$	&	0
&-0.04	&	0	&	0.04	&	0.07
\\
\hline
$S(\widehat\beta)\supseteq S(\beta^*)$&&1             &  Cover     &1             &1            & 1       &1  &.75 &&&&.95    \\
$S(\widehat\beta)= S(\beta^*)$&&.13       &      &  (.996,1)                &                  &             &       & (.72,.77)    &&&&(.93,.96) \\
\hline
\multicolumn{13}{c}{$d_Z=49,d_X=50,n=2000,\pi=0.5$}\\
\hline
$\beta^*_1(=1)$	    &1	&	1.05	&	1.09	&	$\infty$	&	$\infty$	&	$\infty$	&	$\infty$	&	0.31
&	0.99	&	1.03	&	1.08	&	0.03
\\
$\beta^*_2(=-2)$    &-1.97	&	-1.93	&	-1.88	&	$\infty$	&	$\infty$	&	$\infty$	&	$\infty$	&	0.2

&-2.04	&	-1.99	&	-1.95	&	0.07
\\
$\beta^*_3(=-0.5)$ &-0.47	&	-0.42	&	-0.38	&	$\infty$	&	$\infty$	&	$\infty$	&	$\infty$	&	0.2

&-0.54	&	-0.49	&	-0.45	&	0.07
\\
$\beta^*_4(=0.25)$    &0.13	&	0.18	&	0.22	&	$\infty$	&	$\infty$	&	$\infty$	&	$\infty$	&	0.2

	&0.2	&	0.25	&	0.29	&	0.07
\\
$\beta^*_5(=0)$	   &0	&	0.05	&	0.09	&	$\infty$	&	$\infty$	&	$\infty$	&	$\infty$	&	0

&-0.03	&	0.03	&	0.09	&	0.1
\\
$\beta^*_6(=0)$    &0	&	0	&	0	&	$\infty$	&	$\infty$	&	$\infty$	&	$\infty$	&	0

&	-0.04	&	0	&	0.04	&	0.07
\\
\hline
$S(\widehat\beta)\supseteq S(\beta^*)$&&1             &  Cover     &1             &1            & 1       &1  &1 &&&&.50    \\
$S(\widehat\beta)= S(\beta^*)$&&.02       &      &  (.996,1)                &                  &            &      &     &&&&(.47,.53) \\
\hline
\multicolumn{13}{p{135mm}}{\scriptsize \textbf{Notes:} 1000 replications. `SC $s$' use sparsity certificate $s$. `ES' use estimated support. `CB' use $\varPhi=I$. SC/ES use one grid point for $c$. For $d_Z=49$ (resp. 2050) $\underline{r}_n=0.074$ (resp. 0.094). `STIV' uses $c=0.99/\widehat r$.  For SC/ES (resp. CB) `Cover' is the frequency with which $\beta^*$ lies in the bounds defined in \eqref{kbnds} (resp. \eqref{bndOmega}). 0.95 confidence intervals for the coverage are in parentheses (see \cite{wilson}). 
}
\end{tabular}
}
\vspace{-5mm}
\end{table}
\indent We construct 0.95 confidence sets and bands for $\beta^*$. For sets we use $\underline{r}_n$ from Class 3 with $\alpha=0.05$ and set $\widehat r=1.01\underline{r}_n$\footnote{This is possible under Assumption \ref{ass:Nemirovskir}, which permits $\widehat r=\underline{r}_n\sqrt{1+\tau_n}/(1-\tau_n)$ rather than  $\widehat r=\underline{r}_n|\bold{D_Z}\bold{Z}^\top|_\infty$, delivering a smaller value of $\widehat r$, which we find works better in practice.}. We consider sparsity certificates $4,5,6,7,10$. $\mathcal{I}_s$ is a singleton for each sparsity certificate and $d_X,d_Z$ below. We construct the bounds in \eqref{kbnds}, replacing $c>0$ with a grid, the construction of which is discussed below. For computational reasons (to allow sufficiently many replications) we limit the grid to at most two points. Using more points (and/or loss functions) could lead to narrower sets. We follow the same approach to construct the confidence set in \eqref{ESSB} based on an estimated support, taking $c$ equal to the first grid point and $S(\widehat\beta)$ to be the indices of the elements of $\bold{D}_{\bold{X}}^{-1}\widehat\beta$ with absolute value larger than $10^{-4}$. For the confidence bands we use STIV for the pilot and Class 4 to set $\widehat r$. Since the IVs are uncorrelated, the values of $\widehat r$ for classes 3 and 4 are nearly identical. We use the STIV estimator with $c=0.99/\widehat r$ for the pilot and set $\lambda=0.99$ and use Class 3 to estimate $\widehat\Lambda$. \vspace{0.2cm} 
\setlength{\tabcolsep}{4pt}
\begin{table}[t!]
\caption{0.95 confidence sets with $d_X>n$}\label{mc_des2_res}
{\footnotesize
\begin{tabular}{lcccccccccccccc}
\hline
\multicolumn{11}{c}{$d_Z=4100,d_X=4100,n=4000,\pi=0.8$}\\\hline
&\multicolumn{3}{c}{STIV}& SC 4 &SC 5&SC 6&SC 7&SC$^*$ 7&SC 10  & ES  \\
\cmidrule(lr){2-4} \cmidrule(lr){5-11}
				&p2.5 & p50 & p97.5&\multicolumn{7}{c}{Median width/2}  \\
\cmidrule(lr){2-4} \cmidrule(lr){5-11}
$\beta^*_1(=1)$	    &0.87	&	0.91	&	0.94	&	0.65	&	1.08	&	2.66	&$\infty^\dagger$&97.7 &	$\infty$	&	0.22
\\
$\beta^*_2(=-2)$    &-1.96	&	-1.93	&	-1.90	&	0.40	&	0.59	&	1.29	&$\infty^\dagger$&42.84&	$\infty$	&	0.17
\\
$\beta^*_3(=-0.5)$ & -0.46	&	-0.43	&	-0.39&	0.40	&	0.60	&	1.29	&$\infty^\dagger$&42.51&	$\infty$	&	0.17
\\
$\beta^*_4(=0.25)$    &0.14	&	0.18	&	0.21	&	0.40	&	0.59	&	1.29	&$\infty^\dagger$&41.89&	$\infty$	&	0.17
\\
$\beta^*_5(=0)$	   & 0&	0	&	0	&	0.65	&	1.09	&	2.68	&$\infty^\dagger$&97.37&	$\infty$	&	0
\\
$\beta^*_6(=0)$    &0	&	0	&	0	&	0.40	&	0.60	&	1.29	&$\infty^\dagger$&43.32&	$\infty$	&	0
\\
\hline
$S(\widehat\beta)\supseteq S(\beta^*)$&&1             & Cover&1&1&1&1&1&1&.97\\
$S(\widehat\beta)= S(\beta^*)$&&.99             &             &(.98,1)              &       &        &    &       &    & (.94,.99)\\\hline
\multicolumn{11}{c}{$d_Z=4100,d_X=4100,n=4000,\pi=0.5$}\\\hline
$\beta^*_1(=1)$	    &1.02	&	1.05	&	1.08	&	3.52	&	18.64	&	$\infty$	&$\infty$&$\infty$&	$\infty$	&	0.65
\\
$\beta^*_2(=-2)$    &-1.96	&	-1.93	&	-1.90	&	1.56	&	7.02&$\infty$&$\infty$&$\infty$&	$\infty$	&	0.4
\\
$\beta^*_3(=-0.5)$ & -0.46	&	-0.43	&	-0.4&	1.57	&	6.97	&	$\infty$	&$\infty$&$\infty$&	$\infty$	&	0.4
\\
$\beta^*_4(=0.25)$    &0.15	&	0.18	&	0.21&	1.55	&	6.99	&$\infty$&$\infty$&$\infty$&	$\infty$	&	0.4
\\
$\beta^*_5(=0)$	   & 0.02&	0.05	&	0.08	&	3.58	&	19.18	&	$\infty$	&$\infty$&$\infty$&	$\infty$	&	0
\\
$\beta^*_6(=0)$    &0	&	0	&	0	&	1.57	&	7.05	&	$\infty$	&$\infty$&$\infty$&	$\infty$	&	0
\\
\hline
$S(\widehat\beta)\supseteq S(\beta^*)$&&1             & Cover&1&1&1&1&1&1&.97\\
$S(\widehat\beta)= S(\beta^*)$&&.01            &             &(.98,1)              &      &        &    &       &   & (.94,.99)\\
\hline
\multicolumn{11}{p{120mm}}{\scriptsize \textbf{Notes:} 200 replications. `SC $s$' use sparsity certificate $s$. `ES' use estimated support. `CB' use $\varPhi=I$. SC/ES use one grid point for $c$. $\underline{r}_n=0.07$.`STIV' uses $c=0.99/\widehat r$. `Cover' is the frequency with which $\beta^*$ lies in the bounds defined in \eqref{kbnds}. 
$^\dagger$: The frequency of replications with sets of finite width is 0.03. $^*$: Sets using two grid points for $c$.}
\end{tabular}
}
\vspace{-5mm}
\end{table}
\setlength{\tabcolsep}{4pt}
\begin{table}
\caption{0.95 NV-STIV confidence sets for detection of endogenous IVs}\label{mc_des5_res}
\centering
{\footnotesize
\begin{tabular}{lcccccccccccccc}
\hline
\multicolumn{10}{c}{$d_Z=100,d_X=90,n=3000,\pi=0.8$}\\\hline
&\multicolumn{3}{c}{NV-STIV}& SC 4,1 &SC 4,2&SC 4,3&SC 4,5&SC 4,7 & SC 4,10\\
\cmidrule(lr){2-4} \cmidrule(lr){5-10}			
&p2.5 & p50 & p97.5&\multicolumn{6}{c}{Median width/2}  \\
\cmidrule(lr){2-4} \cmidrule(lr){5-10}	
$\theta_{89}(=0.8)$&0.55	&	0.6	&	0.65	&	0.53	&	0.53	&	0.53	&	0.53	&	0.54	&	0.54 
\\
$\theta_{90}(=0)$   & 0	&	0	&	0	&	0.53	&	0.53	&	0.53	&	0.53	&	0.54	&	0.54\\
\hline
$S(\widehat\theta)\supseteq S(\theta^*)$&&1 &Power&.95&.94&.94&.93&.92&.91\\
$S(\widehat\theta)=S(\theta^*)$&&1                   &             &(.93,.96)              &(.92,.95)        & (.92,.95)         & (.91,.94)    & (.9,.93)       & (.89,.92)   \\
                    \hline
\multicolumn{10}{p{120mm}}{\scriptsize \textbf{Notes:} 1000 replications. `SC $s,\widetilde{s}$' use sparsity certificates $s,\widetilde{s}$. SC use one grid point for $c$. $\underline{r}_n=0.07$ is from Class 3 with $\alpha=0.025$, and $\widehat r_1=1.01\underline{r}_{1,n}$. $\underline{r}^{\not\perp}_n=0.06$ is from Class 3 with $\alpha=0.025$. Confidence sets use a grid of 19 points for $\widetilde{c}$. `NV-STIV' uses $\widetilde{c}=0.99$. `Power' is the frequency with which the confidence sets do not include $\theta_{89}=0$.
}
\end{tabular}
}
\vspace{-5mm}
\end{table}

\indent {\bf Rule of Thumb for $c$.} We apply STIV with $c=\widehat r^{-1}$, corresponding to the least shrinkage. As $c$ decreases STIV is almost unchanged, until a point after which $\overline\sigma$ increases discontinuously. We recommend this for a single value of $c$. As $c$ decreases further, STIV is almost unchanged until a point after which there is another increase in $\overline\sigma$. This gives a second grid point, and so on. This rule means that we take the smallest $c$ (yielding the largest sensitivities) for each $\overline\sigma$.\vspace{0.1cm}

\indent {\bf Estimation.}  We consider the challenging setting with $n< d_Z<d_X$. We set $n=750$, $d_X=1750$ and $d_Z=1500$ and $S_I^c=\{1,5,1503,...,1750\}$, hence there are 250 endogenous regressors. Table \ref{mc_est_tab} reports the results. For sufficiently large $c$, STIV performs well in terms of selecting nonzero entries, and does not select those with values of zero. Due to the shrinkage, STIV is biased towards zero with bias decreasing in $c$. 
\vspace{0.1cm}

\indent {\bf Confidence Sets and Bands.} We set $n=2000$, $d_X=50$, $d_Z\in\{2050,49\}$, $S_I^c=\{1,5\}$ and make inference on $\beta^*$. This design is challenging since there are two endogenous regressors and either $d_Z<d_X$ or $n<d_Z$. We limit $d_X$ so as to permit application of all of our methods to the same design over 1000 replications. Below we modify the design to allow for $d_X>n$. 

Table \ref{mc_des1_res} reports the results. Sets based on a sparsity certificate are nested. If $d_Z=2050$, they can be informative on the sign of the first three entries of $\beta^*$. Though robust to identification, the sets can be conservative, are infinite if $d_Z=49$ and have coverage close to 1 if $d_Z=2050$. STIV performs well in selecting the nonzero parameters, resulting in less conservative sets based on estimated support. These are narrower than with sparsity certificate $s=|S(\beta^*)|=4$ as they use information on both the number and identities of relevant regressors. Coverage is below 0.95 when $d_Z=2050$ and $\pi=0.5$ because STIV using the rule of thumb value of $c$ can fail to distinguish $\beta^*_4=0.25$ from zero. In the other designs, the sets can be informative on the signs of the first three entries of $\beta^*$. The bias correction reduces the shrinkage and centers STIV on $\beta^*$. For $d_Z=2050$, there exists a sparse $\Lambda$ verifying \eqref{lInv}, with $|S_I|+|S_I^c|d_Z=4148$ nonzero entries out of $d_Xd_Z=102500$. The bands are narrower than the sets but have coverage slightly below 0.95 due to shrinkage when estimating $\Lambda$. For $d_Z=49$, there does not exist $\Lambda$ verifying \eqref{lInv}, leading to coverage below 0.95, significantly so for $\pi=0.5$.\vspace{0.1cm}

\indent {\bf Confidence Sets with $d_X>n$.} 
We set $n=4000$, $d_X=4100$, $S_I^c=\{1,5\}$ and $d_Z=4100$. Bands are infeasible since $\widehat\Lambda$ requires $d_Xd_Z=4100^2$ second-order cones. Table \ref{mc_des2_res} reports the results. If $\pi=0.8$, sets using a small sparsity certificate are informative on the signs. For $s=7$, the set is infinite if one grid point over $c$ is used but finite with two. STIV performs well in terms of selection, translating into narrower sets based on estimated support. Reducing the strength of the IVs ($\pi=0.5$) increases the width of the sets but coverage remains above 0.95.\vspace{0.1cm}

\indent {\bf Endogenous Instruments.} We take $n=3000$, $d_X=90$, $S_I^c=\{1,5\}$ and $d_Z=100$. There are 10 possibly endogenous IVs with indices $S_{\perp}^c=\{89,90,...,98\}$ and $Z_{89}=\sqrt{1-0.8^2}E+0.8U(\beta^*)$ is endogenous, where $E$ is an independent standard Gaussian. This preserves the variance of $Z$ but implies that $\theta^*$ has one nonzero entry given by $\theta^*_{89}=0.8$. There are as many known exogenous IVs as regressors. We apply the NV-STIV estimator, using STIV for the first stage and taking $\underline{r}_n$ from Class 3 with $\alpha=0.025$ and $\widehat r_1=1.01\underline{r}_n$. For the NV-STIV estimator we take $\underline{r}^{\not\perp}_n$ from Class 3 with $\alpha=0.025$. As both stages use $\alpha=0.025$ we construct 0.95 sets. We use sparsity certificates $s=4$ for $\beta$ and $\widetilde{s}\in[10]$ for $\theta$. The sets are intersected over a grid of $19$ points for $\widetilde{c}$. Table \ref{mc_des5_res} reports results. Due to shrinkage, NV-STIV is centred on 0.6. The endogenous IV is detected with frequency 0.95 for $\widetilde{s}=1$ and 0.91 for $\widetilde{s}=10$.

\subsection{EASI Demand System}\label{sec:appl}
The EASI demand system of \cite{LP} implies the vector of expenditure shares $S\in \R^{d_G}$ for ${d_G}$ goods consumed by a household satisfies 
\begin{align}
&S=\sum_{r=0}^{d_R}b_r T^r+C_1H+C_2 HT+A_0 P+\sum_{h=1}^{d_H}A_hPH_h+BPT+W\label{demeasi},\\
&T=\frac{1}{1-P^\top BP/2}\left(E-P^\top S+P^\top\left(A_0+\sum_{h=1}^{d_H}A_hH_h\right)P/2\right),\label{demy}
\end{align}
where $E\in\mathbb{R}$ is nominal expenditure, $T\in\mathbb{R}$ is deflated expenditure, $P\in\mathbb{R}^{d_G}$ is log-prices, $H\in\mathbb{R}^{d_H}$ is household characteristics, and $W\in\mathbb{R}^{d_G}$ are structural errors. Log-prices are normalized to be zero for a subset of households. The parameters 
are $b_r\in\mathbb{R}^{d_G}$ for $r=0,...,d_R$, $C_1,C_2\in\mathcal{M}_{d_G,d_H}$ and $A_0,...,A_{d_H},B\in\mathcal{M}_{d_G,d_G}$. Theory imposes restrictions such as  
(1) expenditure shares sum to one and (2) Slutsky symmetry, hence 
\begin{align}
&A_0,\dots,A_{d_H}\text{ and }B\text{ are symmetric};\ 1^\top b_0=1,1^\top C_1=1^\top C_2=0,1^\top B=0;\nonumber\\ 
&\forall r\in[d_R],1^\top b_r=0;\ \forall h\in[d_H],1^\top A_h=0.\label{eres}
\end{align}
Because $T$ depends on the parameters, the system \eqref{demeasi} is nonlinear, so difficult to estimate. \cite{LP} propose an approximate system, replacing $T$ with its first-order in prices approximation $D=E-P^\top S$, which is nominal expenditure deflated by a Stone price index. To reduce approximation error, we consider a second-order approximation and inject  
\begin{align}
\forall r\in\N,\ T^r=D^{r-1}\left(D+\frac{r}{2}P^\top\left(A_0+\sum_{h=1}^{d_H}A_hH_h+BD\right)P\right)+O(|P|_2^4)\label{demys}
\end{align}
(derived from \eqref{demy}) into \eqref{demeasi}. An approximation error arises due to the second term in \eqref{demys}, but it is small due to the normalization on log-prices. Our approximation depends on products of parameters, violating linearity. We replace each by a new parameter, restricted using \eqref{eres}. 
\subsubsection{Systems with Approximation Error}\label{approx}
Our results can be applied to estimate the system one equation at a time, ignoring cross-equation restrictions and approximation error. This does not make proper use of the underlying economic theory and would not allow a comparison with \cite{LP}. For this reason, we make some minor modifications to STIV. We allow for an  approximation error by adding an additional (unobserved) term $V(\beta)$ to the structural equation such that $\sigma_{V(\beta)}\leq v_{d_X}\to 0$. The practical implication is a minor modification to the IV-constraint, replacing $\widehat{r}\sigma$ with $\widehat{r}\sigma+(1+\widehat{r})\widehat{v}$ in \eqref{IVC}, where $\widehat{v}$ decays to zero with $n$. This allows for other models with approximation error including nonparametric IV (e.g., Example NP with approximation error) or when a fraction of the data is bracketed (in which case $\widehat{v}$ is random). To allow for a system of $d_G$ equations, the STIV objective function is summed and the IV-constraint is intersected over the equations, and $\alpha$ is replaced by $\alpha/d_G$. The latter allows the structural errors to be dependent across equations, and $d_G$ can depend on $n$. The bias correction and confidence bands are easily modified, and we also allow for approximation error in $\varPhi\beta$ (i.e., the function of interest is approximately linear). Further details and analysis are provided in Section \ref{sSQE}. 

\subsubsection{Implementation and Results}
We use the Canadian data of \cite{LP} for $n=4847$ rental-tenure single-member households with expenditure on rent, recreation and transportation. The $d_G=9$ goods are: food consumed at home, food consumed out, rent, clothing, household operation, household furnishing/equipment, transportation operation, recreation, and personal care. Individual characteristics are: age, gender, a dummy for car nonownership equal to one if real gasoline expenditure (at 1986 prices) is less than \$50, a social assistance dummy equal to one if government transfers are greater than 10 percent of gross income, and a linear time trend. Following \cite{LP}, we use $d_R=5$ for the degree of the expenditure polynomial. Each equation has ${d_X}=1570$ parameters. 
Log-prices are normalized to zero for residents of Ontario in 1986. The approximation error from the second order approximation is likely small because $\E_n[|P|_2^4]=0.0008$. In contrast, $\E_n[|P|_2^2]=0.0268$, and the mean share for 5 goods is less than 0.1, suggesting a large first-order approximation error. 
Since $D=E-P^\top S$ depends on $W$, the $|S_I^c|=963$ regressors which depend on $D$ are endogenous. We construct $d_Z=d_X$ IVs by replacing $D$ by $\overline{D}=E-P^\top\E_n[S]$ (i.e., replacing individual by average shares). \\
\indent The IVs are strong and ${d_Z}={d_X}$ and so we apply Section \ref{CI} to construct uniform 0.9 confidence bands for the Engel curves based on $d_{\varPhi}=9$ grid points. 
In the first step, we apply STIV, adjusting $\widehat r$ according to Class 4, taking $\alpha=0.05/d_G$, $c=0.99/\widehat r$ and $\widehat{v}_g=1/n$ for all $g\in[{d_G}]$. We choose $S_Q$ to exempt the constant, linear, and quadratic parts of the Engel curves ($b_0,b_1,b_2$) and the linear price parameters ($A_0$) from the penalty. It is reasonable to expect that the rest of the parameter be approximately sparse, particularly for the second-order approximation terms. 

\begin{figure}[t!]
\centering
\includegraphics[width=76mm,height=76mm]{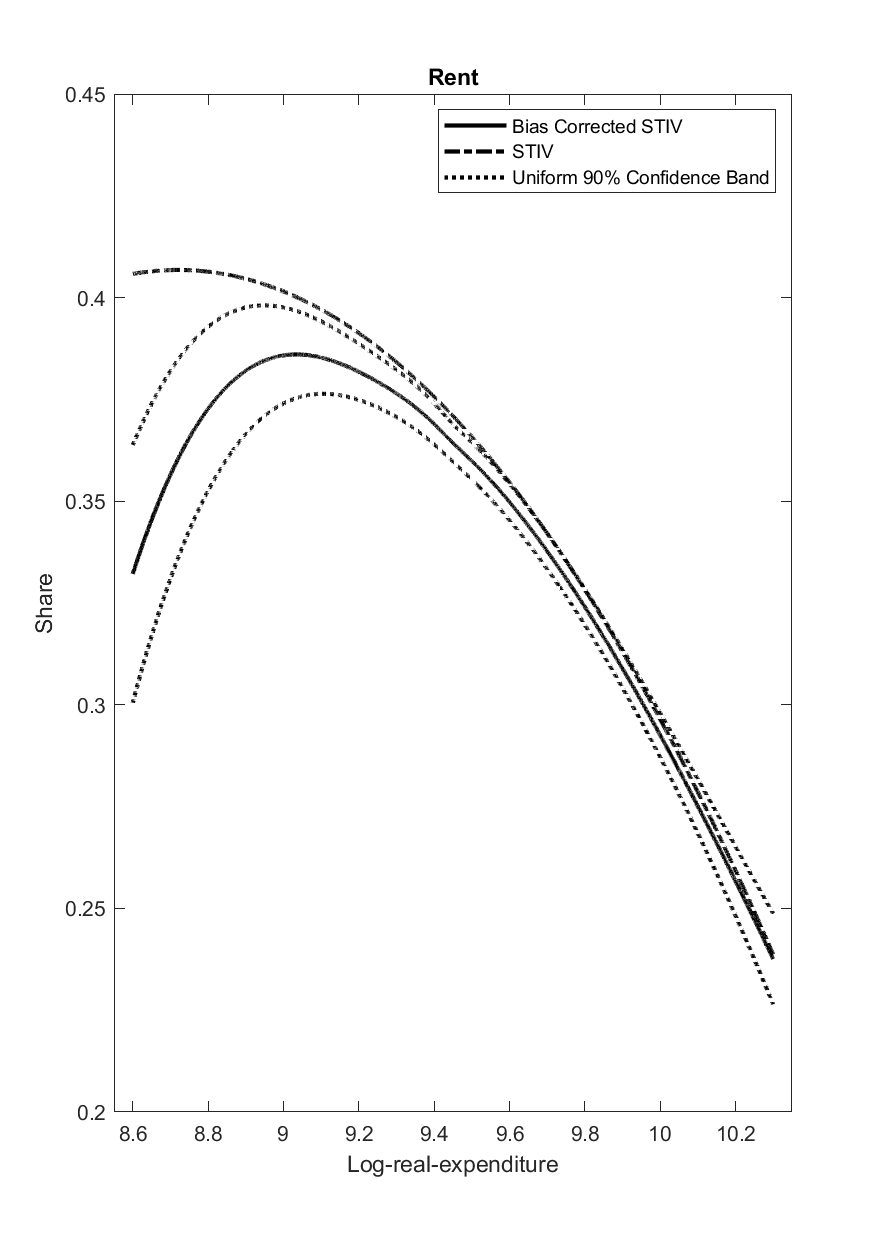}
\vspace{-4mm}
\caption{Engel Curve for Rent}\label{EC0}
\vspace{-1mm}
\end{figure}

For brevity, we do not present $\widehat\beta$ in full because it has $14,130$ elements. Instead, we summarize its support. Of $14,022$ parameters in $S_Q$, only 47 are estimated as nonzero, 22 of which are due to the second-order approximation. To build confidence bands for Engel curves, we obtain $\widehat\Lambda$ using $\underline{r}'_n$ from Class 3 with $\alpha=0.05$ and $\lambda=0.99$. Figure \ref{EC0} depicts the preliminary estimator for rent, its bias corrected counterpart and confidence bands. The second-order approximation yields a different curve to the that of \cite{LP}, which peaks at a higher expenditure level. The bias correction is large as the preliminary estimator lies outside the band. The band is wider at the end points, most likely due to lack of data. Engel curves for the remaining goods are available on request. The bias correction is large for household operation, clothing, personal care and transportation operation. The bands are marginally wider than those of \cite{LP} because we construct uniform bands rather than pointwise intervals and use a more flexible second-order approximation.

 \setcounter{equation}{0}  
 \setcounter{lmm}{0}
 \setcounter{prpstn}{0}
 \setcounter{rmrk}{0}
 \setcounter{dfntn}{0}
 \setcounter{thrm}{0}
 \setcounter{algm}{0}
 \setcounter{section}{1}
  \setcounter{subsection}{0}
 \setcounter{footnote}{0}
 \setcounter{figure}{0}
  \setcounter{table}{0}

%
%
%

\renewcommand{\theequation}{A.\arabic{equation}}
\renewcommand{\thelmm}{A.\arabic{lmm}}
\renewcommand{\thecrllr}{A.\arabic{crllr}}
\renewcommand{\thedfntn}{A.\arabic{dfntn}}
\renewcommand{\theprpstn}{A.\arabic{prpstn}}
\renewcommand{\thermrk}{A.\arabic{rmrk}}
\renewcommand{\thethrm}{A.\arabic{thrm}}
\renewcommand{\theass}{A.\arabic{ass}}
\renewcommand{\thealgm}{A.\arabic{algm}}
\renewcommand{\thesubsection}{A.\arabic{subsection}}
\renewcommand{\thetable}{A.\arabic{table}}
\section*{Appendix}
\subsection{Complements}
The proofs of the results below are in Section \ref{proofAp}. We denote by $F(b)\triangleq ZU(b)$,  $T(L)\triangleq\varPhi-LZX^\top$, 
and  $q_{A}$ (resp. $q_{A|\bold{M}}$) the quantile function of $A$ (resp. of $A$ given $\bold{M}$). When a random vector is a function of an estimator as in $U(\widehat{\beta})$, $\mathbb{E}_n$ is still used to denote $\sum_{i\in[n]}U_i(\widehat{\beta})/n$. 

\subsubsection{Complements on Section \ref{sSTIVintro}}\label{cs3} Proposition \ref{p4} relates the sensitivities (see also Section \ref{s:LBS}). 
\begin{prpstn}\label{p4} Let $S\in[d_X]$, $c>0$, and $\widehat{r}\le1$. For all $S_0\subseteq[d_X]$,  $q\in[1,\infty]$, and $\ell\in\mathcal{L}$,
\begin{enumerate}[\textup{(}i\textup{)}]
\item\label{p410} If $S\subseteq S_0\subseteq[d_X]$, $\widehat\kappa_{\ell,S}\ge\widehat\kappa_{\ell,S_0}$;
\item\label{p411} 
$\widehat\kappa_{\ell^q_{S_0},S}\ge\widehat\kappa_{\ell^q,S}$;
\item\label{p42i} 
$\max\left(\widehat{c}_{\kappa}(S)^{-1/q}\widehat\kappa_{\ell^\infty_{\widehat{S}(S)\cup S_0},S},|S_0|^{-1/q}\widehat\kappa_{\ell^\infty_{S_0},S}\right)\hspace{-.1cm}\le \widehat\kappa_{\ell^q_{S_0},S}\hspace{-.1cm}\le \widehat\kappa_{\ell^\infty_{S_0},S}$,\ $\widehat{c}_{\kappa}(S)^{-1}\widehat\kappa_{\ell^\infty_{\widehat{S}(S)},S}
\le \widehat\kappa_{\ell^1,S}$;
\item\label{p44i} 
$\widehat\kappa_{\ell^\infty_{S_0},S}
=\min_{k\in S_0}\widehat\kappa_{\ell_{k},S}=\min_{k\in S_0}\min_{\substack{\Delta\in \widehat{K}_{S}:\ \Delta_k=1,|\Delta_{S_0}|_{\infty}\le1}}\left|\widehat{\Psi}\Delta\right|_{\infty}$,
\end{enumerate}
where 
$\widehat{S}(S)\triangleq S\cup (S_Q^c\cap S(\widehat{\beta}))$ and $\widehat{c}_{\kappa}(S)\triangleq \gamma(c)(2|S\cap S_Q|+|S_Q^c\cap(S\cup S(\widehat{\beta}))|)$ if $0<c<1$ and else  
$\widehat{S}(S)\triangleq(S\cap S_Q)\cup ((S_Q^c\cup S_I^c)\cap(S\cup S(\widehat{\beta})))$ and 
$\widehat{c}_{\kappa}(S)\triangleq  \gamma(c\widehat{r})(
2|S\cap S_Q|+|S_Q^c\cap(S\cup S(\widehat{\beta}))|+c(1-\min(\widehat{r},1))|S_I^c\cap (S\cup S(\widehat{\beta}))|)$. 
\end{prpstn}

We emphasize 3 more baseline classes $\mathcal{P}$ which we further restrict when need be. Some confidence sets require very mild assumptions on $\mathcal{P}$ while deterministic bounds require working within subsets of these classes. 
The baseline classes are identification robust because they do not restrict the joint distribution of $(Z,X)$. 
Let, for $b\in\R^{d_X}$, 
 $\underline{\bold{D}}(b)$ be the diagonal matrix with positive diagonal elements $1/\widetilde{\sigma}_l(b)$ for $l\in[{d_Z}]$, where $\widetilde{\sigma}_l(b)^2\triangleq\mathbb{E}_n[F_l(b)^2]$, $\underline{\widehat{t}}(b)\triangleq|\underline{\bold{D}}(b)\bold{F}^\top(b)|_\infty/n$,   
$\underline{\mathcal{G}}\triangleq\left\{\underline{\widehat{t}}(\beta)\le \underline{r}_n\right\}$. The value of $\underline{r}_n$ for classes 1-3 is obtained using a union bound and the results in 
\cite{Pi,BGHK2,JSW}.\\  
{\bf Class 2:} 
\emph{
$\exists \mu_4>0:$ 
$\max_{l\in[{d_Z}]}\mathbb{E}[F_l(\beta)^4](\mathbb{E}[F_l(\beta)^2])^{-2}\le \mu_{4}$ and  ${d_Z}<\alpha\exp\left(n/\mu_4\right)/(2e+1)$.} We set  
$\underline{r}_n=\sqrt{2/(n/\log({d_Z}(2e+1)/\alpha)-\mu_4)}$.\\ 
{\bf Class 3:} 
\emph{
There exists $\delta$ in $(0,1]$ and $\mu_{2+\delta}>0$ such that
$$\displaystyle\left|\left(
\left(\E\left[|F_l(\beta)|^{2+\delta}\right]\right)
\left(\E\left[F_l(\beta)^2\right]\right)^{-(2+\delta)/2}\right)_{l\in[{d_Z}]}\right|_{\infty}
\le\mu_{2+\delta},$$ 
and
${d_Z}\le \alpha/(2\Phi(-n^{1/2-1/(2+\delta)}\mu_{2+\delta}^{-1/(2+\delta)}))$}. 
We set 
$\underline{r}_n=-\Phi^{-1}\left(\alpha/(2{d_Z})\right)/\sqrt{n}$.\\
Here $\mathbb{P}(\underline{\mathcal{G}})\ge1-\alpha-\alpha^B_n$ and $\alpha^B_n\triangleq \alpha C_1\mu_{2+\delta}\left(1+\sqrt{n}\underline{r}_n\right)^{2+\delta}n^{-\delta/2}$, where $C_1$ is an unknown universal constant, is a finite sample bound on coverage error. For classes 1 and 2, $\mathbb{P}(\underline{\mathcal{G}})\ge1-\alpha$, so we set $\alpha^B_n=0$. For classes 3-4, \eqref{cover} is modified to replace $1-\alpha$ by $1-\alpha-\alpha^B_n$.

\indent 
We use concentration arguments which involve 
$C_{{\rm N}}(m)\triangleq e(2\log(m)-1)$ for $m\ge3$ (Theorem 2.2 in \cite{DvdGVW}). For random $A\in\R^{d_A}$ and $B\in\R^{d_B}$ and sequences $M_A$, $M_A'$, and $M_{AB^{\top}}$ which can depend respectively on $d_A(d_A+1)/2$,   $d_A$, and $d_Ad_B$, 
denote by
\begin{align}
\mathcal{E}_{A}&\triangleq\left\{\left|D_A(\E_n-\E)\left[AA^{\top}\right]D_A\right|_{\infty}\ge \tau_n\right\},\ \mathcal{E}_{AB^{\top}}\triangleq\left\{\left|D_A(\E_n-\E)\left[AB^{\top}\right]D_B\right|_{\infty}\ge \tau_n\right\},\notag\\
\mathcal{E}_{A}'&\triangleq\left\{\min_{l\in[d_A]}\left(\bold{D}_{\bold{A}}^{-1}\right)_{l,l}\left(D_A\right)_{l,l}\le\sqrt{1-\tau_n}\ \text{or}\ \max_{l\in[d_A]}\left(\bold{D}_{\bold{A}}^{-1}\right)_{l,l}\left(D_A\right)_{l,l}\ge\sqrt{1+\tau_n}\right\}\notag,
\end{align}
\begin{enumerate}[\textup{(N.}i\textup{)}]
\item\label{Ni} $\E\left[\left|D_A\left(AA^{\top}-\E\left[AA^{\top}\right]\right)D_A\right|_{\infty}^2\right]\le M_A$,
\item\label{Niii} $\E\left[\left|D_A\left(AB^{\top}-\E\left[AB^{\top}\right]\right)D_B\right|_{\infty}^2\right]\le M_{AB^{\top}}$,
\item\label{Nii} $\E\left[\left|\left(A_{l}^2/\E\left[A_l^2\right]-1\right)_{l=1}^{d_A}\right|_{\infty}^2\right]\le M_A'$,
\end{enumerate}
$\alpha_{n}(A)\triangleq C_{{\rm N}}(d_A(d_A+1)/2)M_A/(n\tau_n^2)$, $\alpha_{n}(A)'\triangleq C_{{\rm N}}(d_A)M_A'/(n\tau_n^2)$, and $\alpha_{n}(AB^{\top})\triangleq C_{{\rm N}}(d_Ad_B)M_{AB^{\top}}/(n\tau_n^2)$. When $A$ depends on $\beta$ or $\Lambda$ but we omit it from the definition of $M_A$ and $\alpha_n(A)$, it means that the same sequence is used for all $\beta\in\mathcal{I}$, hence restricting $\mathcal{P}$. When $A$ in (N.\ref{Nii}) is a matrix, everything holds for the vectorization. 
\begin{lmm}\label{lconc}
Under (N.\ref{Ni}) for $A$ and $M_A$, $\mathbb{P}\left(\mathcal{E}_{A}\right)\le \alpha_{n}(A)$, under (N.\ref{Nii}) for $A$ and $M_A'$, $\mathbb{P}\left(\mathcal{E}_{A}'\right)\le \alpha_{n}(A)'$, and under (N.\ref{Niii}) for $A$ and $B$ and $M_{AB^{\top}}$, $\mathbb{P}\left(\mathcal{E}_{AB^{\top}}\right)\le \alpha_{n}(AB^{\top})$.
\end{lmm}

Taking $\widehat{r}=\underline{r}_n\left|\bold{D}_{\bold{Z}}\bold{Z}^{\top}\right|_{\infty}$ 
yields $\underline{\mathcal{G}}\subseteq\mathcal{G}$. Assumption \ref{ass:Nemirovskir} permits to work with the smaller $\widehat{r}=\underline{r}_n\sqrt{1+\tau_n}/(1-\tau_n)$.  The union bound used for $\underline{r}_n$ in classes 1-3 does not account for dependence over $l\in[d_Z]$ of $F_l(\beta)$, and so $\widehat{r}$ can be larger than necessary. 
To account for dependence, we consider Class 4 presented in Section \ref{s:LBS} under which 
$\widehat{r}=\left(q_{G|\bold{Z}}(1-\alpha)+2\zeta_n\right)/\sqrt{n},$ 
where 
$G\triangleq|\sqrt{n}\bold{D}_{\bold{Z}}\mathbb{E}_n[ZE]|_{\infty}$ and 
$\zeta_n\ge2\max\left(\tau_n/(1-\tau_n),(1/\sqrt{1-\tau_n}-1)\right) 
\log\left(2d_Z/\alpha_n\right)$. Section \ref{s:LBS} also points to useful results for dependent data. 

We now provide probabilistic conditions under which we can replace random quantities appearing in the right-hand sides in Proposition  \ref{t1} by deterministic ones. These are $\widehat{r}$, $\widehat{\sigma}(\beta)$, and the sensitivities. For classes 1-3 
we set 
\begin{equation}\label{GA1}
r_n\triangleq \underline{r}_n B_Z/\sqrt{1-\tau_n}\quad\text{and}\quad
\mathcal{G}_{A1}\triangleq \left\{\widehat{r}\le r_n\right\}\cap\mathcal{E}_{Z}'^c\cap\mathcal{E}_{X}'^c\cap\mathcal{E}_{ZX^{\top}}^c\cap\mathcal{E}_{U(\beta)}^c
\end{equation}
and, for Class 4, $r_n$ is defined in Section \ref{s:LBS} and $\mathcal{G}_{A1}\triangleq \left\{\widehat{r}\le r_n\right\}\cap\mathcal{E}_{Z}^c\cap\mathcal{E}_{X}'^c\cap\mathcal{E}_{ZX^{\top}}^c\cap\mathcal{E}_{U(\beta)}^c$. 
for all $n\in\N$, $\mathbb{P}\left(\mathcal{G}_{A1}\right)\ge1-\alpha^{A1}_n$.
We further restrict the class $\mathcal{P}$ and add: 
\begin{ass}\label{ass:Nemirovski} 
Let ${d_X},{d_Z}\ge3$, $\alpha_n$, $M_{U}$, $M_{ZX^{\top}}$, $M_{X}'$, $M_{Z}'$, 
and $B_Z$ (that can depend on $n$ and $d_Z$) positive. $\mathcal{P}$ is such that 
(N.\ref{Nii}) holds for $X$ and $M_{X}'$, (N.\ref{Niii}) holds for $Z$ and $X$ and $M_{ZX^{\top}}$.  
For $\mathcal{P}$ from class 1-3, we maintain (N.\ref{Ni}) holds for $U(\beta)$ and $M_{U}$, (N.\ref{Nii}) $Z$ and $M_{Z}'$, and 
$\mathbb{P}\left(\left|D_Z\bold{Z}^{\top}\right|_{\infty}> B_Z\right)\le \alpha_n$.
Moreover, 
\begin{equation}
\alpha^{A1}_n\triangleq\alpha^B_n+\alpha^C_n+\alpha_n(X)'+\alpha_n(ZX^{\top})\to0,\label{alphapsi}
\end{equation}
where $\alpha^C_n\triangleq\alpha_n+\alpha_n(Z)'+\alpha_n(U)$ for classes 1-3 and 
is defined in Section \ref{s:LBS} for Class 4. 
\end{ass}

If $D_ZZ$ is sub-Gaussian, $B_Z$ can be proportional to $\sqrt{\log(C nd_Z/\alpha_n)}$, where the constants $C$ and of proportionality depend on tail parameters of the sub-Gaussian distribution. Section \ref{s:LBS} presents the adjustments for classes 1-3 with assumptions \ref{ass:Nemirovski} and \ref{ass:Nemirovskir}. 

The population counterparts of $\widehat{K}_S$ and $\widehat{\overline{K}}_{S}$ 
replace 
$\widehat{r}$ in $\widehat{g}$ by $r_n$, which we denote by $g$:
\begin{align}
K_S&\triangleq\left\{\Delta\in\mathbb{R}^{d_X}:
1_{n}\left|\Delta_{S^c\cap S_Q}\right|_1\le \left|\Delta_{S\cap S_Q}\right|_1+c g(\Delta)\right\},\label{cone1}\\
\overline{K}_{S}&\triangleq\left\{\Delta\in\mathbb{R}^{d_X}:
1_{n}\left|\Delta_{S^c\cap S_Q}\right|_1
\le 2\left(\left|\Delta_{S\cap S_Q}\right|_1+cg(\Delta)\right)+\left|\Delta_{S_Q^c}\right|_1\right\}\notag,
\end{align}
\begin{lmm}\label{thrm:DLBsensitivities} 
On the event $\mathcal{G}_{A1}$, we have, for all $c>0$, 
\begin{align}
&\sigma_{U(\beta)}^2(1-\tau_n)\le\widehat{\sigma}(\beta)^2\le \sigma_{U(\beta)}^2(1+\tau_n),\notag\\
&
\forall S\subseteq[{d_X}],\ \ell\in\mathcal{L},\ \widehat\kappa_{\ell,S}\ge\frac{\kappa_{\ell,S}}{1+\tau_n}\left(1-\frac{\tau_n}{\kappa_{\ell^1,S}}\right),\ 
\widehat{\overline{\kappa}}_{\ell,S}\ge\frac{\overline{\kappa}_{\ell,S}}{1+\tau_n}\left(1-\frac{\tau_n}{\overline{\kappa}_{\ell^1,S}}\right),
\label{thrm:DLBsensitivities1b}\\
&\text{if}\ |S\cap S_Q|\leq s,\ \forall \ell\in\mathcal{L},\ \widehat\kappa_{\ell}(s)\ge\kappa_{\ell}(s)
\triangleq\frac{\kappa_{\ell}^0(s)}{1+\tau_n}\min_{S:|S\cap S_Q|\leq s}\left(1-\frac{\tau_n}{\kappa_{\ell^1,S}}\right),\label{thrm:DLBsensitivities1bomega}
\end{align}
where
$\kappa_{\ell}^0(s)$ 
is the population analogue of $\widehat\kappa_\ell(s)$. Under Assumption \ref{ass:Nemirovski}, for classes 1-4, we have $\mathbb{P}(\widehat{r}\le r_n)\ge1-\alpha_n^C$. 
\end{lmm} 

\begin{prpstn}\label{p5} 
We have, for all $S,S_0\subseteq[d_X]$, $q\in[1,\infty]$, and $k\in[d_X]$, 
\begin{enumerate}[\textup{(}i\textup{)}]
\item\label{p50i} $\kappa_{\ell^q_{S_0},S}\ge\kappa_{\ell^q,S}$,
\item\label{p5ii} $\max(c_{\kappa}(S)^{-1/q}\kappa_{\ell^\infty_{\overline{S}\cup S_0},S},|S_0|^{-1/q}\kappa_{\ell^\infty_{S_0},S})\le\kappa_{\ell^q_{S_0},S}\le\kappa_{\ell^\infty_{S_0},S}$,
$c_{\kappa}(S)^{-1}\kappa_{\ell^\infty_{\overline{S}},S}\le \kappa_{\ell^1,S}$,
\item\label{p5iib} $\kappa_{\ell^q_S,S}\le u_{\kappa}|S|^{1-1/q}\kappa_{\ell^1,S}$ under Condition IC,
\item\label{p501i} $\kappa_{\ell^\infty_{S_0},S}=\min_{k\in S_0}\kappa_{\ell_{k},S}$,
\item\label{p5i} 
$\kappa_{\ell^\infty_{S_0},S}\ge
\min_{k\in S_0}\max_{\lambda:|\lambda|_1\le1}\left(\lambda^{\top}\Psi_{\cdot,k}-(c_{\kappa}(S)-1)\max_{k'\ne k}|\lambda^{\top}\Psi_{\cdot,k'}|\right)$ if $\overline{S}\subseteq S_0$,
\item\label{p5iii} 
$\kappa_{\ell_{k},S}\ge \max_{\lambda:|\lambda|_1\le1}(\lambda^{\top}\Psi_{\cdot,k}+\max_{k'\ne k}|\lambda^{\top}\Psi_{\cdot,k'}|)(1+c_{\kappa}(S)\max_{k'\ne k}|\lambda^{\top}\Psi_{\cdot,k'}|/\kappa_{\ell^\infty_{\overline{S}},S})^{-1}$ and $\kappa_{\ell_{k},S}\ge \kappa_{\ell^{\infty}_{\overline{S}\cup\{k\}},S}$, 
\end{enumerate}
where 
$\overline{S}\triangleq(S\cap S_Q)\cup S_Q^c$ and $c_{\kappa}(S)\triangleq ((1+1_{n})|S\cap S_Q|+1_{n}|S_Q^c|)/\max(0,1_{n}-c)$ if $0<c<1_{n}$, else $\overline{S}\triangleq (S\cap S_Q)\cup S_Q^c\cup S_I^c$ and $c_{\kappa}(S)\triangleq((1+1_{n})|S\cap S_Q|+1_{n}|S_Q^c|+c(1-\min(r_n,1))|S_I^c|)/\max(0,1_{n}-c\min(r_n,1))$.\\
Moreover, if $\Delta\in K_S$ then 
\begin{equation}\label{galboundd}
|\Delta|_1\le c_{\kappa}(S)|\Delta_{\overline{S}}|_{\infty}. 
\end{equation}
The above statements hold if we replace $\kappa$ by $\overline{\kappa}$ and 
 $c_{\kappa}(S)$ by $c_{\overline{\kappa}}(S)$, the definition of which is the same but replacing $(1+1_n)$ by $(2+1_n)$ and $c$ by $2c$. We also have $\overline{\kappa}_{h,S}\ge \overline{\kappa}_{\ell^1,S}$.
\end{prpstn}
By item \eqref{p5iii} in Proposition \ref{p5} under Condition IC, we have
\begin{equation}\label{eellk}
\kappa_{\ell_{k},S}\ge \max\left(\sup_{\eta\in(0,1)}\eta\max_{\lambda\in \overline{S}_k(S,\eta)} 
\lambda^{\top}\Psi_{\cdot,k},\kappa_{\ell^{\infty}_{S\cup\{k\}},S}\right),
\end{equation}
where, for all $\eta\in(0,1)$ and $k\in[d_X]$, 
$$\overline{S}_k(S,\eta)\triangleq\left\{
\lambda\in\R^{d_Z}:|\lambda|_1\le1,1+\max_{k'\ne k}|\lambda^{\top}\Psi_{\cdot,k'}|/(\lambda^{\top}\Psi_{\cdot,k})\ge \eta 
\left(1+\max_{k'\ne k}\left|\lambda^{\top}\Psi_{\cdot,k'}\right|/\kappa_{\ell^1,S}\right)\right\}.$$

\subsubsection{Complements on Section \ref{CI}}\label{Lambdaest}
Estimation of $\Lambda$ is more computationally intensive than STIV because there are $d_{\varPhi}d_Z$ second-order cones (STIV has 1). For STIV we use the MOSEK solver, but if $d_{\varPhi}$ and $d_Z$ are very large, MOSEK can fail. For this reason we use an iterative procedure, which alternates between updating $\widehat\Lambda$ and $\widehat\nu$. Updating $\widehat\Lambda$ is more computationally demanding, so we apply FISTA with partial smoothing (\cite{BT09,BT12}). Details are in Section \ref{sFISTA}.

We now analyze $\widehat\Lambda$, which is a special case of the C-STIV estimator presented in Section \ref{scstiv}, applied to a system of $d_{\varPhi}$ equations. 
We also allow for approximation error, as in Section \ref{sSQE}. 
We denote by $G(b)\triangleq ZW(b)$ and 
\begin{align*}
\underline{\mathcal{G}}'\triangleq\left\{
\max_{f\in[d_{\varPhi}],k\in[{d_X}]}\frac{\left|\mathbb{E}_n\left[
T_{f,k}(\Lambda)\right]\right|}
{\mathbb{E}_n\left[T_{f,k}(\Lambda)^2\right]^{1/2}}\le \underline{r}'_n\right\}.
\end{align*}
The cones used to establish the rate of convergence of $\widehat{\Lambda}$ are sets of $\Delta'\in\mathcal{M}_{d_{\varPhi},{d_Z}}$ such that
\begin{align}
K_{S}'&\triangleq\left\{
\Delta':1_{n}(1-\lambda)\left|\Delta_{S^c}'\right|_1\le
(1+\lambda)
 \left|\Delta_S'\right|_1\right\},\ 
\overline{K}_{S}'\triangleq\left\{
\Delta': 1_{n}(1-\lambda)\left|\Delta_{S^c}'\right|_1\le
(2+\lambda) 
 \left|\Delta_S'\right|_1\right\},
 \notag
\end{align}
where $|\Delta_S'|_1$ follows the obvious modification to our notation 
in which one sums the absolute values of the entries $(k,l)\in S\subseteq [d_{\varPhi}]\times [d_Z]$ of $\Delta'$. We use $\kappa',\overline{\kappa}'$ to denote the population sensitivities using the cones above defined identically to $\kappa,\overline{\kappa}$, replacing $|\Psi\Delta|_{\infty}$ by $|\Delta'\Psi|_\infty$. Since $\widehat\Lambda$ can have more than one column, 
we use the operator norm from $\ell^{p}$ to $\ell^{q}$ which we denote by $|\cdot|_{p,q}$. 
We denote the population sensitivities for those losses 
by $\kappa_{\ell^{(p,q)},S}',\overline{\kappa}_{\ell^{(p,q)},S}'$. 
We also define $\overline{\rho}^{ZX}$ as $\widehat{\rho}^{ZX}$ replacing $\bold{D}_{\bold{Z}}$ (resp. $\bold{D}_{\bold{X}}$) by $D_{Z}$ (resp. $D_{X}$) and $\mathbb{P}(\beta,\Lambda)$ the distribution of 
$(X,Z,U(\beta),T(\Lambda), 
\Lambda ZU(\beta))$. 
\begin{ass}\label{assCI0} 
Let $M_{T}'$, 
$M_{G}$, $M_{2,G}$, $M_{EZX^\top}$, $q_2>0$, $(B_n)_{n\in\N}$ such that $B_n\ge1$, 
$(\rho^{ZX}_n)_{n\in\N}$, and $j\in[3]$ and a prior value of the parameter of Class $j$,  
such that, for all $(\beta,\Lambda,\mathbb{P})$ such that $(\beta,\Lambda)\in\mathcal{I}_{\varPhi}$ 
:  $\mathbb{P}(\beta,\Lambda)\in\mathcal{P}_{\varPhi}$, $q_1\in[2]$ and $n\in\N$, 
\begin{enumerate}[(i)]
\item\label{d.g.p.1} 
Assumption \ref{ass:Nemirovski} holds and 
$\mathbb{P}\left(\overline{\rho}^{ZX}>\rho^{ZX}_n
\right)\le \alpha_n$; 
\item\label{d.g.p.4} (N.\ref{Nii}) holds for $T(\Lambda)$ and $M_{T}'$;
\item\label{d.g.p.4a} (N.\ref{Ni}) holds for $G(\beta)$ and $M_{G}$; 
\item\label{d.g.p.4aa} (N.\ref{Ni}) holds for $\Lambda G(\beta)$ and $M_{2,G}$; 
\item\label{d.g.p.4b} $\left|\left(\max\left(\E\left[\left(\left(D_{\Lambda G(\beta)}\Lambda\right)_{f,\cdot} G(\beta)\right)^{2+q_1}\right],
\E\left[\left(\left(D_{\Lambda G(\beta)}\Lambda\right)_{f,\cdot} G(\beta)E\right)^{2+q_1}\right]
\right)\right)_{f=1}^{d_{\varPhi}}\right|_{\infty}\hspace{-.4cm}\le B_n^{q_1}$;
\item\label{d.g.p.4b2} $\max\left(\E\left[\left(\left|D_{\Lambda G(\beta)}\Lambda G(\beta)\right|_{\infty}/B_n\right)^{q_2}\right],
\E\left[\left(\left|D_{\Lambda G(\beta)}\Lambda G(\beta)E\right|_{\infty}/B_n\right)^{q_2}\right]\right)
\le 2$;
\item The distribution of $T(\Lambda)$ belongs to Class $j$ replacing $\alpha$ by $\alpha_n$; 
\item (N.\ref{Niii}) holds for $EZX^\top$ and $4M_{ZX^\top}$;
\item $\alpha^{A2}_n\to0$, 
where
\begin{align}
&\hspace{-1.2cm}\alpha^{A2}_n\triangleq2\zeta'_n+\zeta''_n+\varphi(d_Z,\tau_n)+\iota(d_{\varPhi},n)+C_{{\rm N}}(d_\varPhi(d_\varPhi+1)/2)M_{2,G}/(n\tau_n^2)+\alpha^{S}_n+\alpha^{BC}_n,\label{alphaPhi}\\
&\hspace{-1.1cm}\vspace{-.1cm}(\zeta'_n)^2\triangleq 3\alpha_n +\alpha_n(ZX^{\top})+\alpha^{S}_n+\alpha^{BC}_n
+\iota(d_{\varPhi},n),\quad 
\zeta''_n\triangleq \alpha_n+\alpha^{S}_n+\alpha^{BC}_n+\iota(d_{\varPhi},n),\notag\\
&\hspace{-1.1cm}\vspace{-.1cm}\iota(d,n)\triangleq C_2\left(\left(B_n^2\left(\log(dn)\right)^7/n\right)^{1/6}+\left(B_n^2(\log(dn))^3n^{-1+2/q_2}\right)^{1/3}\right)\quad \forall d\in\N,\notag
\end{align}

$\forall x\in(0,1),\varphi(d,x)\triangleq C_1x^{1/3}\max\left(1,\log(2{d}/x)\right)^{2/3}$, $C_1$ is constant and $C_2$ can depend on $q_2$, 
$\alpha^{BC}_n=\mathbb{P}(\mathcal{E}^c)$, $\mathcal{E}\triangleq\underline{\mathcal{G}}'\cap\mathcal{G}_{A1}\cap\{\overline{\rho}^{ZX}\le\rho^{ZX}_n\}\cap\mathcal{E}_{T}'^{c}\cap\mathcal{E}_{G(\beta)}^c\cap\mathcal{E}_Z^c$.  
\end{enumerate}
\end{ass}
\begin{prpstn}\label{tCLambda}
When $\mathcal{P}_{\varPhi}$ is such that Assumption \ref{assCI0}  \eqref{d.g.p.1}, 
\eqref{d.g.p.4}, \eqref{d.g.p.4a}, and (N.\ref{Ni}) holds for $Z$ and $M_{Z}$, then  
for all  $(\beta,\Lambda,\mathbb{P})$ such that $(\beta,\Lambda)\in\mathcal{I}_{\varPhi}$ and 
all solution $(\widehat\Lambda,\widehat\nu)$ of \eqref{Mhat} with $\lambda\in\left(0,1\right)$, we have, 
on $\mathcal{E}$, 
\begin{enumerate}[\textup{(}i\textup{)}]  
\item\label{tCLambdai}
For all $\ell\in\mathcal{L}$, $\ell\left(\left(\widehat{\Lambda}-\Lambda\right) D_{Z}^{-1}\right)\le 
\frac{2\underline{r}'_n\sqrt{1+\tau_n}\Sigma\left(\Lambda\right)\Gamma_{\kappa}'(S(\Lambda))}{1_{n}\kappa_{\ell,S(\Lambda)}'}$,\\
$\widehat{\nu}\le (1+\tau_n)\Sigma\left(\Lambda\right)\left(1+\frac{2\underline{r}'_n\rho^{ZX}_n\Gamma_{\kappa}'\left(S(\Lambda)\right)}{\lambda\kappa_{\ell^1_{S(\Lambda)},S(\Lambda)}'}\right)$,\\
where $
 \Gamma_{\kappa}'(S)\triangleq(1+\tau_n)\gamma\left(
\tau_n/\kappa_{\ell^{(\infty,\infty)},S}'+ \underline{r}'_n\rho^{ZX}_n(1+\tau_n)/(\lambda\kappa_{\ell^1_{S},S}')\right)$;
\item\label{tCLambdaii}
$\left|\left(\widehat{\Lambda}-\Lambda\right) D_{Z}^{-1}\right|_{1}\le \frac{2}{1_{n}}\min_{S\subseteq[d_{\varPhi}]\times[{d_Z}]}\max\left(\frac{\underline{r}'_n\sqrt{1+\tau_n}\Sigma\left(\Lambda\right)\Gamma_{\overline{\kappa}}'(S)}{\overline{\kappa}_{\ell^1,S}'},\frac{3+\lambda}{1-\lambda} \left|\Lambda_{S^c} D_{Z}^{-1}\right|_1\right),$
\begin{align*}
\hspace{-1cm}\widehat\nu\le (1+\tau_n)\Bigg(\Sigma(\Lambda)+\frac{\rho^{ZX}_n}{\lambda}\min_{S\subseteq[{d_X}]}\max&\left(2\Sigma(\Lambda)\left(\gamma\left(\frac{\underline{r}'_n\rho^{ZX}_n(1+\tau_n)}{\lambda\overline{\kappa}_{h,S}'}\gamma\left(\frac{\tau_n}{\overline{\kappa}_{\ell^1,S}'}\right)\right)-1\right),\right.\\
&\quad\left.\frac{3\left|\Lambda_{S^c} D_{Z}^{-1}\right|_1}{2\sqrt{1+\tau_n}}\right)\Bigg),
\end{align*}
where $\Gamma_{\overline{\kappa}}'(S)$, given by replacing $\kappa_{\ell^{(\infty,\infty)},S}',\kappa_{\ell^1_{S},S}'$ by $\overline{\kappa}_{\ell^{(\infty,\infty)},S},\overline{\kappa}_{h,S}'$.
\end{enumerate}
\end{prpstn}

We denote by $v^{\Lambda,\beta}_n$ 
and $v^{\Sigma(\Lambda)}_n$ the upper bounds on the right of \eqref{tCLambdai} and \eqref{tCLambdaii} (taking $|D_{\Lambda G(\beta)}\cdot|_{\infty,\infty}\sigma_{W(\beta)}$ for $\ell$, $|D_{\Lambda Z}\cdot|_{\infty,\infty}$ in Section \ref{conhomban}, and multiplying both sides by $|D_{\Lambda G(\beta)}|_{\infty}$ in case \eqref{tCLambdaii}) 
which can depend on $(\beta,\Lambda)$. 
For coverage guarantees 
we use: 
\begin{ass}\label{assCI1a} $\mathcal{P}_{\varPhi}$ is such that Assumption \ref{assCI0} holds and, 
for all $(\beta,\Lambda)\in\mathcal{I}_{\varPhi}$, we have
\begin{enumerate}[\textup{(}i\textup{)}]
\item\label{evD} $2v_n^D<1$;
\item\label{ezetaa} $\zeta_n\ge\max\left(2v_n^G,2v_n^T,
4v_n^D\log\left(2d_\varPhi/\alpha_n\right)
/(1-2v_n^D),v_n^R\right)$; 
\end{enumerate}
where $v_n^G\triangleq\sqrt{n} v^{\Lambda,\beta}_nr_n^E
+(|D_{\Lambda G(\beta)}\Lambda D_{Z}^{-1}|_{\infty,\infty}+v^{\Lambda,\beta}_n/\sigma_{W(\beta)})(v^{\beta}_n\tau_n+\sqrt{n}r_n^Ev_{d_X})$, $v_n^T\triangleq \sqrt{n}v^{\Lambda,\beta}_nr_n\sqrt{1-\tau_n^2}$, 
 $r_n^E\triangleq \underline{r}_n^E
2\log\left(2n/\alpha_n\right)B_Z\sqrt{1+\tau_n}$, 
$\underline{r}_n^E$ is obtained like $\underline{r}_n$ for Class 1 replacing $\alpha$ by $\alpha_n$ and $d_Z$ by $2d_Z$, and 
\begin{align*}
v^{D}_n&\triangleq v^{\Lambda,\beta}_nB_Z\sqrt{1+\tau_n}
+(\left|D_{\Lambda G(\beta)}\Lambda D_{Z}^{-1}\right|_{\infty,\infty}\hspace{-.2cm}+v^{\Lambda,\beta}_n/\sigma_{W(\beta)})(v^{\beta}_n\rho^{ZX}_n+B_Z\sqrt{1+\tau_n}v_{d_X})+\tau_n,\\
v^R_n&\triangleq\sqrt{n}\Big(\left|D_{\Lambda G(\beta)}\right|_{\infty}(\underline{r}'_n v^{\Sigma(\Lambda)}_n v^{\beta}_n\sqrt{1+\tau_n}+|\overline{V}(\beta)|_{\infty})\\
&\hspace{1.4cm} \left.+\left(\left|D_{\Lambda G(\beta)}\Lambda D_{Z}^{-1}\right|_{\infty,\infty}+v^{\Lambda,\beta}_n/\sigma_{W(\beta)}\right)
v_{d_X}\sqrt{1+\tau_n}\right)/(1-v^{D}_n).
\end{align*} 
\end{ass}
The coverage result that we obtain is more general than stated in Theorem \ref{tCI}. It is for approximately linear functions $\varPhi\beta+\overline{V}(\beta)$ for $\overline{V}(\beta)\in\R^{d_{\varPhi}}$ and is stated as
\vspace{-.1cm}
$$\mathbb{P}\left(\varPhi\beta+\overline{V}(\beta)\in \widehat{C}_{\varPhi}\right)\ge 1-\alpha-\alpha^{A2}_n.$$
\vspace{-.4cm}

We provide analysis under conditional homoskedasticity 
in Section \ref{conhomban}.  

\subsubsection{Complements on Section \ref{sec:endiv}}\label{nvstiv}
We denote by $\Sigma_{\not\perp}(\beta,\theta)\triangleq\max_{l\in S_{\perp}^c}((D_Z)_{l,l}\sigma_{T_l(\beta,\theta)})$,  $T(\beta,\theta)\triangleq ZU(\beta)-\theta$, 
and
\vspace{-.2cm}
\begin{align}
\underline{\mathcal{G}}_{\not\perp}\triangleq&\left\{
\max_{l\in S_{\perp}^c}\frac{\left|\mathbb{E}_n\left[T_l(\beta,\theta)\right]\right|}
{\mathbb{E}_n\left[T_l(\beta,\theta)^2\right]^{1/2}}\le \underline{r}^{\not\perp}_n\right\}.\label{G2}
\end{align}
\vspace{-.2cm}

We modify Assumption \ref{ass:Nemirovski} by replacing $\left\{\left|D_Z\bold{Z}^{\top}\right|_{\infty}> B_Z\right\}$ by $\left\{\left|(D_Z\bold{Z}^{\top})_{S_\perp,\cdot}\right|_{\infty}> B_Z\right\}$ and adding $\mathbb{P}\left(\widehat{\rho}^{ZX}_{S_{\perp}^c}>\rho^{ZX}_{S_{\perp}^c n}\right)\le \alpha_n$, (N.\ref{Nii}) for $(T(\beta,\theta))_{l\in S_{\perp}^c}$ and $M_{T}'$, and $\alpha_n(T)\to 0$. Theorems \ref{th:nonvalid} and \ref{th:nonvalid1} use  
vectors of functions 
which have $s,\widetilde{s}$ as arguments and we denote the evaluation using 4 arguments. 
The population sensitivities and their lower bounds replace $|\Psi\Delta|_\infty$ (resp. $r_n$) by $|(\Psi\Delta)_{S_\perp}|_\infty$ (resp. $r_n^{\perp}$), $\kappa_{\ell}(s)$ is defined in Lemma \ref{thrm:DLBsensitivities}, and $\kappa^{\Psi}(s)$ is defined similarly from 
 $\widehat\kappa^\Psi(s)$. 
\begin{thrm}\label{th:nonvalid1}
Let $s\in[d_Q]$, $\widetilde{s}\in\left[d_Z-d_{\perp}\right]$, $\widetilde{c}\in\left(0,1/\underline{r}^{\not\perp}_n\right)$, and $c>0$. If $\mathcal{P}_{\not\perp}$ is such that Assumption \ref{ass:Nemirovski} holds then, for all $(\beta,\theta,\mathbb{P})$ such that $\left(\beta,\theta\right)\in \mathcal{I}_{s,\widetilde{s}}$, and any NV-STIV estimator,
on $\mathcal{G}_{\perp}\cap\underline{\mathcal{G}}_{\not\perp}\cap\mathcal{G}_{A1}$, 
\vspace{-.3cm}
\begin{align}
\left| D_{Z}\left(\widehat{\theta}-\theta\right)_{S_{\perp}^c}\right|_\infty&\le \sqrt{1+\tau_n}\widetilde{\underline{\omega}}\left(s,\left|S\left(\theta\right)\right|,\beta,\theta\right),\label{eq:th:nonvalid:1}
\end{align}
\vspace{-.3cm}

\noindent where $\widetilde{\underline{\omega}}(s,\widetilde{s},\beta,\theta)=
2
\gamma\left(\underline{r}^{\not\perp}_n\widetilde{s}/\widetilde{c}\right)
\left(r^{\not\perp}_n\Sigma_{\not\perp}\left(\beta,\theta\right)/1_n+\delta^{\Psi}(s,\beta)\right)$, $\delta^{\Psi}(s,\beta)\triangleq \delta(s,\beta)/\kappa^{\Psi}(s)$, and  
$\delta(s,\beta)\triangleq 2r^{\perp}_n\gamma\left(r^{\perp}_n/\kappa_g(s)\right)\sigma_{U(\beta)}\sqrt{1+\tau_n}
(1+2r^{\perp}_n\Gamma_\kappa(S(\beta))/(c\kappa_{\ell^1_{S(\beta)\cap S_Q},S(\beta)}))$.
\end{thrm}
By  (N.\ref{Niii}) for $Z,X$ (which is part of Assumption \ref{ass:Nemirovski}) and the computations in 
\eqref{eNPIV0}, on $\mathcal{G}_{A1}$, $\kappa^{\Psi}(s)\ge \kappa_{\ell^1}(s)(1-\tau_n)/(|\Psi|_{\infty}+\tau_n)$. Recall also that $\kappa_g(s)=\kappa_{\ell^1}(s)$ under Condition IC, otherwise it depends on the LP used to compute $\widehat\kappa_g(s)$. 
For the second statement in Theorem \ref{th:nonvalid}, $\widetilde{\omega}\left(s,\widetilde{s},\beta,\theta\right)$ is obtained by replacing $\widehat{\delta}$ and $\widehat{\delta}^{\Sigma}$ in the definition of $\widehat{\widetilde{\omega}}\left(\widetilde{s}\right)$ by their deterministic upper bounds. For $\widehat{\delta}^{\Sigma}$ we use $\rho^{ZX}_{S_\perp^c,n} \delta(s,\beta)/\kappa_{\ell^1}(s)$. 
The deterministic upper bounds on $\widehat{\delta}$ and $\widehat{\delta}^{\Sigma}$ hold on $\mathcal{G}_{\perp}\cap \underline{\mathcal{G}}_{\not\perp}\cap\mathcal{G}_{A1}$ and are obtained using Lemma \ref{sigmabounds}.

\subsubsection{Systems of Equations with Approximation Errors}\label{sSQE}
To allow for approximation error we use $U(\beta)=W(\beta)+V(\beta)$ and suppose that \eqref{einstr} holds with $W(\beta)$ in place of $U(\beta)$ and $V(\beta)$ is a small approximation error, for which we assume that 
$\sigma_{V(\beta)}\leq v_{d_X}$, for $v_{d_X}$ decaying to zero with $d_X$.  
The assumptions previously made on $\left(X,Z,U(\beta)\right)$ are made on $\left(X,Z,W(\beta)\right)$ and $\mathcal{I}$ is modified accordingly and incorporates $\sigma_{V(\beta)}\leq v_{d_X}$. 
The model with approximation errors allows for the structural equation
\begin{equation}Y=f(\widetilde{X})+W,\quad\mathbb{E}[W|\widetilde{Z}]=0,\label{npiv}\end{equation}
where $f\in\mathcal{S}$, 
and for functions $(g_k^X)_{k\in\N}$ and a decreasing  sequence $(v_{d_X})_{d_X\in\N}$, 
\begin{equation}
\forall d_X\in\N,\sup_{g\in\mathcal{S}}\inf_{b\in\R^{d_X}}\E\left[\left(g\left(\widetilde{X}\right)-\sum_{k=1}^{d_X}g_k^X\left(\widetilde{X}\right)b_k\right)^2\right]\le v_{d_X}^2.\label{eboundNP2}
\end{equation}
The rate of decay of $(v_{d_X})_{d_X\in\N}$ is usually taken slow so $\mathcal{S}$ can be large.
It corresponds to minimum smoothness but $f$ can lie in a class of smoother functions. 
 The model with approximation error involves 
 $X=(g^X_1(\widetilde{X}),\dots,g^X_{{d_X}}(\widetilde{X}))$, $V(\beta)=f(\widetilde{X})-\sum_{k=1}^{{d_X}}g_k^X(\widetilde{X})\beta_k$, 
and $\sigma_{V(\beta)}\le v_{{d_X}}$. 
$V(\beta)$ is the error made by approximating the function in the high-dimensional space, and $v_{d_X}=o(n^{-1/2})$ 
for ${d_X}$ large enough. For well chosen classes $\mathcal{S}$ and functions $(g_k^X)_{k\in\N}$, the vector $\beta\in\mathcal{I}$ is approximately sparse. 
We use IVs which are 
functions of $\widetilde{Z}$. 
We consider a system where $Y,U(\beta),V(\beta),W(\beta)\in\mathbb{R}^{d_G}$, $U_g(b)\triangleq Y_g-X^{\top}b_{\cdot,g}$, 
$\beta\in\mathcal{M}_{d_X,d_G}$, $S_Q,S_I\subseteq [d_X]\times[d_G]$, and $\sigma_{V_g(\beta)}\le v_{g,d_X}$ for all $g\in[d_G]$. This is the setup of Section \ref{sec:appl} where $X$ is used in all equations. Else, a simple modification applies. We now define: 

\begin{dfntn} For $c,\widehat{v}>0$, the
E-STIV estimator $\left(\widehat{\beta},\widehat{\sigma}\right)$ is
any solution of 
\begin{equation}\label{ESTIV}
\min_{b\in \widehat{\mathcal{I}}_{E}(\widehat{r},\sigma),\sigma\ge0}\left(\left|\bold{D}_{\bold{X}}^{-1}b_{S_Q}\right|_1+c|\sigma|_1\right),
\end{equation}
where, setting for all $b\in\mathcal{M}_{d_X,d_G}$ and $g\in[d_G]$, $\widehat{\sigma}_g(b)^2\triangleq\E_n[U_g(b)^2]$, 
\begin{equation*}
\widehat{\mathcal{I}}_{E}(\widehat{r},t)\triangleq\left\{b\in\mathcal{B},\forall g\in[{d_G}],\left| \bold{D}_{\bold{Z}}\mathbb{E}_n[ZU_g(b)]\right|_{\infty}\le \widehat{r}\sigma_g +(\widehat{r}+1)\widehat{v}_g,\widehat{\sigma}_g(b)\le \sigma_g\right\}.
\end{equation*} 
\end{dfntn}
For a 
nonparametric model \eqref{npiv} one can take $\widehat{v}=\sqrt{1+\tau_n}v_{d_X}$. 
The E-STIV can also be used 
when, for $i\in B\subseteq[n]$, the outcomes are bracketed. Then, we let for $i\in B^c$ $y_i$ be the observed outcome and $v_i=0$, while, for $i\in B$, $y_i$ is the midpoint of the bracket. 
One has $|v_i|\le e_i$, where $e_i$ are half-widths of the brackets, and we let $\widehat{v}=n^{-1}\sum_{i\in B}e_i^2$ and $v_{d_X}^2=\E_n\E\left[\indic{\{i\in B\}}e_i^2\right]$.  
With $d_G$ equations, we allow for cross-equation restrictions, and the number of equations $d_G$ can depend on $n$. The E-STIV estimator is used 
in Section \ref{sec:appl}. 

To allow for approximation error, we modify $\mathcal{P}$ so that $W(\beta)$ plays the role of $U(\beta)$. For simplicity we only analyze classes 1-3. We choose $r_n$ and $\underline{r}_n$ as in Section \ref{cs3} replacing $\alpha$ by $\alpha/{d_G}$ and use
$
\mathcal{G}\triangleq\{\max_{g\in[{d_G}],\ l\in[{d_Z}]}\left|\mathbb{E}_n\left[Z_{l}W_g(\beta)\right]\right|
\mathbb{E}_n[Z_{l}^2]^{-1/2}\mathbb{E}_n\left[W_g(\beta)^2\right]^{-1/2}\le r_n \}$ 
and 
$\mathcal{G}_{A1}$ is defined in Section \ref{cs3} replacing $\mathcal{E}_{U(\beta)}^c$ by $\mathcal{E}_{V(\beta)}^c\cap\mathcal{E}_{W(\beta)}^c$, where  the probability $\alpha_n(U)$ in the definition of $\alpha^C_n$ is replaced by $2{d_G}\alpha_n(U)$.  
The population sensitivities are obtained replacing 
$\left|\Psi\Delta\right|_{\infty}$ by  
$
\sum_{g=1}^{d_G}\left|\Psi\Delta_{\cdot,g} \right|_{\infty}$, $\widehat K_S$ and $\widehat{\overline{K}}_{S}$ by 
$
K_S\triangleq\{\Delta\in\mathbb{R}^{d_X}:
1_{n}|\Delta_{S^c\cap S_Q}|_1\le |\Delta_{S\cap S_Q}|_1+c g(\Delta)\}$ and 
$\overline{K}_{S}$ where the right-hand side is $2(|\Delta_{S\cap S_Q}|_1+cg(\Delta))+|\Delta_{S_Q^c}|_1$, $g(\Delta)\triangleq r_n(\beta)|\Delta_{S_I}|_1+|\Delta_{S_I^c}|_1$, and $r_n(\beta)\triangleq\max_{g\in[{d_G}]}\min(r_n\hspace{-.05cm}+\hspace{-.05cm}(r_n\hspace{-.05cm}+\hspace{-.05cm}1)\max(0,1_{n}\sigma_{W_g(\beta)}/v_{g,d_X}-1)^{-1}\hspace{-.1cm},1)$ 
replaces $r_n$ in 
the sensitivities. 
\begin{prpstn}\label{ThoracleSE} 
For all $(\beta,\mathbb{P})$ such that $\beta\in\mathcal{I}$, assuming as well
 $\E_n[v_{g,d_X}^2]\le \widehat{v}_g^2$ on $\mathcal{G}_{A1}$ and all solution $\left(\widehat{\beta},\widehat\sigma\right)$ of \eqref{ESTIV}, the following hold on $\mathcal{G}\cap\mathcal{G}_{A1}$  
\begin{enumerate}[\textup{(}i\textup{)}]  
\item\label{ThoracleSEi}
For a sparse matrix $\beta$, for all $\ell\in\mathcal{L}$, we have
\begin{align*}
&\hspace{-.5cm}
 \ell\left( D_X^{-1}\left(\widehat{\beta}-\beta\right)\right)
\le 
\frac{2r_n}{1_{n}\kappa_{l,S(\beta)}}\left(\sum_{g=1}^{d_G}\sigma_{W_g(\beta)}+\left(r_n+2\right)v_{g,d_X}\right)\Gamma_{\kappa}(S(\beta));
\end{align*}
\item\label{ThoracleSEii}
For all $S,S_0\in[d_X]^{d_G}$, 
and $q\in[1,\infty]$, we have
\begin{align*}
&\hspace{-.7cm}
\left| D_X^{-1} \left(\widehat{\beta}-\beta\right)_{S_0}\right|_q
\le \max\left(
\frac{r_n}{1_{n}\overline{\kappa}_{\ell^q_{S_0},S}}\left(\sum_{g=1}^{d_G}\sigma_{W_g(\beta)}+\left(r_n+2\right)v_{g,d_X}\right)
\Gamma_{\overline{\kappa}}(S),\frac{6}{1_{n}}\left| D_X^{-1}\beta_{S^c\cap S_Q}\right|_1\right).
\end{align*}
\end{enumerate}
\end{prpstn}
In a model with $V_{g}(\beta)=0$, we take $\widehat{v}_g=0$ and can derive the same results as for the STIV estimator, including the confidence sets. 
The confidence bands of Section \ref{CI} are easily adapted. For the equation $g$ confidence band $\widehat{C}_{\varPhi,g}$, we use E-STIV for $\widehat\beta$ and replace $\widehat{\varPhi\beta}$ by $\widehat{\varPhi\beta}_{\cdot,g}$ in equation \eqref{bndOmega}. Assumption \ref{assCI1a} uses
$
\max_{g\in[d_G]}| D_{X}^{-1}(\widehat{\beta}-\beta)_{\cdot,g}|_1\le v^{\beta}_n$ 
and replaces quantities on the right-hand side which are specific to equation $g$ by the maximum over $g\in[d_G]$.  
Theorem \ref{tCI} is modified to replace $\mathbb{P}(\varPhi\beta\in \widehat{C}_{\varPhi})\ge 1-\alpha-\alpha^{A2}_n$ by $\mathbb{P}(\varPhi\beta_{\cdot,g}+\overline{V}_g(\beta)\in \widehat{C}_{\varPhi,g})\ge 1-\alpha-\alpha^{A2}_n\quad \forall g\in[{d_G}]$. The term $\varPhi\beta_{\cdot,g}+\overline{V}_g(\beta)$ for $\overline{V}_g(\beta)\in\R^{d_{\varPhi}}$ is the equation $g$ approximately linear function (see the proof of Theorem \ref{tCI}). The proof of Theorem \ref{tCI} allows for approximation error and extension to systems is straightforward. Proposition \ref{ThoracleSEO} considers losses useful for a system of nonparametric IV equations and rates of estimation of $\sigma_{W_g(\beta)}$ for the confidence bands under conditional homoskedasticity in Section \ref{conhomban}.

\subsection{Proofs of the Results in the Main Text}\label{s93} 
\noindent{\bf Proof of Proposition \ref{t1}.} First prove the first inequality. 
Take $\beta\in\mathcal{I}$ and set $\widehat{\Delta}\triangleq\bold{D}_{\bold{X}}^{-1}(\widehat\beta -\beta)$.
By definition of $\widehat{\mathcal{I}}$ and $\widehat{\sigma}(\beta)=\mathbb{E}_n[U(\beta)^2]$, 
on $\mathcal{G}$, we have 
$\beta\in\widehat{\mathcal{I}}\left(\widehat{r},\widehat{\sigma}(\beta)\right)$. Also, on $\mathcal{G}$, 
\begin{align}
\left|\widehat{\Psi}\widehat{\Delta} \right|_{\infty} &\le \left|\bold{D}_{\bold{Z}}\mathbb{E}_n[ZU(\widehat\beta)]\right|_{\infty}
+\left|\bold{D}_{\bold{Z}}\mathbb{E}_n[ZU(\beta)]
\right|_{\infty}\le
\widehat{r}\left(\widehat{\sigma}+\widehat{\sigma}(\beta)\right).\label{stop1}
\end{align}
Also,
$(\widehat\beta,\widehat{\sigma})$ minimizes the criterion
$\left|\bold{D}_{\bold{X}}^{-1}\beta\right|_1 +c\sigma$. 
Thus, on $\mathcal{G}$, we have
\begin{equation}\label{eq:main}
\left|\bold{D}_{\bold{X}}^{-1}\widehat{\beta}_{S_Q}\right|_1 +c\widehat{\sigma}\le
|\bold{D}_{\bold{X}}^{-1}\beta_{S_Q}|_1+c\widehat{\sigma}(\beta).
\end{equation}
This implies, on $\mathcal{G}$,
\begin{align}
\left|\widehat{\Delta}_{S(\beta)^c\cap S_Q}\right|_1&=\sum_{k\in S(\beta)^c\cap S_Q}
\left|\mathbb{E}_n[X_k^2]^{1/2}\widehat\beta_k\right|\label{in0}\\
&\le\sum_{k\in S(\beta)\cap S_Q}\left(
\left|\mathbb{E}_n[X_k^2]^{1/2}\beta_k\right|
-\left|\mathbb{E}_n[X_k^2]^{1/2}\widehat\beta_k\right|\right)+c\left(\widehat{\sigma}(\beta)-
\widehat{\sigma}
\notag
\right)\\
&
\le\left|\widehat{\Delta}_{S(\beta)\cap S_Q}\right|_1+c\left(\widehat{\sigma}(\beta)-
\widehat\sigma\left(\widehat\beta\right)
\right).\nonumber
\end{align}
The last inequality holds because by construction $\widehat\sigma(\widehat\beta)\le\widehat\sigma$. Since $b\to\sqrt{\widehat{\sigma}(b)}$ is convex and $$w_*\triangleq- \mathbb{E}_n[XU(\beta)]\mathbb{E}_n[U(\beta)^2]^{-1/2}\indic{\left\{\mathbb{E}_n[U(\beta)^2]\ne0\right\}}\in \partial \widehat{\sigma}(\cdot)(\beta).$$  
we have
$\widehat{\sigma}(\beta)-\widehat\sigma\left(\widehat\beta\right)\le w_*^{\top}\left(\beta-\widehat\beta\right)=\left(\bold{D}_{\bold{X}}w_*\right)^{\top}\bold{D}_{\bold{X}}^{-1}\left(\beta-\widehat\beta\right)=-\left(\bold{D}_{\bold{X}}w_*\right)^{\top}\widehat{\Delta}$.\\ 
Now, for all $k\in S_I$, we have 
$\left|\left(\bold{D}_{\bold{X}}w_*\right)_k\right|\le \widehat{r}$ on $\mathcal{G}$. This is because these regressors serve as their own IV and, on $\mathcal{G}$, $\beta\in\widehat{\mathcal{I}}\left(r_n,\widehat{\sigma}(\beta)\right)$. Also, for all rows of index $k$ in the set $S_I^c$,
$$
\left|\left(\bold{D}_{\bold{X}}w_*\right)_k\right|\le
\left|\mathbb{E}_n[X_{k}U(\beta)]\right|\mathbb{E}_n[X_k^2]^{-1/2}\mathbb{E}_n[U(\beta)^2]^{-1/2}\le 1
$$
due to the Cauchy-Schwarz inequality. Finally, we obtain 
\begin{align}\label{dubovic}
\widehat{\sigma}(\beta)
-\widehat\sigma\left(\widehat\beta\right)\le \widehat{r}\left|\widehat{\Delta}_{S_I}\right|_1+\left|\widehat{\Delta}_{S_I^c}\right|_1.
\end{align} 
Combining \eqref{dubovic} with (\ref{in0}), on $\mathcal{G}$ we have
$\widehat{\Delta}\in \widehat{K}_{S(\beta)}$. Using
\eqref{stop1} and \eqref{dubovic}, we find
\begin{align}
\left|\widehat{\Psi}\widehat{\Delta} \right|_{\infty}\le
\widehat{r}\left(\widehat\sigma+\widehat{\sigma}\left(\widehat\beta\right)
+\widehat{\sigma}(\beta)
-\widehat\sigma\left(\widehat\beta\right)\right)\le \widehat{r}\left(2\overline{\sigma}+  \widehat{r}\left|\widehat{\Delta}_{S_I}\right|_1+\left|\widehat{\Delta}_{S_I^c}\right|_1\right).\label{in1b}
\end{align}
The definition of the sensitivities yield, on $\mathcal{G}$, $
\left|\widehat{\Psi}\widehat{\Delta} \right|_{\infty}\le \widehat{r}\left(2\overline{\sigma}+\widehat{r}\frac{\left|\widehat{\Psi}\widehat{\Delta}
\right|_{\infty}}{\widehat\kappa_{\widehat{g},S(\beta)}}\right)$, hence
\begin{equation}\label{eq:ubc}
\left|\widehat{\Psi}\widehat{\Delta} \right|_{\infty}\le 2\widehat{r}\overline{\sigma}\gamma\left(\widehat{r}/\widehat\kappa_{\widehat{g},S(\beta)}\right).
\end{equation}
\eqref{eq:ubc} and the definition of the sensitivities yield the first upper bound. For the second, we use \eqref{stop1} and item (i) in Lemma \ref{sigmabasicbounds}. We now prove the second inequality. Take $\beta\in\mathcal{I}$ and  $S\subseteq[{d_X}]$.
Acting as in \eqref{in0}, on $\mathcal{G}$,
\begin{align*}
&\sum_{k\in S^c\cap S_Q}
\left|\mathbb{E}_n[X_k^2]^{1/2}\widehat\beta_k\right|+
\sum_{k\in S^c\cap S_Q}
\left|\mathbb{E}_n[X_k^2]^{1/2}\beta_k\right|\\
&\le\sum_{k\in S\cap S_Q}\left(
\left|\mathbb{E}_n[X_k^2]^{1/2}\beta_k\right|
-\left|\mathbb{E}_n[X_k^2]^{1/2}\widehat\beta_k\right|\right)+2\sum_{k\in S^c\cap S_Q}
\left|\mathbb{E}_n[X_k^2]^{1/2}\beta_k\right|+c\left(\widehat{\sigma}(\beta)-\widehat\sigma\left(\widehat\beta\right)\right)\\
&\le\left|\widehat{\Delta}_{S\cap S_Q}\right|_1+2\left|\bold{D}_{\bold{X}}^{-1}\beta_{S^c\cap S_Q}\right|_1+c\widehat{r}\left|\widehat{\Delta}_{S_I}\right|_1
+c\left|\widehat{\Delta}_{S_I^c}\right|_1.
\end{align*}
This yields
$|\widehat{\Delta}_{S^c\cap S_Q}|_1\le|\widehat{\Delta}_{S\cap S_Q}|_1+2|\bold{D}_{\bold{X}}^{-1}\beta_{S^c\cap S_Q}|_1+c\widehat{r}|\widehat{\Delta}_{S_I}|_1
+c|\widehat{\Delta}_{S_I^c}|_1$. We show the first inequality by considering two cases.\\
Case 1: $2|\bold{\widehat{
D}_X}^{-1}\beta_{S^c\cap S_Q}|_1\le|\widehat{\Delta}_{S\cap S_Q}|_1+c\widehat{r}|\widehat{\Delta}_{S_I}|_1
+c|\widehat{\Delta}_{S_I^c}|_1+|\widehat{\Delta}_{S_Q^c}|_1$, then $\widehat{\Delta}\in \widehat{\overline{K}}_{S}$.  From this and the definition of
$\widehat{\overline{\kappa}}_{\ell^q_{S_0},S}$, we get the upper bound corresponding to the first term in the minimum. To obtain the second term 
we use the first upper bound in item (ii) in Lemma \ref{sigmabasicbounds}.\\
Case 2: 
$2|\bold{\widehat{
D}_X}^{-1}\beta_{S^c\cap S_Q}|_1>|\widehat{\Delta}_{S\cap S_Q}|_1+c\widehat{r}|\widehat{\Delta}_{S_I}|_1
+c|\widehat{\Delta}_{S_I^c}|_1+|\widehat{\Delta}_{S_Q^c}|_1$, so
$|\widehat{\Delta}|_1=|\widehat{\Delta}_{S^c\cap S_Q}|_1+|\widehat{\Delta}_{S\cap S_Q}|_1+|\widehat{\Delta}_{S_Q^c}|_1\le 6|\bold{
\widehat{D}_X}^{-1}\beta_{S^c\cap S_Q}|_1$.\\
In conclusion, $|\widehat{\Delta}_{S_0}|_q$ is smaller than the maximum of the two bounds.\hfill $\square$\vspace{0.2cm}

\noindent{\bf Proof of Proposition \ref{prop:LBsensitivities}.}  
We use $|\Delta_{S\cap S_Q}|_1\le \min(s,|\widehat{S}\cap S_Q|)|\Delta_{\widehat{S}\cap S_Q}|_{\infty}$. The last constraint gives rise to 
the union of sets involving the linear constraint $|\Delta_{S\cap S_Q}|_1\le \min(s,|\widehat{S}\cap S_Q|)|\Delta_{S}|$, hence the second minimum. We conclude from the definition of the sensitivities, the cones $\widehat K_{\widehat{S}}$, 
and the fact that minimizing on a larger set yields lower bounds on the sensitivities. 
\hfill $\square$\vspace{0.2cm}

\noindent{\bf Proof of Theorem \ref{t1b}.} \eqref{t1bi} and \eqref{t1biii} follow from the second bounds in Proposition \ref{t1} and  
Lemma \ref{thrm:DLBsensitivities}. 
Part \eqref{t1bii} follows from  \eqref{t1bi} and \eqref{t1biii} with $\ell_k(\Delta)$ 
and the fact that the assumption on $|\beta_k|$ implies: $\widehat\beta_k \ne 0$ for $k\in S(\beta)$ (resp., $S_*$ as defined at the end of Section \ref{s54}). 
\hfill $\square$\vspace{0.2cm}

\noindent{\bf Proof of Theorem \ref{th:threshold}.} Fix $s$ and $\beta$ in
$\mathcal{I}_s$ and work on $\mathcal{G}\cap\mathcal{G}_{A1}$. Using lemmas \ref{thrm:DLBsensitivities}, \ref{thrm:DLBsensitivitiesO}, and \ref{sigmabounds} \eqref{estvar}, 
we obtain
$\widehat\omega_{k}(s) \le \omega_k(s)$. 
The following two cases can occur.\\ 
First, if $k\in S(\beta)^c$
(so that $\beta_k=0$) then, using the bound in \eqref{CC} for $\ell$ defined by $\ell(\Delta)=|\Delta_k|$ we obtain
$\mathbb{E}_n[X_k^2]^{1/2}|\widehat\beta_k|\le\widehat\omega_k(s)$, which implies
$\widehat{\beta}_{k}^{\widehat\omega}=0$.\\
Second, if $k\in S(\beta)$, then 
again by \eqref{CC},  
we get
$||\widehat\beta_k|-|\beta_k||\le|\widehat\beta_k-
\beta_k|\le\widehat\omega_k(s)/\sqrt{(1-\tau_n)\E[X_k^2]}\le\omega_k(s)/\sqrt{(1-\tau_n)\E[X_k^2]}$. Since $|\beta_k|>
2\omega_k(s)/\sqrt{(1-\tau_n)\E[X_k^2]}$ for $k\in S(\beta)$, we obtain $|\widehat\beta_k|>\omega_k(s)/\sqrt{(1-\tau_n)\E[X_k^2]}\ge
\widehat{\omega}_k(s)/\E_n[X_k^2]^{1/2}$, so that
$\widehat{\beta}_{k}^{\widehat\omega}=\widehat\beta_k$.\hfill$\square$\vspace{0.2cm}





\noindent{{\bf Proof of Theorem \ref{tCI}.}} 
The elements relative to assumptions and estimation of $\Lambda$ are in Appendix \ref{Lambdaest}, some of which are used below. 
Take $(\beta,\Lambda)\in\mathcal{I}_{\varPhi}$ and let
$$\mathcal{E}_{2,E}\triangleq\left\{\forall l\in[d_Z],|\mathbb{E}_n[G_l(\beta)E]|> \underline{r}_n^E\mathbb{E}_n[(G_l(\beta)E)^2]^{1/2},|\mathbb{E}_n[Z_lV(\beta)E]|> \underline{r}_n^E\mathbb{E}_n[(Z_lV(\beta)E)^2]^{1/2}\right\}.$$
We use  
$\max\left(\sqrt{1+\tau_n}-1,1-\sqrt{1-\tau_n}\right)=1-\sqrt{1-\tau_n}\le \tau_n$ and, for all $a\in\R^{d_Z}$, $b\in\R^{d_X}$,
 \begin{align*}
&
\mathbb{E}_n\left[\left(
a^{\top}
ZX
^{\top}b\right)^2\right]
=\left(D_{Z}^{-1}a\right)^{\top}
\mathbb{E}_n\left[D_{Z}ZX^{\top}bb^{\top}XZ^{\top}D_{Z}\right]D_{Z}^{-1}a
\le \left|D_{Z}^{-1}a\right|_1^2\left|D_{X}^{-1}b\right|_1^2(\overline{\rho}^{ZX})^2,\\
&\sqrt{n}\bold{D}_{\bold{F}(\widehat{\beta})\widehat{\Lambda}^{\top}}\left(\widehat{\varPhi\beta}-\varPhi\beta-\overline{V}(\beta)\right)=R+\sqrt{n}\bold{D}_{\bold{F}(\widehat{\beta})\widehat{\Lambda}^{\top}}\widehat{\Lambda}\mathbb{E}_n[G(\beta)],\\
&R=\sqrt{n}\bold{D}_{\bold{F}(\widehat{\beta})\widehat{\Lambda}^{\top}}\left(\varPhi-\widehat{\Lambda}\mathbb{E}_n[ZX^{\top}]\bold{D}_{\bold{X}}\right)\widehat\Delta-\sqrt{n}\bold{D}_{\bold{F}(\widehat{\beta})\widehat{\Lambda}^{\top}}\overline{V}(\beta)+\sqrt{n}\bold{D}_{\bold{F}(\widehat{\beta})\widehat{\Lambda}^{\top}}\widehat{\Lambda} \mathbb{E}_n\left[ZV(\beta)\right],\\
&T_{\varPhi}\triangleq\left|\sqrt{n}\bold{D}_{\bold{F}(\widehat{\beta})\widehat{\Lambda}^{\top}}\widehat{\Lambda} \mathbb{E}_n\left[G(\beta)\right]\right|_{\infty},\quad T_{\varPhi1}\triangleq\left|\sqrt{n}D_{\Lambda G(\beta)}\widehat{\Lambda}\mathbb{E}_n\left[G(\beta)\right]\right|_{\infty},\\
&T_{\varPhi0}\triangleq\left|\sqrt{n}D_{\Lambda G(\beta)}\Lambda \mathbb{E}_n\left[G(\beta)\right]\right|_{\infty},\quad
G_{\varPhi1}\triangleq\left|\sqrt{n}D_{\Lambda G(\beta)}\widehat{\Lambda}\mathbb{E}_n\left[F(\widehat\beta)E\right]\right|_{\infty},\\
&G_{\varPhi0}\triangleq\left|\sqrt{n}D_{\Lambda G(\beta)}\Lambda  \mathbb{E}_n\left[G(\beta)E\right]\right|_{\infty}. 
\end{align*}
We now work on $\mathcal{G}\cap\mathcal{E}$ 
of probability at least $1-\alpha^{S}_n-\alpha^{BC}_n$. 
For all $f\in[d_{\varPhi}]$, we have 
\begin{align*}
&
\left(D_{\Lambda G(\beta)}\right)_{f,f}\left|\mathbb{E}_n\left[(\widehat{\Lambda}_{f,\cdot}F(\widehat\beta))^2\right]^{1/2}-\mathbb{E}\left[(\Lambda_{f,\cdot}G(\beta))^2\right]^{1/2}\right|\\
&\le \mathbb{E}_n\left[\left(\left(D_{\Lambda G(\beta)}\right)_{f,f}\left(\widehat{\Lambda}_{f,\cdot}Z(X^{\top}(\widehat{\beta}-\beta)+V(\beta)+W(\beta))-\Lambda_{f,\cdot}ZW(\beta)\right)\right)^2\right]^{1/2}+\tau_n\\
&\le\mathbb{E}_n\left[\left(\left(D_{\Lambda G(\beta)}\right)_{f,f}(\widehat{\Lambda}_{f,\cdot}-\Lambda_{f,\cdot})ZW(\beta)\right)^2\right]^{1/2}\\
&\quad +\mathbb{E}_n\left[
\left(\left(D_{\Lambda G(\beta)}\right)_{f,f}(\widehat{\Lambda}_{f,\cdot}-\Lambda_{f,\cdot}+\Lambda_{f,\cdot})Z(X^{\top}(\beta-\widehat{\beta})+V(\beta))\right)^2\right]^{1/2}+\tau_n\le v^{D}_n.
\end{align*}
We have obtained $\left|\bold{D}_{\bold{F}(\widehat{\beta})\widehat{\Lambda}^{\top}}D_{\Lambda G(\beta)}^{-1}\right|_{\infty}\le 1/(1-v^{D}_n)$. 
On $\mathcal{G}\cap\mathcal{E}$, $|R|_{\infty}\le v_n^R$ and 
$\left|T_{\varPhi}-T_{\varPhi1}\right|\le T_{\varPhi1}v^D_n/(1-v^D_n)$, 
$\left|T_{\varPhi1}-T_{\varPhi0}\right|\le 
v_n^T$, 
 so 
$\left|T_{\varPhi}-T_{\varPhi0}\right|\le (T_{\varPhi0}+v_n^T)v^D_n/(1-v^D_n)+ v_n^T$.\\ 
Also, on $\mathcal{G}\cap\mathcal{E}\cap 
\left\{|\bold{E}|_{\infty}\le2\log\left(2n/\alpha_n\right)\right\}\cap\mathcal{E}_{2,E}^c\cap\mathcal{E}_{EZX^{\top}}^c$  of probability at least $1-\alpha^{S}_n-\alpha^{BC}_n-2\alpha_n-\alpha_n(EZX^{\top})$, $\left|G_{\varPhi}-G_{\varPhi1}\right|\le G_{\varPhi1}v^D_n/(1-v^D_n)$ and by convexity
\begin{align*}
\left|G_{\varPhi1}-G_{\varPhi0}\right|&\le
\left|\sqrt{n}D_{\Lambda G(\beta)}\left(\widehat{\Lambda}-\Lambda\right) \mathbb{E}_n\left[ZW(\beta)E\right]\right|_{\infty}\\
&\quad+\left|\sqrt{n}D_{\Lambda G(\beta)}\left(\widehat{\Lambda}-\Lambda+\Lambda\right)
\left( \mathbb{E}_n\left[ZX^{\top}\right](\beta-\widehat{\beta})+ \mathbb{E}_n\left[V(\beta)E\right]\right)\right|_{\infty}\le v_n^G,
\end{align*} 
hence $\left|G_{\varPhi}-G_{\varPhi0}\right|\le
(G_{\varPhi0}+v_n^G)v^D_n/(1-v^D_n)+v_n^G$. By Assumption \ref{assCI1a} \eqref{ezetaa},  we have $(\zeta_n-v_n^T)(1-v_n^D)/v_n^D-v_n^T\ge \zeta_n(1-2v_n^D)/(2v_n^D)$ and the same replacing $v_n^T$ by $v_n^G$.\\ 
Using \eqref{ezetaa} and $2\log\left(2d_\varPhi/\alpha_n\right)\ge q_{N_{\varPhi0}}(1-\alpha_n)$ where $N_{\varPhi0}\triangleq\left|D_{\Lambda G(\beta)}\Lambda E_{G(\beta)}\right|_{\infty}$ and $E_{G(\beta)}$
is a Gaussian vector of covariance $\mathbb{E}[G(\beta)G(\beta)^{\top}]$, by \eqref{GA} (which hold with obvious modifications), we get $\mathbb{P}\left(|T_\varPhi-T_{\varPhi 0}|>\zeta_n \right)\leq\zeta''_n$ and $\mathbb{P}\left(\mathbb{P}\left(|G_{\varPhi}-G_{\varPhi0}|>\zeta_n| \bold{Z}\widehat\Lambda^\top\right)>\zeta'_n\right)<\zeta'_n$. 
The second inequality uses the Markov inequality and the law of iterated expectations. We conclude like in Section \ref{conhomban} and the proof of Class 4 (see Section \ref{s:LBS}).\hfill $\square$\vspace{0.2cm}

\noindent{\bf Proof of Theorem \ref{th:nonvalid}.} It is given together with the proof of Theorem \ref{th:nonvalid1}.\hfill $\square$

\newpage
%
%
%

 \setcounter{equation}{0}  
 \setcounter{lmm}{0}
 \setcounter{crllr}{0}
 \setcounter{prpstn}{0}
 \setcounter{rmrk}{0}
 \setcounter{dfntn}{0}
 \setcounter{thrm}{0}
 \setcounter{ass}{0}
 \setcounter{algm}{0}
 \setcounter{section}{1}
  \setcounter{subsection}{0}
 \setcounter{page}{1}
 \setcounter{footnote}{0}
 \setcounter{figure}{0}
 \setcounter{table}{0}

\renewcommand{\thepage}{O-\arabic{page}}
\renewcommand{\theequation}{O.\arabic{equation}}
\renewcommand{\thelmm}{O.\arabic{lmm}}
\renewcommand{\thecrllr}{O.\arabic{crllr}}
\renewcommand{\thedfntn}{O.\arabic{dfntn}}
\renewcommand{\theprpstn}{O.\arabic{prpstn}}
\renewcommand{\thermrk}{O.\arabic{rmrk}}
\renewcommand{\thethrm}{O.\arabic{thrm}}
\renewcommand{\theass}{O.\arabic{ass}}
\renewcommand{\thealgm}{O.\arabic{algm}}
\renewcommand{\thesubsection}{O.\arabic{subsection}}
\renewcommand{\thetable}{O.\arabic{table}}
\section*{Online Appendix}
\subsection{Complements}
\subsubsection{Complements on Section \ref{cs3}}\label{s:LBS} 
The following propositions
establish lower bounds on $\widehat{\kappa}_{\ell^q,S}$ 
when $Z=X$, 
$d_Q=d_X$, $\mathcal{B}=\R^{d_X}$. 
Let $S\subseteq [d_X]$ and $c<1/\widehat{r}$. We have 
\begin{equation*}
\widehat{K}_S\subseteq C_S\triangleq\left\{\Delta\in\mathbb{R}^{d_X}:  (1-c\widehat{r})|\Delta_{S^c}|_1\le (1+c\widehat{r})|\Delta_S|_1\right\}. 
\end{equation*}
We define the
following generalizations of the restricted eigenvalue (RE) constants
$$
\kappa_{{\rm RE},S} \triangleq \inf_{\Delta\in\R^{d_X}\setminus\{0\}:\
\Delta\in C_S} \frac{|\Delta^{\top}\widehat{\Psi}\Delta|}{|\Delta_{S}|_2^2},
\quad \quad
\kappa'_{{\rm RE},S} \triangleq \inf_{\Delta\in\R^K\setminus\{0\}:\
\Delta\in C_S}
\frac{|S|\,|\Delta^{\top}\widehat{\Psi}\Delta|}{|\Delta_{S}|_1^2}.
$$
\begin{prpstn}\label{p3} For any $S\subseteq [d_X]$, we have
$$\kappa_{\ell^1,S}
\ge\frac{1-c\widehat{r}}{2}\kappa_{\ell^1_S,S} \ge \frac{(1-c\widehat{r})^2}{4|S|}\kappa'_{{\rm RE},S} \ge
\frac{(1-c\widehat{r})^2}{4|S|}\kappa_{{\rm RE},S}.$$
\end{prpstn}
\noindent{\bf Proof.} For $\Delta$ such that
$|\Delta_{S^c}|_1\le\frac{1+c\widehat{r}}{1-c\widehat{r}}|\Delta_{S}|_1$ we have
$|\Delta|_1\le\frac{2}{1-c\widehat{r}}|\Delta_{S}|_1$.  Thus, one obtains
$$
\frac{|\Delta^{\top}\widehat{\Psi}\Delta|}{|\Delta_S|_1^2} \le
\frac{|\Delta|_1|\widehat{\Psi}\Delta|_\infty}{|\Delta_S|_1^2} 
\le
\frac2{1-c\widehat{r}}\frac{|\widehat{\Psi}\Delta|_\infty}{|\Delta_S|_1}
\le
\frac4{(1-c\widehat{r})^2}\frac{|\widehat{\Psi}\Delta|_\infty}{|\Delta|_1}\,.
$$
Taking the infimum over $\Delta$'s proves the first two inequalities of the proposition. The second
inequality uses the fact that from H\"older's inequality $|\Delta|_1^2\le|S||\Delta_{S}|_2^2$.
\hfill $\square$\vspace{0.2cm}

We now obtain bounds for sensitivities $\kappa_{\ell^q,S}$ with $1< q\le
2$. For any $s\in[d_X]$, we consider the
restricted eigenvalue constant: $ \kappa_{{\rm RE}}(s) \triangleq
\min_{|S|\le s} \kappa_{{\rm RE},S}$.
\begin{prpstn}\label{p3_1} For any $s,m\in[d_X]$ such that $s+m\le d_X$ and $q\in(1,2]$, we have
$$\kappa_{\ell^q,S}
\ge C(q)\left(1+\frac{m}{s}\right)^{1/2-1/q} s^{-1/q}\kappa_{{\rm RE}}(s+m), \quad \forall \ S: \ |S|\le
s,
$$ where $C(q)=
\frac{1-c\widehat{r}}{2}
\left(1 + \frac{1+c\widehat{r}}{1-c\widehat{r}}\left(q-1\right)^{-1/q}\right)^{-1}$.
\end{prpstn}
\noindent{\bf Proof.} For $\Delta\in\R^{d_X}$ and a set $S\subseteq
[d_X]$, let $S_1$ 
be the subset of indices in
$[d_X]$ corresponding to the $m$ largest in absolute value
components of $\Delta$ outside of $S$. Define $S_+=S\cup S_1$. If
$|S|\le s$ we have $|S_+|\le s+m$. It is easy to see that the $k^{\text{th}}$ 
largest absolute value of elements of $\Delta_{S^c}$ satisfies
$|\Delta_{S^c}|_{(k)}\le |\Delta_{S^c}|_1/k$. Thus,
$$
|\Delta_{S^c_+}|_q^q =\sum_{j\in S^c_+}|\Delta_j|^q=\sum_{k\ge
s+1}|\Delta_{S^c}|_{(k)}^q \le |\Delta_{S^c}|_1^q\sum_{k\ge
s+1}\frac1{k^q} \le \frac{|\Delta_{S^c}|_1^q}{(q-1)s^{q-1}}\, .
$$
For $\Delta\in C_S$, this implies
$$
|\Delta_{S^c_+}|_q \le \frac{|\Delta_{S^c}|_1}{(q-1)^{1/q}s^{1-1/q}}
\le \frac{c_0|\Delta_{S}|_1}{(q-1)^{1/q}s^{1-1/q}} \le
\frac{c_0|\Delta_{S}|_q}{(q-1)^{1/q}}\, ,
$$
where $c_0=\frac{1+c\widehat{r}}{1-c\widehat{r}}$. Therefore, using that $|\Delta_{S}|_q\le |\Delta_{S_+}|_q$  we get, for $\Delta\in C_S$,
$$|\Delta|_q    \le |\Delta_{S_+}|_q + |\Delta_{S^c_+}|_q\le (1+c_0(q-1)^{-1/q})|\Delta_{S_+}|_q$$
so 
\begin{equation}\label{eq:prop_RE}
|\Delta|_q  \le
(1+c_0(q-1)^{-1/q})(s+m)^{1/q-1/2}|\Delta_{S_+}|_2.
\end{equation}
Using (\ref{eq:prop_RE}) and 
$|\Delta|_1\le\frac{2}{1-c\widehat{r}}|\Delta_{S}|_1
\le
\frac{2\sqrt{s}}{1-c\widehat{r}}|\Delta_{S}|_2\le \frac{2\sqrt{s}}{1-c\widehat{r}}|\Delta_{S_+}|_2
$ for $\Delta\in C_S$, we get
$$
\frac{|\Delta^{\top}\widehat{\Psi}\Delta|}{|\Delta_{S_+}|_2^2} \le 
\frac{|\Delta|_1|\widehat{\Psi}\Delta|_\infty}{|\Delta_{S_+}|_2^2}
\le \frac{2\sqrt{s}|\widehat{\Psi}\Delta|_\infty}{(1-c\widehat{r})|\Delta_{S_+}|_2}
\le
\left(1+\frac{m}{s}\right)^{1/q-1/2}
\frac{s^{1/q}|\widehat{\Psi}\Delta|_\infty}{C(q)|\Delta|_q}.$$
Using $|S_+|\le s+m$ we have proved the result.\hfill
$\square$\vspace{0.2cm}

We conclude this section by mentioning that, without endogeneity, the sensitivity $\widehat\kappa_{\ell^q,S}$ shares similarities with the characteristic introduced independently in \cite{YZ}. It differs in the definitions of $\widehat{K}_{S}$ and $\widehat{\Psi}$ 
and in that it does not involve scaling by $|S|^{1/q}$. Moreover \cite{DF} shows that previously introduced measures are computationally infeasible to verify but that the sensitivities that we introduce in this paper have desirable average-case perspective relative to NP-hardness in addition to them being weaker and more general than the others. 

To tighten the bounds in Table \ref{TabConst}, one can specify a small set $U\subseteq[d_X]$ and include the additional constraint $\mu_j=\eta_j\Delta_j$, $\forall j\in U$ in the LPs of Table \ref{TabConst}, where $\eta_j=\pm 1$ is the sign of $\Delta_j$. Since the signs are unknown, one replaces $\min_{k\in[d_X]}$ by $\min_{k\in[d_X],\eta_j=\pm 1\forall j\in U}$ in Table \ref{TabConst}. This augments the number of LPs by a factor of $2^{|U|}$. In our simulations we take $U=S_I^c$ to construct lower bounds based on a sparsity certificate. The design is such that $|U|=2$. If constructing lower bounds using $\widehat{S}$ of small cardinality, we use $U=\widehat{S}$. 

Other bounds can be derived from Proposition \ref{p4} and the following Proposition. Similar bounds can be obtained for the sensitivities  based on $\widehat{\overline{K}}_{S}$ (see \cite{GRT}).
\begin{prpstn}\label{p4b} Let $S\in[d_X]$, $c>0$, and $\widehat{r}\le1$. We have
\begin{align}
\hspace{-.5cm}
\widehat\kappa_{\ell^1,S}\ge \max&\left(\left(\frac{2}{\widehat\kappa_{\ell^1_{S\cap S_Q},S}}+\frac{1}{\widehat\kappa_{\ell^1_{S_Q^c},S}}+
\frac{c}{\widehat\kappa_{\widehat{g},S}}\right)^{-1}\notag\hspace{-.2cm},
\left(\gamma(c\widehat{r})\left(
\frac{2}{\widehat\kappa_{\ell^1_{S\cap S_Q},S}} +\frac{1}{\widehat\kappa_{\ell^1_{S_Q^c},S}}+\frac{c(1-\widehat{r})}{\widehat\kappa_{\ell^1_{S_I^c}, S}}  
\right)\right)^{-1}\hspace{-.2cm},\right.\\
&\left.\quad \left(\gamma(c)\left(\frac{2}{\widehat\kappa_{\ell^1_{S\cap S_Q},S}}+ \frac{1}{\widehat\kappa_{\ell^1_{S_Q^c},S}}\right)\right)^{-1}
\right);\label{bk1}\\
\hspace{-.5cm}
\widehat\kappa_{\widehat{g},S}\ge
\max&\left(\left(\gamma(c\widehat{r})\left(
\frac{2\widehat{r}}{\widehat\kappa_{\ell^1_{S\cap S_Q},S}}+
\frac{\widehat{r}}{\widehat\kappa_{\ell^1_{S_Q^c},S}}+
\frac{1-\widehat{r}}{\widehat\kappa_{\ell^1_{S_I^c}, S}}\right)\right)^{-1}\hspace{-.3cm},\left(\frac{\widehat{r}}{\widehat\kappa_{\ell^1_{S_I}, S}}+
\frac{1}{\widehat\kappa_{\ell^1_{S_I^c}, S}}\right)^{-1}\hspace{-.3cm},\widehat\kappa_{\ell^1,S}\right)\label{p46i2}
\end{align}
\end{prpstn}


In the cases where $S\subseteq\widehat{S}\subseteq S(\widehat{\beta})$ which we consider, we can use
$\widehat{c}_{\kappa}(S,S(\widehat{\beta}))\le \widehat{c}_{\kappa}(\widehat{S},S(\widehat{\beta}))$ and $\widehat{S}(S,S(\widehat{\beta}))\subseteq
(\widehat{S}\cap S_Q)\cup ((S_Q^c\cup S_I^c)\cap S(\widehat{\beta}))$, when $1\le c<\min(\widehat{r},1)^{-1}$, and $\widehat{S}(S,S(\widehat{\beta}))\subseteq
(\widehat{S}\cap S_Q)\cup (S_Q^c\cap S(\widehat{\beta}))$, when $c<1$. 

When  $\left|S\cap S_Q\right|\le s$, we have $\widehat{c}_{\kappa}(S,S(\widehat{\beta}))\le\widehat{c}_{\kappa}(s)\triangleq\min(\widehat{c}_{>,\kappa}(s),c_{<,\kappa}(s))$, where
$$\widehat{c}_{>,\kappa}(s)\triangleq 
\gamma(c\widehat{r})\left(2s+\left|S_Q^c\right|+c(1-\widehat{r})\left(\left|S_I^c\cap S_Q^c\right|+
\min\left(\left|S_I^c\cap S_Q\right|,s+\left|S_I^c\cap S_Q\cap S(\widehat{\beta})\right|\right)\right)\right)$$ and
$c_{<,\kappa}(s)\triangleq 
(2 s+|S_Q^c|)\gamma(c)$ 
and $\widehat{S}(S,S(\widehat{\beta}))\subseteq \overline{S}$. 
To compute a lower bound on $\widehat\kappa_{\ell^1}(\widehat{S})$, one can rely on \eqref{p44i} in Proposition \ref{p4} to obtain a lower bound on $\widehat\kappa_{\ell^\infty_{\widehat{S}(\widehat{S},S(\widehat{\beta}))},\widehat{S}}$ and multiply it by $\widehat{c}_{\kappa}(\widehat{S},S(\widehat{\beta}))^{-1}$. To compute a lower bound on $\widehat\kappa_{\ell^1}(s)$, one can use $\widehat{c}_{\kappa}(s)^{-1}\widehat\kappa_{\ell^\infty}(s)$.\\
\indent The lower bounds in Proposition \ref{prop:LBsensitivities} can be adapted to the sensitivities $\widehat{\overline{\kappa}}$ 
using sets $\widehat{\overline{B}}$ instead of $\widehat B$ involving the restrictions 
$-\mu\le\Delta\le \mu,-\nu1\le\widehat{\Psi}\Delta\le \nu1$  and 
$(1-2c\widehat{r})\sum_{k\in S_I}\mu_k+(1-2c)\sum_{k\in S_{I}^c}\mu_k\le 3\sum_{k\in \widehat{S}\cap S_Q}\mu_k+2\sum_{k\in S_Q^c}\mu_k$ for $\widehat{\overline{B}}(\widehat{S})$ and $(1-2c\widehat{r})\sum_{k\in S_I}\mu_k+(1-2c)\sum_{k\in S_{I}^c}\mu_k\le 3s\mu_j+2\sum_{k\in S_Q^c}\mu_k$ for $\widehat{\overline{B}}(j)$. 

We now present the adjustments for classes 1-3 with assumptions \ref{ass:Nemirovski} and \ref{ass:Nemirovskir} under 
independence between IVs and structural errors.
\begin{ass}\label{ass:Nemirovskir} Let 
${d_X},{d_Z}\ge3$, $M_{ZU^{\top}}\geq 0$. For all $(\beta,\mathbb{P})$ such that 
$\beta\in\mathcal{I}$, 
(N.\ref{Ni}) holds for $U(\beta)$ and $M_{U}$ and (N.\ref{Nii}) holds for $Z$ and $M_{Z}'$, we have
$Z$ and $U(\beta)$ are independent, and $$\E\left[\left|\left(\left(F_l(\beta)\right)^2/\left(\E\left[Z_l^2\right]\E\left[U(\beta)^2\right]\right)-1\right)_{l\in[d_Z]}\right|_{\infty}^2\right]\le M_{ZU^{\top}}.$$
\end{ass}
This is a condition of type (N.\ref{Ni}). Assumption \ref{ass:Nemirovskir} permits to work with $\widehat{r}=\underline{r}_n\sqrt{1+\tau_n}/(1-\tau_n)$, which is smaller than $\widehat{r}=\underline{r}_n\left|\bold{D}_{\bold{Z}}\bold{Z}^{\top}\right|_{\infty}$ as in the main text,  
and have 
$$\inf_{(\beta,\mathbb{P}):\ \beta\in\mathcal{I}} \mathbb{P}\left(\mathcal{G}\right)\geq\inf_{(\beta,\mathbb{P}):\ \beta\in\mathcal{I}} \mathbb{P}\left(\underline{\mathcal{G}}\right)-\alpha^C_n,$$
where $\alpha^C_n\triangleq \alpha_n(U)+ \alpha_n(Z)'+\alpha_n(ZU^{\top})$, because 
$$\underline{\mathcal{G}}\cap\left\{\max_{l\in[{d_Z}]}\mathbb{E}_n\left[\left(F_l(\beta)\right)^2\right]/\left(
\mathbb{E}_n[Z_{l}^2]\mathbb{E}_n[U(\beta)^2]\right)\le (1+\tau_n)/(1-\tau_n)^2\right\}\subseteq\mathcal{G}.$$ 
Combining Assumption \ref{ass:Nemirovskir} with any of classes 1-3 yields an upper bound on the coverage error, also denoted by $\alpha^B_n$, which is the one above 
plus $\alpha^C_n$.\\ 
We now present Class 4.\\
\noindent {\bf Class 4:} \emph{${d_Z}\ge3$, $M_{U}, M_Z,M_Z',q_2>0$, and a sequence $(B_n)_{n\in\N}$ such that $B_n\ge1$. For all $(\beta,\mathbb{P})$: $\beta\in\mathcal{I}$ and $q_1\in[2]$, 
\begin{enumerate}[\textup{(C4.}i\textup{)}]
\item\label{Bi} $\mathbb{E}\left[\left.U(\beta)^2\right|Z\right]=\sigma_{U(\beta)}^2$;
\item\label{Bii} (N.\ref{Ni}) holds for $U(\beta)$, $Z$ and $M_{U}$, $M_{Z}$; 
\item\label{Bv} $\left|\left(\max\left(\E\left[\left(\left(D_Z\right)_{l,l}F_l(\beta)/\sigma_{U(\beta)}\right)^{2+q_1}\right],
\E\left[\left(\left(D_Z\right)_{l,l}Z_lE\right)^{2+q_1}\right]
\right)\right)_{l=1}^{d_Z}\right|_{\infty}
\le B_n^{q_1}$;
\item\label{Bvi} $\max\left(\E\left[\left(\left|D_ZF(\beta)\right|_{\infty}/(B_n\sigma_{U(\beta)})\right)^{q_2}\right],
\E\left[\left(\left|D_ZZE\right|_{\infty}/B_n\right)^{q_2}\right]\right)\le 2$;
\end{enumerate}
where $E$ is standard normal independent of $Z$.\\ 
}
For the corresponding $\widehat{r}$, for all $n$, we have $\mathbb{P}\left(\mathcal{G}\right)\ge 1-\alpha-\alpha^B_n$, 
where $\alpha^B_n\triangleq 2\zeta'_n+(\zeta'_n)^2
+\varphi(d_Z,\tau_n)+\iota(d_Z,n)+\alpha_n(Z)+\alpha_n(U)$ and 
$(\zeta'_n)^2\triangleq\alpha_n+\iota(d_Z,n)+\alpha_n(Z)$. 
Also $\mathbb{P}(\widehat{r}\le r_n)\ge 1-\alpha^C_n$, where
$r_n\triangleq \left(
2\log(2d_Z/(\alpha-\zeta'_n-\varphi(d_Z,\tau_n)))+3\zeta_n
\right)/\sqrt{n}$ and 
$\alpha^C_n\triangleq \alpha_n(Z)
+\zeta'_n$.
\noindent {\bf Proof.} 
Let $\beta\in\mathcal{I}$. 
Define 
\begin{align*}
T\triangleq\left|\frac{\sqrt{n}}{\widehat{\sigma}(\beta)}\bold{D}_{\bold{Z}} \mathbb{E}_n[F(\beta)]\right|_{\infty},\ T_0\triangleq\left|\frac{\sqrt{n}}{\sigma_{U(\beta)}}D_{Z} \mathbb{E}_n[F(\beta)]\right|_{\infty},\ G_0\triangleq\left|\sqrt{n}D_Z\mathbb{E}_n[ZE]\right|_{\infty}.
\end{align*}
$T_0$, $G_0$, and $N_0\triangleq|D_Z\bold{E}_Z|_{\infty}$, where $\bold{E}_Z$ is a Gaussian vector of covariance $\mathbb{E}[ZZ^{\top}]$, have same covariance matrix, indeed
$$
\mathbb{E}\left[D_ZZZ^{\top}D_Z\frac{U(\beta)^2}{\sigma_{U(\beta)}^2}\right]=
\mathbb{E}\left[D_ZZZ^{\top}D_Z\mathbb{E}\left[\left.\frac{U(\beta)^2}{\sigma_{U(\beta)}^2}\right|Z\right]\right]
=\mathbb{E}[D_ZZZ^{\top}D_Z].$$
Let us show that, for all $\alpha\in(0,1)$, 
$
\left|\mathbb{P}\left(T\le q_{G|\bold{Z}}(\alpha)\right)-\alpha\right|
\le\alpha^B_n$. 
Using (C4.\ref{Bv}), (C4.\ref{Bvi}) and Proposition 2.1 in \cite{CCK2}, we obtain
\begin{equation}\label{GA}
\sup_{t\in\R}\max\left(
\left|\mathbb{P}\left(T_0\le t\right)-\mathbb{P}\left(N_0\le t\right)\right|
,\left|\mathbb{P}\left(G_0\le t\right)-\mathbb{P}\left(N_0\le t\right)\right|\right)
\le\iota(d_Z,n).
 \end{equation}
Indeed, by (C4.\ref{Bi}), the law of iterated expectations, and independence between $E$ and $Z$,  for all $\forall l\in[{d_Z}]$, $\E\left[\left(\left(D_Z\right)_{l,l}F_l(\beta)/\sigma_{U(\beta)}\right)^{2}\right]=\E\left[\left(\left(D_Z\right)_{l,l}Z_lE\right)^{2}\right]=1$, so Condition M1 in \cite{CCK2} holds. 
By the arguments in the proof of lemmas \ref{thrm:DLBsensitivities}, 
and Lemma 3.2 in \cite{CCK}, denoting by $q_{G_0|\bold{Z}}(\alpha)$ the $\alpha$ quantile of $G_0$ given $\bold{Z}$,
\begin{equation}\label{eiota}
\hspace{-.3cm}\min\left(\mathbb{P}\left(q_{G_0|\bold{Z}}(\alpha)\le q_{N_0}(\alpha+ \varphi(d_Z,\tau_n))\right),\mathbb{P}\left(q_{N_0}(\alpha)\le q_{G_0|\bold{Z}}(\alpha+ \varphi(d_Z,\tau_n))\right)\right)\ge 1-\alpha_n(Z).
\end{equation}
For all $\alpha\in(0,1)$, 
we have, by \eqref{GA}-\eqref{eiota},
\begin{align*}
&\alpha- \varphi(d_Z,\tau_n)-\alpha-\iota(d_Z,n)-\alpha_n(Z)\le\mathbb{P}\left(T_0\le q_{G_0|\bold{Z}}(\alpha)\right)-\alpha,\\ 
&\mathbb{P}\left(T_0\le q_{G_0|\bold{Z}}(\alpha)\right)-\alpha\le\alpha+ \varphi(d_Z,\tau_n)-\alpha+\iota(d_Z,n)+\alpha_n(Z)
\end{align*}
so 
\begin{equation}\label{Q1}
\left|\mathbb{P}\left(T_0\le q_{G_0|\bold{Z}}(\alpha)\right)-\alpha\right|
\le  \varphi(d_Z,\tau_n)+\iota(d_Z,n)+\alpha_n(Z).
\end{equation}
On $\mathcal{E}_Z'^c$, we have 
$|G-G_0|\le \left(1/\sqrt{1-\tau_n}-1\right)G_0,$
hence, by the Markov inequality, the law of iterated expectations, and the second bound in \eqref{GA}, 
\begin{align}
\mathbb{P}\left(\mathbb{P}\left(\left.|G-G_0|>\zeta_n\right|\bold{Z}\right)>\zeta'_n\right)&<\zeta'_n.\label{l331}
\end{align}
On $\mathcal{E}_Z'^c\cap\mathcal{E}_U^c$, we have $|T-T_0|\le \tau_n T_0/(1-\tau_n)$,
 hence, 
by the first bound in \eqref{GA},
\begin{equation}\label{l332}
\mathbb{P}\left(|T-T_0|>\zeta_n\right)\le(\zeta_n')^2+\alpha_n(U).
\end{equation}
Using Lemma 3.3 in \cite{CCK} and \eqref{l331} in the first display,  and \eqref{l332} and \eqref{Q1} in the second,
\begin{align*}
\hspace{-.3cm}\mathbb{P}\left(T-2\zeta_n\ge q_{G|\bold{Z}}(1-\alpha)\right)&<\mathbb{P}\left(T-\zeta_n\ge q_{G_0|\bold{Z}}(1-\alpha-\zeta'_n)\right)+\zeta'_n\\
&\le\alpha+2\zeta'_n+(\zeta'_n)^2+\varphi(d_Z,\tau_n)+\iota(d_Z,n)+\alpha_n(Z)+\alpha_n(U). 
\end{align*}
The bound $r_n$ on $\widehat{r}$ follows from Lemma 3.3 in 
\cite{CCK}, \eqref{eiota}, and for all $\alpha\in(0,1)$, $q_{N_0}(\alpha)\le2\log\left(2d_Z/\alpha\right)
$. 
\hfill$\square$\vspace{0.2cm}

Other classes can also be used. To account for dependent data one can use \cite{CSWX} and \cite{ZW} for results involving respectively self-normalization and the bootstrap (see also the references therein).

\subsubsection{Additional Material for Section \ref{sec:main}}\label{mainonline}
Hermite polynomials are orthonormal in $L^2(\mu)\triangleq\{f:\int_{\R}f(x)^2\exp(-x^2/2)dx<\infty\}$ equipped with $(f,g)_{L^2(\mu)}\triangleq\left( \int_{\R}f(x)g(x)\exp(-x^2/2)dx\right)/\sqrt{2\pi}$ defined for $f,g\in L^2(\mu)$, hence $D_X=D_Z=I$. Basic properties of these polynomials yield for $l\in[d_Z]$ and $k\in[d_Z]$
\begin{align*}
\Psi_{l,k}&=\E[h_{l-1}(\widetilde{Z})h_{k-1}(\widetilde{X}/\sqrt{\pi^2+\sigma^2})]=\E[h_{l-1}(\widetilde{Z})\E[h_{k-1}(\widetilde{X}/\sqrt{\pi^2+\sigma^2})|\widetilde{Z}]]\\
&=\E\left[h_{l-1}(\widetilde{Z}) \left(\frac{\pi}{\sqrt{\pi^2+\sigma^2}}  \right)^{k-1}      h_{l-1}(\widetilde{Z})\right]=\left(\frac{\pi}{\sqrt{\pi^2+\sigma^2}}  \right)^{k-1}\indic(l=k). 
\end{align*}

We now comment Assumption SV$(q)$ 
and the results (Theorem 1 and 2 and Corollary 1 and 2) in \cite{BCHN}.  
\cite{BCHN} introduce the sparse singular values for the matrix $\widehat{\Psi}$ and the bounds are for the sensitivities rather than for the population ones as in this paper to obtain deterministic bounds.  
 Assumption SV$(q,(\underline{\delta}_n)_{n\in\N},(\overline{\delta}_n)_{n\in\N},(l_n)_{n\in\N},\eta_0)$ \eqref{CSi} is from Theorem 1. Choosing $l_n=4\overline{\delta}_n^2u_{\kappa}^2/\underline{\delta}_n^2$ (resp. $l_n=\log(n)$ and $\underline{\delta}_n$ and $\overline{\delta}_n$ constant, in which case, for $n$ large enough, $l_n\ge4\overline{\delta}_n^2u_{\kappa}^2/\underline{\delta}_n^2$) corresponds to the choice made in Corollary 1 (resp. Corollary 2). As is apparent in Theorem 1 and Corollary 1, 
an unpleasant feature of Assumption SV$(q,(\underline{\delta}_n)_{n\in\N},(\overline{\delta}_n)_{n\in\N},(l_n)_{n\in\N},\eta_0)$ \eqref{CSi} is that $\underline{\delta}_n,\overline{\delta}_n$ depend on $l_n$ which depends on $\underline{\delta}_n,\overline{\delta}_n$. Theorem 2 gives yet another bound on the C-STIV (see Section \ref{scstiv}), however it is neither a rate nor a result that can be used to form a confidence set. 
It has a high-level assumption that a random counterpart of  \eqref{cnontrivial} (a more complicated one due to a maximum) is an event of probability converging to 1 while we use the function $\gamma$. In contrast with $\mathcal{G}\cap\mathcal{G}^\Psi$, the probability of the event depends on the dependence between $X$ and $Z$, and can be close to zero. Condition IC is used because C-STIV does not allow $c>1$. 

Using ideas in the proofs of   
Proposition \ref{p3_1} and Theorem 1 in \cite{BCHN}, and the notations introduced for Assumption SV$(q,(\underline{\delta}_n)_{n\in\N},(\overline{\delta}_n)_{n\in\N},(l_n)_{n\in\N},\eta_0)$, Proposition \ref{plb} below gives a lower bound on $\kappa_{\ell^1,S}$ which is tighter than (but in the same spirit as) the one derived from SV$(1,(\underline{\delta}_n)_{n\in\N},(\overline{\delta}_n)_{n\in\N},(l_n)_{n\in\N},\eta_0)$ \eqref{CSi}. From it we can easily derived results on rates of convergence as in Corollary \ref{crllrrate}.  
The lower bound depends explicitly on $S$ beyond its cardinality. We can bound the other sensitivities from it using \eqref{kbndIC}. It yields
\begin{equation}\label{lbsup}\kappa_{\ell^1,S}\ge\sup_{\eta\in(0,1)}\frac{\eta}{u_{\kappa}|S|}\max_{m\in\mathcal{M}(\eta)}\frac{1}{\sqrt{1+m/|S|}}\min_{\substack{S_1\subseteq S^c\\|S_1|=m}}\sigma_{\min}(S,S_1),
\end{equation}
where, for all $S,S_1\in[d_X]$ such that $S_1\subseteq S^c$, 
\begin{align*}
\sigma_{\min}(S,S_1)&=
\max_{\substack{R\subseteq[d_Z]\\|R|= |S|+ |S_1|}}
\sigma_{\min}(\Psi_{R,S\cup S_1}),\\ 
\sigma_{\max}(S,S_1)&=\min_{R(S,S_1)\in\mathcal{R}(S,S_1)}\max_{\substack{S_2\subseteq (S\cup S_1)^c\\
|S_2|= |S_1|}}
\sigma_{\max}(\Psi_{R(S,S_1),S_2}),
\end{align*}
$\mathcal{R}(S,S_1)$ is the collection of sets $R\subseteq[d_Z]$ which achieve the maximum in $\sigma_{\min}(S,S_1)$, and
for all $\eta\in(0,1)$, $\mathcal{M}(\eta)\subseteq [d_X-|S|]$ is the set of 
 $m$, 
such that, for all $S_1\subseteq S^c$ with $|S_1|=m$, there exists $R(S,S_1)\in\mathcal{R}(S,S_1)$ such that 
$$(u_{\kappa}-1)\sqrt{\frac{|S|}{m}}\max_{\substack{S_2\subseteq (S\cup S_1)^c\\
|S_2|= |S_1|}}
\sigma_{\max}(\Psi_{R(S,S_1),S_2})\le (1-\eta)\sigma_{\min}(S,S_1).$$
\begin{prpstn}\label{plb}
Under Condition IC,  we have
$$\kappa_{\ell^1,S}\ge\frac{1}{u_{\kappa}|S|}\max_{m\in[d_X-|S|]}\frac{1}{\sqrt{1+m/|S|}}\min_{\substack{S_1\subseteq S^c\\|S_1|=m}}\left(\sigma_{\min}(S,S_1)-(u_{\kappa}-1)\sqrt{\frac{|S|}{m}}\sigma_{\max}(S,S_1)\right).$$
\end{prpstn}
\noindent{\bf Proof.} 
Take $\Delta\in K_S$, $S_1$ the set of $m$ largest entries of $\Delta$ of index in $S^c$, $S_2$ the subsequent $m$ largest in $S^c$, and so forth, $U(S,S_1)$ (resp. $V(S,S_1)$) the matrix formed by stacking the left-singular vectors (resp. the right-singular vectors), and $\lambda\in\R^{|S|+m}$ such that $U(S,S_1)^{\top}\lambda=V(S,S_1)^{\top}\delta_{S\cup S_1}$, where $\delta_{S\cup S_1}$ is the restriction of $\Delta$ to $S\cup S_1$. Let $R(S,S_1)\in\mathcal{R}(S,S_1)$ which minimizes the expression in the definition of $\sigma_{\max}(S,S_1)$. 
By the inverse 
triangle inequality, 
\begin{align*}
\sigma_{\min}(S,S_1)|\Delta_{S\cup S_1}|_2^2&\le
\left|\lambda^{\top}\Psi_{R(S,S_1),S\cup S_{1}}\delta_{S\cup S_{1}}\right|\\
&\le \sum_{j\ge 2}
|\lambda^{\top}\Psi_{R(S,S_1),S_j}\delta_{S_j}|+|\lambda|_1\left|\Psi\Delta\right|_{\infty}\\
&\le \sum_{j\ge 2}
|\lambda^{\top}\Psi_{R(S,S_1),S_j}\delta_{S_j}|+\sqrt{|S|+m}|\Delta_{S\cup S_1}|_2
\left|\Psi\Delta\right|_{\infty}.
\end{align*}
For $j\ge 2$, we have 
\begin{align*}
|\lambda^{\top}\Psi_{R(S,S_1),S_j}\delta_{S_j}|
&\le |\lambda^{\top}\Psi_{R(S,S_1),S_j}|_2|\Delta_{S_j}|_2\\
&\le \frac{1}{\sqrt{m}} |\lambda^{\top}\Psi_{R(S,S_1),S_j}|_2|\Delta_{S_{j-1}}|_1\le \frac{1}{\sqrt{m}}\sigma_{\max}(S,S_1)|\Delta_{S_{j-1}}|_1|\Delta_{S\cup S_1}|_2,
\end{align*}
so, using $\Delta\in K_S$ in the last inequality, 
$$\sum_{j\ge 2}|\lambda^{\top}\Psi_{R(S,S_1),S_j}\delta_{S_j}|\le \frac{u_{\kappa}-1}{\sqrt{m}}\sigma_{\max}(S,S_1)|\Delta_{S}|_1|\Delta_{S\cup S_1}|_2.$$
This yields 
$$
\sigma_{\min}(S,S_1)|\Delta_{S\cup S_1}|_2-\frac{u_{\kappa}-1}{\sqrt{m}}\sigma_{\max}(S,S_1)|\Delta_{S}|_1\le \sqrt{|S|+m}
\left|\Psi\Delta\right|_{\infty}$$
and by $|\Delta_{S}|_1/\sqrt{|S|}\le |\Delta_{S}|_2\le |\Delta_{S\cup S_1}|_2$ and $|\Delta|_1\le u_{\kappa}|\Delta_S|_1$, 
$$\frac{1}{u_{\kappa}|S|\sqrt{1+m/|S|}}\left(\sigma_{\min}(S,S_1)-(u_{\kappa}-1)\sqrt{\frac{|S|}{m}}\sigma_{\max}(S,S_1)\right)\le \frac{\left|\Psi\Delta\right|_{\infty}}{|\Delta|_1}.$$
We obtain the result because in the above expression $S_1$ depends on $\Delta$ and $m$ is arbitrary.
\hfill$\square$

\subsubsection{Bounds for $\widehat\sigma(\beta)$, $\widehat\sigma(\widehat\beta)$ and $\widehat\sigma$ and Nonparametric IV}
\noindent We use Lemma \ref{sigmabasicbounds} to prove Proposition \ref{t1} and Lemma \ref{sigmabounds} to prove Theorem \ref{t1b}.
\begin{lmm}\label{sigmabasicbounds} For all $(\beta,\mathbb{P})$ such that $\beta\in\mathcal{I}$, any STIV estimator and $c>0$, we have, on $\mathcal{G}$,
\begin{enumerate}[\textup{(}i\textup{)}]
\item for sparse vectors 
$$\widehat{\sigma}+\widehat{\sigma}(\beta)\le 2\widehat{\sigma}(\beta)\gamma\left(\frac{\widehat{r}}{c\widehat\kappa_{\ell^1_{S(\beta)\cap S_Q},S(\beta)}}  \right),$$
\item for arbitrary vectors 
$$
\widehat{\sigma}+\widehat{\sigma}(\beta)\le 2\min_{S\subseteq[{d_X}]}\max\left(\widehat{\sigma}(\beta)\gamma\left(\frac{\widehat{r}}{c\widehat{\overline{\kappa}}_{\widehat{h},S}}\right),\widehat{\sigma}(\beta)+\frac{3}{2c}\left|\bold{D}_{\bold{X}}^{-1}\beta
_{S^c\cap S_Q}\right|_1\right).
$$
\end{enumerate}
\end{lmm}
\noindent{\bf Proof.} 
By
\eqref{eq:main} and the definition of
$\widehat\kappa_{\ell^1_{S(\beta)\cap S_Q},S(\beta)}$,
\begin{align}
c\widehat{\sigma}&\le|\widehat{\Delta}_{S(\beta)\cap S_Q}|_1+c
\widehat{\sigma}(\beta)\le |\widehat{\Psi}\widehat{\Delta}|_{\infty}/\widehat\kappa_{\ell^1_{S(\beta)\cap S_Q},S(\beta)}
+c\widehat{\sigma}(\beta),\label{stopsigma}
\end{align}
and, by adding $c\widehat{\sigma}(\beta)$ to both sides and \eqref{stop1}, we obtain the first term in the minimum.\\
To deal with approximately sparse vectors, we use that in Case 1 in Proposition \ref{t1}
\begin{align}
\widehat{\sigma}&\le\frac{1}{c}\left(\left|\bold{D}_{\bold{X}}^{-1}\beta_{S_Q}\right|_1
-\left|\bold{D}_{\bold{X}}^{-1}\widehat{\beta}_{S_Q}\right|_1\right)+\widehat{\sigma}(\beta)\notag\\
&\le\frac{1}{c}\min\left(\left|\widehat{\Delta}_{S_Q}\right|_1,\left|\widehat{\Delta}_{S\cap S_Q}\right|_1+\left|\bold{D}_{\bold{X}}^{-1}\beta
_{S^c\cap S_Q}\right|_1\right)+
\widehat{\sigma}(\beta)\label{refdessous}\\
&\le\frac{1}{c}\min\left(\left|\widehat{\Delta}_{S_Q}\right|_1,\frac12\left(3\left|\widehat{\Delta}_{S\cap S_Q}\right|_1+c\widehat{r}\left|\widehat{\Delta}_{S_I}\right|_1+c\left|\widehat{\Delta}_{S_I^c}\right|_1+\left|\widehat{\Delta}_{S_Q^c}\right|_1\right)\right)+
\widehat{\sigma}(\beta)\label{stopsigman2}\\
&\le\frac{\left|\widehat{\Psi}\widehat{\Delta} \right|_{\infty} }{c\widehat{\overline{\kappa}}_{\widehat{h},S}}+
\widehat{\sigma}(\beta),\label{stopsigman}
\end{align}
which, with \eqref{stop1} yields the first upper bound. 
We obtain the second one using the inequality in Case 2 in Proposition \ref{t1} and \eqref{refdessous}.\hfill $\square$
\begin{lmm}\label{sigmabounds}
We have $\widehat{\sigma}(\widehat\beta)\le\widehat\sigma$ and, under the assumptions of Theorem \ref{t1b}, on $\mathcal{G}\cap\mathcal{G}_{A1}$,  
\begin{enumerate}[\textup{(}i\textup{)}]
\item\label{estvar} for sparse vectors 
\begin{align}
&
\sqrt{1-\tau_n}\sigma_{U(\beta)}\left(1-\frac{2r_n\Gamma_\kappa(S(\beta))}{\kappa_{g,S(\beta)}}\right)
\le\widehat{\sigma}(\widehat\beta)\le\widehat\sigma\le\sqrt{1+\tau_n}\sigma_{U(\beta)}
\left(1+\frac{2r_n\Gamma_\kappa(S(\beta))}{c\kappa_{\ell^1_{S(\beta)\cap S_Q},S(\beta)}}\right)
,\notag
\end{align}
\item\label{estvar2} for arbitrary vectors 
{\footnotesize\begin{align}
&\hspace{-.4cm}
\sqrt{1-\tau_n}\left(\sigma_{U(\beta)}-
\min_{S\subseteq[{d_X}]}\max\left(\sigma_{U(\beta)}\frac{2r_n\Gamma_{\overline{\kappa}}(S)}{\overline{\kappa}_{g,S}}
,\frac{2}{c1_{n}}\left| D_X^{-1}\beta_{S^c\cap S_Q}\right|_1\right)\right)\le\widehat{\sigma}(\widehat\beta)\notag\\
&\widehat\sigma\le\sqrt{1+\tau_n}\left(\sigma_{U(\beta)}+\frac{1}{c}\min_{S\subseteq[{d_X}]}\max\left(
\sigma_{U(\beta)}\frac{2r_n\Gamma_{\overline{\kappa}}(S)}{\overline{\kappa}_{h,S}}
,3\left| D_X^{-1}\beta_{S^c\cap S_Q}\right|_1\right)\right).\notag
\end{align}}
\end{enumerate}
\end{lmm}
\noindent {\bf Proof.} 
The first upper bound in \eqref{estvar} comes from \eqref{dubovic} and Lemma \ref{thrm:DLBsensitivities}. The last one comes from 
the first inequality in \eqref{stopsigma} and Lemma \ref{thrm:DLBsensitivities}.\\ 
Similarly, the first upper bound in \eqref{estvar2} comes from \eqref{dubovic}, Lemma \ref{thrm:DLBsensitivities}, and the fact that in Case 2 (see the proof of Proposition \ref{t1})  $\widehat{g}(\widehat\Delta)\le 2|(\bold{\widehat{
D}_X}^{-1}\beta)_{S^c\cap S_Q}|_1/c$. The last one comes from 
the first inequality in \eqref{stopsigma}, \eqref{refdessous}-\eqref{stopsigman2}, and Lemma \ref{thrm:DLBsensitivities}.\\ 
We now present the following complement to Proposition \ref{ThoracleSE}. We use $\Psi_X=D_X\mathbb{E}[XX^{\top}]D_X$. 
\begin{prpstn}\label{ThoracleSEO} 
For all $(\beta,\mathbb{P})$ such that $\beta\in\mathcal{I}$, assuming as well
 $\E_n[v_{g,d_X}^2]\le \widehat{v}_g^2$ on $\mathcal{G}_{A1}$ and all solution $\left(\widehat{\beta},\widehat\sigma\right)$ of \eqref{ESTIV}, the following hold on $\mathcal{G}\cap\mathcal{G}_{A1}$  
(on $\mathcal{G}\cap\mathcal{G}_{A1}\cap \mathcal{E}_X'^c$ for the second inequality of \eqref{ThoracleSEii})
\begin{enumerate}[\textup{(}i\textup{)}]  
\item\label{ThoracleSEOi}
For a sparse matrix $\beta$, for all $\ell\in\mathcal{L}$, we have
\begin{align*}
&\hspace{-.5cm}
\mathbb{E}_n\left[\left(X^{\top}\left(\widehat{\beta}_{\cdot,g}-\beta_{\cdot,g}\right)\right)^2\right]\le\frac{2r_n}{1_n\sqrt{\kappa_{\ell^1,S(\beta)}}}\left(\sum_{g=1}^{d_G}\sigma_{W_g(\beta)}+\left(r_n+2\right)v_{g,d_X}\right)\Gamma_{\kappa}(S(\beta));
\end{align*}
\item\label{ThoracleSEOii}
$
\mathbb{E}_n\left[\left(X^{\top}\left(\widehat{\beta}_{\cdot,g}-\beta_{\cdot,g}\right)\right)^2\right]^{1/2}\le 
\left|D_X^{-1}\left(\widehat{\beta}_{\cdot,g}-\beta_{\cdot,g}\right)\right|_{1}\sqrt{|\Psi_X|_{\infty}+\tau_n}$. 
\end{enumerate}
In all cases, we have
\begin{equation}\label{esigma}
|\widehat{\sigma}_g(\widehat\beta)-\sigma_{W_g(\beta)}|\le\mathbb{E}_n\left[\left(X^{\top}\left(\widehat{\beta}_{\cdot,g}-\beta_{\cdot,g}\right)\right)^2\right]^{1/2}+
\sigma_{W_g(\beta)}\tau_n
+\sqrt{1+\tau_n}v_{g,d_X}.
\end{equation}
\end{prpstn}
\subsubsection{C-STIV}\label{scstiv}
The C-STIV estimator applies the self-tuning for every moment like the NV-STIV estimator, but estimates $\theta$ and $\beta$ simultaneously using all $d_Z$ moments.  
This results in a conic program with multiple conic constraints. The model, restrictions, and sets of $s,\widetilde{s}$-identifiable parameters considered in this section are defined in 
Section \ref{sec:endiv}. We simply use $\mathcal{I}$ for $\mathcal{I}_{[d_Q],[d_Z-d_\perp]}$. 
We maintain either of class 1-3 replacing $F(\beta)$ by $T(\beta,\theta)\triangleq F(\beta)-\theta$. 
\begin{dfntn} For $c\in(0,1)$, a 
C-STIV estimator $(\widehat{\beta}, \widehat{\theta}, \widehat{\sigma})$ is
any solution of 
$$\min_{\left(b,t\right)\in\widehat{\mathcal{I}}_C( \underline{r}_n,\sigma),\sigma\ge0}
\left|\bold{D}_{\bold{X}}^{-1}b_{S_Q}\right|_1+\left|\bold{D}_{\bold{Z}}t_{ S_{\perp}^c}\right|_1+c\sigma,$$
where, for $r,\sigma>0$, 
\begin{align*}
&\widehat{\mathcal{I}}_C(r,\sigma)\triangleq\left\{\left(b,t\right)\in\mathcal{B}\times\Theta:\ \left|\bold{D}_{\bold{Z}}\left(\E_n[ZU(b)]-t\right)\right|_{\infty}\le r\sigma,\widehat{\Sigma}\left(b,t\right) \le\sigma\right\},\\
&\widehat{\Sigma}(b,t)\triangleq\max_{l\in[{d_Z}]}\widehat{\sigma}_l(b,t),\quad \widehat{\sigma}_l(b,t)^2\triangleq (\bold{D}_{\bold{Z}})_{l,l}^2\mathbb{E}_n[(Z_{l}U(b)-t_l)^2].\end{align*}
\end{dfntn}
When $d_{\perp}=d_Z$, C-STIV is an alternative to STIV, 
and its analysis is similar. The main difference is that there are $d_Z$ second-order conic constraints, whereas STIV has one. 
This makes C-STIV harder to compute, especially when $d_Z$ is large. However, we use a smoothing approach to approximate it which make its computation faster, as explained in Section \ref{sFISTA}. We use the same procedure to compute the BC-STIV estimator defined in \eqref{Mhat}, which is a C-STIV estimator for a system of $d_{\varPhi}$ equations, and involves $d_{\varPhi}d_X$ second-order cones. 
The C-STIV in the formulation of this paper imposes $c<1$. 

To obtain rates of convergence, the class is restricted similarly to Assumption \ref{ass:Nemirovski}, replacing $\mathbb{P}\left(\left|D_Z\bold{Z}^{\top}\right|_{\infty}> B_Z\right)\le \alpha_n$ by $\mathbb{P}\left(\widehat{\rho}^{ZX}>\rho^{ZX}_n\right)\le \alpha_n$, where $\rho^{ZX}_n$ depends on $n$ via ${d_Z}$ and $d_X$, 
and (N.\ref{Ni}) holds for $U(\beta)$ and $M_{U}$ by (N.\ref{Nii}) for $T(\beta,\theta)$ and $M_{T}'$. 
The event $\mathcal{G}_{A1}$ is modified accordingly. For simplicity, we continue to refer to these as 
Assumption \ref{ass:Nemirovski} and $\mathcal{G}_{A1}$ 
and use 
\begin{align*}
\alpha^C_n&=\alpha_n+\alpha_n(T)' +\alpha_n(Z)',\ 
\alpha^{A1}_n=\alpha_n(X)'+\alpha_n(ZX^{\top})
+\alpha^B_n+\alpha^C_n.
\end{align*} 
By Assumption \ref{ass:Nemirovski}, for all 
$n\in\N$, $\mathbb{P}(\mathcal{G}_{A1})\ge 1-\alpha^{A1}_n\to1$. 
Table \ref{TabCorresp} provides the C-STIV analogues of STIV objects. The cones become $\widehat{K}_{S,\widetilde{S}}$ and $\widehat{\overline{K}}_{S,\widetilde{S}}$, where $S\subseteq[{d_X}]$ and $\widetilde{S}\subseteq[{d_Z}]$. 
The population sensitivities $\kappa$ and $\overline{\kappa}$ are obtained by replacing $\widehat{\Psi}$, $\widehat{K}_{S,\widetilde{S}}$, $\widehat{\overline{K}}_{S,\widetilde{S}}$ by $\Psi$, $K_{S,\widetilde{S}}$, and $\overline{K}_{S,\widetilde{S}}$ in the definition of the $\widehat\kappa$ and $\widehat{\overline{\kappa}}$. 
The sensitivities, their population counterparts, and lower bounds depend either on two sets $S$ and $\widetilde{S}$ or on two sparsity certificates $s$ and $\widetilde{s}$. Computable lower bounds on the sensitivities are obtained by LP using the sets $\widehat{B}(k,l)$ and $\widehat{B}(S,\widetilde{S})$, and bounds on the population sensitivities are obtained identically but replacing $\widehat\rho^{ZX}$ by $\rho_n^{ZX}$. 
\begin{prpstn}\label{prop:cstiv} 
On the event $\mathcal{G}_{A1}$, we have, for all $c>0$, 
\begin{align}
&(D_Z)_{l,l}\sigma_{T_l(\beta,\theta)}1_n\le\widehat{\sigma}_l(b,t)\le(D_Z)_{l,l}\sigma_{T_l(\beta,\theta)}/1_n,\ \forall l\in[d_Z];\notag\\
&\forall \left(b,t\right)\in\R^{{d_X}+{d_Z}},\ell\in\mathcal{L},\sqrt{1-\tau_n}\ell\left( D_X^{-1}b,D_Zt\right)\le \ell\left(\bold{D}_{\bold{X}}^{-1}b,\bold{D}_{\bold{Z}}t\right)\le \ell\left( D_X^{-1}b,D_Zt\right)/\sqrt{1-\tau_n};\notag\\
&\forall S\subseteq[{d_X}],\widetilde{S}\subseteq[{d_Z}],l\in[{d_Z}],\ 
\widehat\kappa_{\ell,S,\widetilde{S}}\ge\kappa_{\ell,S,\widetilde{S}}\left(1-\frac{\tau_n}{\kappa_{\ell^1_{[{d_X}],\varnothing},S,\widetilde{S}}}\right)1_n;\notag\\
&\hspace{4.5cm}\widehat{\overline{\kappa}}_{\ell,S,\widetilde{S}}\ge\overline{\kappa}_{\ell,S,\widetilde{S}}\left(1-\frac{\tau_n}{\overline{\kappa}_{\ell^1_{[{d_X}],\varnothing},S,\widetilde{S}}}\right)1_n;\notag\\
&\text{if}\ |S\cap S_Q|\leq s,|\widetilde{S}|\le\widetilde{s},\ \forall \ell\in\mathcal{L},\ \widehat\kappa_{\ell}(s,\widetilde{s})\ge\kappa_{\ell}(s,\widetilde{s})
\triangleq\kappa_{\ell}^0(s,\widetilde{s})\min_{\substack{S:|S\cap S_Q|\leq s\\ \widetilde{S}:|\widetilde{S}|\leq \widetilde{s}}}\left(1-\frac{\tau_n}{\kappa_{\ell^1_{[{d_X}],\varnothing},S,\widetilde{S}}}\right)1_n.\notag
\end{align}
 The results for C-STIV are the same as those for STIV using Table \ref{TabCorresp} for correspondance. 
 \end{prpstn} 
 
\noindent {\bf Proof.} 
Take $(\beta,\theta)\in\mathcal{I}$.
Set $\widehat{\Delta}\triangleq\bold{D}_{\bold{X}}^{-1}(\widehat{\beta} -\beta)$ and
$\widetilde{\Delta}\triangleq\bold{D}_{\bold{Z}}(\widehat{\theta} -\theta)$. We now work on $\underline{\mathcal{G}}$. 
Clearly $\left(\beta,\theta\right)$ belongs to $\widehat{\mathcal{I}}_C( \underline{r}_n,\widehat{\Sigma}\left(\beta,\theta\right))$ and by the arguments in the proof of Proposition \ref{t1}
\begin{align}
\left|\widehat{\Psi}\widehat{\Delta}+\widetilde{\Delta}\right|_{\infty} &\le
 \underline{r}_n\left(\widehat{\sigma}+\widehat{\Sigma}\left(\beta,\theta\right)\right)\label{stop2C}\\
\left|\widehat{\Delta}_{S(\beta)^c\cap S_Q}\right|_1+\left|\widetilde{\Delta}_{S\left(\theta\right)^c}\right|_1&\le
\left|\widehat{\Delta}_{S(\beta)\cap S_Q}\right|_1+\left|\widetilde{\Delta}_{S\left(\theta\right)}\right|_1+c\left(\widehat{\Sigma}\left(\beta,\theta\right)-\widehat{\Sigma}\left(\widehat\beta,\widehat{\theta}\right)\right)\nonumber.
\end{align}
Each function 
$\widehat{\sigma}_l$ 
is convex and
 $$w_{l*}\triangleq- \left(
\begin{array}{c}
\overline{w}_{l} \\
\widetilde{w}_{l}
\end{array}
\right)
\indic{\left\{\mathbb{E}_n\left[T_l(\beta,\theta)^2\right]\ne0\right\}}\in\partial \widehat{\sigma}_l\left(\beta,\theta\right),$$
where 
$$\overline{w}_{l}\triangleq\frac{\mathbb{E}_n\left[X Z_{l}T_l(\beta,\theta)\right]}{\mathbb{E}_n\left[Z_{l}^2\right]^{1/2}\mathbb{E}_n\left[T_{l}(\beta,\theta)^2\right]^{1/2}},\quad   
\widetilde{w}_l\triangleq
\left(\begin{array}{c}
 0\\ 
\frac{\mathbb{E}_n\left[T_l(\beta,\theta)\right]}{\mathbb{E}_n\left[Z_{l}^2\right]^{1/2}\mathbb{E}_n\left[T_l(\beta,\theta)^2\right]^{1/2}}\\
 0  
\end{array}\right).$$

By the Cauchy-Schwarz inequality, for all $k\in [d_X]$, 
$\left(\bold{D}_{\bold{X}}\right)_{k,k}|(\overline{w}_{l})_k|\le \widehat{\rho}^{ZX}$. Taking $w_*=(\overline{w}^{\top},\widetilde{w}^{\top})^{\top}$ as one of the $w_{l*}$ for which $\widehat{\sigma}_l\left(\beta,\theta\right)=\widehat{\Sigma}\left(\beta,\theta\right)$ yields an element of $\partial{\widehat{\Sigma}}\left(\beta,\theta\right)$. 
By definition of the subdifferential $\partial{\widehat{\Sigma}}\left(\beta,\theta\right)$, we have
\begin{align}
\widehat{\Sigma}\left({\beta},\theta\right)- \widehat{\Sigma}\left(\widehat{\beta},\widehat{\theta}\right) &\le w_*^{\top}
\left(\begin{array}{ccc}
\beta-\widehat{\beta}\\
\theta-\widehat{\theta}    
\end{array}\right)\le \left|\bold{D}_{\bold{X}}\overline{w}\right|_\infty \left|\widehat{\Delta}\right|_1+\left|\bold{D}_{\bold{Z}}^{-1}\widetilde{w}\right|_\infty \left|\widetilde{\Delta}_{S_{\perp}^c}\right|_1\notag\\
&\le \widehat{\rho}^{ZX}\left|\widehat{\Delta}\right|_1+ \underline{r}_n\left|\widetilde{\Delta}_{S_{\perp}^c}\right|_1\label{dubovicC}.
\end{align}
As a result, we have $(\widehat{\Delta},\widetilde{\Delta})\in \widehat{K}_{S(\beta),S\left(\theta\right)}$. Using \eqref{stop2C} and \eqref{dubovicC}, we find
\begin{align}
\left|\widehat{\Psi}\widehat{\Delta}+\widetilde{\Delta}\right|_{\infty} &\le  \underline{r}_n\left(2\overline{\sigma}+\widehat{\rho}^{ZX}\left|\widehat{\Delta}\right|_1+ \underline{r}_n\left|\widetilde{\Delta}_{S_{\perp}^c}\right|_1\right).\label{in1b}
\end{align}
Using the definition of the sensitivities, we obtain
\begin{align*}
&\left|\widehat{\Psi}\widehat{\Delta}+\widetilde{\Delta}\right|_{\infty}\le  \underline{r}_n\left(2\overline{\sigma}+ \frac{\left|\widehat{\Psi}\widehat{\Delta}+\widetilde{\Delta}\right|_{\infty}}{\widehat\kappa_{\widehat{g},S(\beta),S\left(\theta\right)}}\right)\le 2 \underline{r}_n\overline{\sigma}\gamma\left(\frac{\underline{r}_n}{\widehat\kappa_{\widehat{g},S(\beta),S\left(\theta\right)}}\right),\\
&
c\widehat{\sigma}\le|\widehat{\Delta}_{S(\beta)\cap S_Q}|_1+\left|\widetilde{\Delta}_{S\left(\theta\right)}\right|_1+c
\widehat{\Sigma}\left(\beta,\theta\right)\le
\frac{\left|\widehat{\Psi}\widehat{\Delta}+\widetilde{\Delta}\right|_{\infty}}{\widehat\kappa_{\ell^1_{S(\beta)\cap S_Q,S\left(\theta\right)},S(\beta),S\left(\theta\right)}}
+c\widehat{\Sigma}\left(\beta,\theta\right).
\end{align*}
For nonsparse vectors, $S\subseteq[{d_X}]$ and $\widetilde{S}\subseteq[{d_Z}]$,  we obtain
\begin{align*}
\left|\widehat{\Delta}_{S^c\cap S_Q}\right|_1+\left|\widetilde{\Delta}_{\widetilde{S}^c}\right|_1&\le\left|\widehat{\Delta}_{S\cap S_Q}\right|_1+\left|\widetilde{\Delta}_{\widetilde{S}}\right|_1+c\left(\widehat{\rho}^{ZX}\left|\widehat{\Delta}\right|_1+ \underline{r}_n\left|\widetilde{\Delta}_{S_{\perp}^c}\right|_1\right)\\
&\quad +2\left|\bold{D}_{\bold{X}}^{-1}\beta_{S^c\cap S_Q}\right|_1+2\left|\bold{D}_{\bold{Z}}\theta_{\widetilde{S}^c}\right|_1.
\end{align*}
We again consider two cases.\\ 
First, if 
$2|\bold{D}_{\bold{X}}^{-1}\beta_{S^c\cap S_Q}|_1+2|\bold{D}_{\bold{Z}}\theta_{\widetilde{S}^c}|_1\le
|\widehat{\Delta}_{S\cap S_Q}|_1+|\widetilde{\Delta}_{\widetilde{S}}|_1+c(\widehat{\rho}^{ZX}|\widehat{\Delta}|_1+ \underline{r}_n|\widetilde{\Delta}_{S_{\perp}^c}|_1)+|\widehat{\Delta}_{S_Q^c}|_1$, then $(\widehat{\Delta},\widetilde{\Delta})\in\widehat{\overline{K}}_{S,\widetilde{S}}$. Also, we have
\begin{align}
\widehat{\sigma}&\le\frac{1}{c}\left(\left|\bold{D}_{\bold{X}}^{-1}\beta_{S_Q}\right|_1
-\left|\bold{D}_{\bold{X}}^{-1}\widehat{\beta}_{S_Q}\right|_1+\left|\bold{D}_{\bold{Z}}\theta_{ S_{\perp}^c}\right|_1
-\left|\bold{D}_{\bold{Z}}\widehat{\theta}_{ S_{\perp}^c}\right|_1\right)+\widehat{\Sigma}\left(\beta,\theta\right)\notag\\
&\le\frac{\left|\widehat{\Psi}\widehat{\Delta}+\widetilde{\Delta} \right|_{\infty} }{c\widehat{\overline{\kappa}}_{\widehat{h},S,\widetilde{S}}}+
\widehat{\Sigma}\left(\beta,\theta\right).\notag\end{align}
Second, if $2|\bold{D}_{\bold{X}}^{-1}\beta_{S^c\cap S_Q}|_1+2|\bold{D}_{\bold{Z}}\theta_{\widetilde{S}^c}|_1>
|\widehat{\Delta}_{S\cap S_Q}|_1+|\widetilde{\Delta}_{\widetilde{S}}|_1+c(\widehat{\rho}^{ZX}|\widehat{\Delta}|_1+ \underline{r}_n|\widetilde{\Delta}_{S_{\perp}^c}|_1)+|\widehat{\Delta}_{S_Q^c}|_1$, then
we have
\begin{align*} 
\left|\widehat{\Delta}\right|_1+\left|\widetilde{\Delta}_{S_{\perp}^c}\right|_1&=\left|\widehat{\Delta}_{S^c\cap S_Q}\right|_1+\left|\widehat{\Delta}_{S\cap S_Q}\right|_1+\left|\widehat{\Delta}_{S_Q^c}\right|_1+\left|\widetilde{\Delta}_{\widetilde{S}^c}\right|_1+\left|\widetilde{\Delta}_{\widetilde{S}}\right|_1\\
&\le 6\left(\left|\bold{D}_{\bold{X}}^{-1}\beta_{S^c\cap S_Q}\right|_1+\left|\bold{D}_{\bold{Z}}\theta_{\widetilde{S}^c}\right|_1\right).
\end{align*}
For the deterministic lower bounds on the sensitivities we use that, on $\mathcal{G}_{A1}$, denoting by $\overline{\Delta}= D_X^{-1}\bold{D}_{\bold{X}}\widehat{\Delta}$ and $\overline{\widetilde{\Delta}}=D_Z\bold{D}_{\bold{Z}}^{-1}\widetilde{\Delta}$, we have
\begin{align*}
\left|\widehat{\Psi}\widehat{\Delta}+\widetilde{\Delta}\right|_{\infty}&\ge\min_{l\in[{d_Z}]}\left(\bold{D}_{\bold{Z}}D_Z^{-1}\right)_{l,l}
\left|D_Z\mathbb{E}_n\left[ZX^{\top}\right] D_X\overline{\Delta}+\overline{\widetilde{\Delta}}\right|_{\infty}\\
&\ge\frac{1}{\sqrt{1+\tau_n}}\left(\left|D_Z\E\left[ZX^{\top}\right] D_X\overline{\Delta}+
\overline{\widetilde{\Delta}}\right|_{\infty}-\left|D_Z(\E_n-\E)\left[ZX^{\top}\right] D_X\overline{\Delta}\right|_{\infty}
\right)\\
&\ge\frac{1}{\sqrt{1+\tau_n}}\left(\left|\Psi\overline{\widehat{\Delta}}+\overline{\widetilde{\Delta}}\right|_{\infty}-\tau_n\left|\overline{\Delta}\right|_1\right). 
\end{align*}
The remaining arguments are similar to those already seen.\hfill $\square$
\begin{table}

{\footnotesize
  \caption{Table of correspondence for the results on the C-STIV}
\label{TabCorresp}  
\begin{adjustwidth}{-1cm}{}
\begin{tabular}{cc}
\hline
STIV&C-STIV\\
\hline
$\overline{\sigma},\widehat{r},r_n$&$\overline{\sigma}=(\widehat{\sigma}+\widehat{\Sigma}(\widehat\beta,\widehat{\theta}))/2, \underline{r}_n,\underline{r}_n$\\
$\widehat{K}_{S}$&\hspace{-1cm}$\widehat{K}_{S,\widetilde{S}}\triangleq\left\{\begin{array}{l}
(\Delta,\widetilde{\Delta})
:\Delta_{S^c\cap S(\widehat{\beta})^c}=0,\widetilde{\Delta}_{\widetilde{S}^c\cap S(\widehat{\theta})^c}=0,
|\Delta_{S^c\cap S_Q}|_1+|\widetilde{\Delta}_{\widetilde{S}^c}|_1
\le |\Delta_{S\cap S_Q}|_1
+|\widetilde{\Delta}_{\widetilde{S}}|_1+c\widehat{g}(\Delta,\widetilde\Delta)
\end{array}
\right\}$\\
$\widehat{\overline{K}}_{S}$&$\widehat{\overline{K}}_{S,\widetilde{S}}\triangleq\left\{
(\Delta,\widetilde{\Delta}):
|\Delta_{S^c\cap S_Q}|_1+|\widetilde{\Delta}_{\widetilde{S}^c}|_1 
\le 2(|\Delta_{S\cap S_Q}|_1
+|\widetilde{\Delta}_{\widetilde{S}}|_1+c\widehat{g}(\Delta,\widetilde\Delta))+|\Delta_{S_Q^c}|_1
\right\}$\\
$K_{S}$&$K_{S,\widetilde{S}}\triangleq\left\{
\begin{array}{l}
(\Delta,\widetilde{\Delta})
:
(1-\tau_n-c\rho^{ZX}_n)|\Delta|_1+(1-\tau_n- \underline{r}_n)|\widetilde{\Delta}_{S_{\perp}^c}|_1\le2|\Delta_{S\cap S_Q}|_1+|\Delta_{S_Q^c}|_1+2|\widetilde{\Delta}_{\widetilde{S}}|_1
\end{array}
\right\}$\\
$\overline{K}_{S}$&\hspace{-.3cm}$\overline{K}_{S,\widetilde{S}}\triangleq\left\{
\begin{array}{l}
(\Delta,\widetilde{\Delta}):
(1-\tau_n-2c\rho^{ZX}_n)|\Delta|_1
+(1-\tau_n-\underline{r}_n)|\widetilde{\Delta}_{S_{\perp}^c}|_1\le3|\Delta_{S\cap S_Q}|_1+2|\Delta_{S_Q^c}|_1+3|\widetilde{\Delta}_{\widetilde{S}}|_1
\end{array}
\right\}$\\
$|\widehat{\Psi}\Delta|_{\infty}$&$|\widehat{\Psi}\Delta +\widetilde{\Delta}|_{\infty}$\\
$|\bold{D}_{\bold{X}}^{-1}\beta_{S^c\cap S_Q}|_1$&$|\bold{D}_{\bold{X}}^{-1}\beta_{S^c\cap S_Q}|_1+|\bold{D}_{\bold{Z}}\theta_{\widetilde{S}^c}|_1$\\
$| D_X^{-1}\beta_{S^c\cap S_Q}|_1/1_n$&$(| D_X^{-1}\beta_{S^c\cap S_Q}|_1+|D_Z\theta_{\widetilde{S}^c}|_1)/(1-\tau_n)$\\
$\widehat{g}(\Delta)$&$\widehat{g}(\Delta,\widetilde\Delta)\triangleq\widehat{\rho}^{ZX}|\Delta|_1+ \underline{r}_n|\widetilde{\Delta}_{S_{\perp}^c}|_1$\\
$\widehat{h}(\Delta)$&$\widehat{h}(\Delta,\widetilde\Delta)\triangleq\min(|\Delta_{S_Q}|_1+|\widetilde{\Delta}_{S_{\perp}^c}|_1,(3|\Delta_{S\cap S_Q}|_1+3|\widetilde{\Delta}_{\widetilde{S}}|_1+c
\widehat{g}(\Delta,\widetilde\Delta)
+|\Delta_{S_Q^c}|_1)/2)$\\
$\ell^q_{S_0}$&$\ell^q_{S_0,\widetilde{S}_0}(\Delta)\triangleq |\Delta_{S_0}|_q+|\widetilde{\Delta}_{\widetilde{S}_0}|_q$\\
$\widehat\kappa_{\ell_k,S}$&$\displaystyle\widehat\kappa_{\ell_{k},0,S,\widetilde{S}}\triangleq\min_{(\Delta,\widetilde{\Delta})\in \widehat{K}_{S,\widetilde{S}}:\ 
|\Delta_k|=1}|\widehat{\Psi}\Delta +\widetilde{\Delta}|_{\infty},\quad \displaystyle\widehat\kappa_{0,\ell_l,S,\widetilde{S}}\triangleq\min_{(\Delta,\widetilde{\Delta})\in \widehat{K}_{S,\widetilde{S}}:\ 
|\widetilde{\Delta}_l|=1}|\widehat{\Psi}\Delta +\widetilde{\Delta}|_{\infty}$\\
$\widehat{B}(\widehat{S})$&\hspace{-.7cm}$\widehat{B}(\widehat{S},\widehat{\widetilde{S}})\triangleq\left\{
\begin{array}{l}
-\mu\le\Delta\le \mu,\ -\widetilde{\mu}\le\widetilde{\Delta}\le \widetilde{\mu},
\mu_{\widehat{S}^c\cap S(\widehat{\beta})^c}=0,\mu_{\widehat{\widetilde{S}}^c\cap S(\widehat{\theta})^c}=0,\ -\nu1\le\widehat{\Psi}\Delta+\widetilde{\Delta}\le \nu1\\
(1-c\widehat{\rho}^{ZX})(\sum_{j\in {S_I}^c}\mu_j)+
(1-c \underline{r}_n)(\sum_{l\in S_{\perp}^c}\widetilde{\mu}_l)
\le 2
(\sum_{j\in \widehat{S}\cap S_Q}\mu_j+\sum_{l\in \widehat{\widetilde{S}}}\widetilde{\mu}_l)+\sum_{j\in S_Q^c}\mu_j
\end{array}\right\}$\\
$\widehat{B}(k)$&
$\widehat{B}(k,l)\triangleq\left\{\begin{array}{l}
-\mu\le\Delta\le \mu,\ -\widetilde{\mu}\le\widetilde{\Delta}\le \widetilde{\mu},\ -\nu1\le\widehat{\Psi}\Delta+\widetilde{\Delta}\le \nu1\\
(1-c)(\sum_{j\in S_I}\mu_j)+(1-c\widehat{\rho}^{ZX})(\sum_{j\in S_I^c}\mu_j)+
(1-c \underline{r}_n)(\sum_{l\in S_{\perp}^c}\widetilde{\mu}_l)\\
\le 2(s\mu_k+\widetilde{s}\widetilde{\mu}_l)+\sum_{j\in S_Q^c}\mu_j
\end{array}\right\}$\\
$\Gamma_{\kappa}(S)$&$\Gamma_{\kappa}(S,\widetilde{S})\triangleq\gamma(\tau_n/\kappa_{\ell^1_{[{d_X}],\varnothing},S,\widetilde{S}}+\underline{r}_n/(c 1_n\kappa_{\ell^1_{S\cap S_Q,\widetilde{S}}, S,\widetilde{S}}))/(1_n\sqrt{1-\tau_n})$\\
$\Gamma_{\overline{\kappa}}(S)$&$\Gamma_{\overline{\kappa}}(S,\widetilde{S})\triangleq
\gamma(\tau_n/\overline{\kappa}_{\ell^1_{[{d_X}],\varnothing},S,\widetilde{S}}+\underline{r}_n/(c 1_n\overline{\kappa}_{h,S,\widetilde{S}}))/(1_n\sqrt{1-\tau_n})$\\
$\sigma_{U(\beta)}$&$\displaystyle \Sigma(\beta,\theta)\triangleq\max_{l\in[{d_Z}]}(D_Z)_{l,l}\sigma_{T_l(\beta,\theta)}$\\
$\widehat{\beta}^{\widehat\omega}$&$
\widehat{\beta}_{k}^{\widehat\omega}\triangleq\widehat{\beta}_k\indic\left\{
                 |\widehat{\beta}_k|>\widehat\omega_k(s,\widetilde{s})/\mathbb{E}_n[X_k^2]^{1/2}
\right\},\quad \widehat{\theta}_{l}^{\widehat\omega}\triangleq\widehat{\theta}_l\indic\left\{
                  |\widehat{\theta}_l|>\widehat{\widetilde{\omega}}_l(s,\widetilde{s})\mathbb{E}_n[Z_l^2]^{1/2}
\right\}$\\
$ \widehat\omega_k(s)$&$\widehat{\omega}_k(s,\widetilde{s})\triangleq
2 \underline{r}_n\overline{\sigma} \gamma(\underline{r}_n/\widehat\kappa_{\widehat{g}}(s,\widetilde{s}))/\widehat\kappa_{\ell_{k},0}(s,\widetilde{s}),\quad \widehat{\widetilde{\omega}}_l(s,\widetilde{s})\triangleq 2 \underline{r}_n\overline{\sigma}\gamma(\underline{r}_n/\widehat\kappa_{\widehat{g}}(s,\widetilde{s})/\widehat\kappa_{0,\ell_l}(s,\widetilde{s})$\\
$\omega_k(s)$&$\displaystyle\omega_k(s,\widetilde{s})\triangleq\frac{2\underline{r}_n\Sigma(b,t)
\gamma(\underline{r}_n/\kappa_g(s,\widetilde{s}))}{1_n\kappa_{\ell_k,0}(s,\widetilde{s})}
\left(1+\frac{2\underline{r}_n\Gamma_\kappa(S(b),S(t))}{c\kappa_{\ell^1_{S(b)\cap S_Q,S(t)},S(b),S(t)}}\right)$\\
&$\widetilde{\omega}_l(s,\widetilde{s})$ obtained by replacing $\kappa_{\ell_k,0}(s,\widetilde{s})$ by $\kappa_{0,\ell_l}(s,\widetilde{s})$\\
\hline
\end{tabular}
 \end{adjustwidth}}
\end{table}
 
\paragraph{FISTA with Partial Smoothing.}\label{sFISTA}
The C-STIV estimator $(\widehat\beta,\widehat\theta,\widehat\sigma)$ is a solution to a conic program with $d_Z$ cones, and the BC-STIV estimator $(\widehat\Lambda,\widehat\upsilon)$ in \eqref{Mhat} is a solution of a conic program with $d_{\varPhi}d_Z$ cones. If $d_Z$ is large, conic programs are not computationally tractable, so we apply an iterative procedure based on partial smoothing. We present the algorithm for C-STIV, though it can be applied to BC-STIV with minor modifications. Start by noting
\begin{equation*}(\widehat{\beta},\widehat\theta)\in\argmin_{(b,t)\in\mathcal{B}\times\Theta}\left(\frac{1}{c}\left(|\bold{D}_{\bold{X}}^{-1}b_{S_Q}|_1+|\bold{D_Z}t_{S_\perp^c}|_1\right)+\mathcal{O}(b,t)\right),\widehat{\sigma}=\mathcal{O}(b,t),
\end{equation*}
where
$\mathcal{O}(b,t)^2\triangleq
\max_{l\in[d_Z]}\left(\bold{D_Z}\right)_{l,l}^2\max\left(\widehat{\sigma}_l(b,t)^2,\left(\mathbb{E}_n\left[Z_{l}U(b)\right]-t_l\right)^2/\underline{r}_n^2\right)$. We now use ideas from \cite{owen} and \cite{ZZ}. 
Because $2u=\min_{\sigma>0}\left\{\sigma+u^2/\sigma\right\}$, $(\widehat{\beta},\widehat\theta)$ can be obtained by solving  
\begin{equation*}(\widehat{\beta},\widehat\theta,\widehat\sigma)\in\argmin_{(b,t,\sigma)\in\mathcal{B}\times\Theta\times(0,\infty)}\left(\frac{2}{c}\left(|\bold{D}_{\bold{X}}^{-1}b_{S_Q}|_1+|\bold{D_Z}t_{S_\perp^c}|_1\right)+\sigma+\frac{\mathcal{O}(b,t)^2}{\sigma}\right).\end{equation*}
The objective function is convex because $f(x,y)=x^2/y$ is convex on $\R\times(0,\infty)$. Hence, when 
$\mathcal{B}\times\Theta$ is a product, 
a solution of C-STIV can be obtained by the following iterations 
\begin{algm}\label{algm1b} Initialize at $(\widehat\beta^{(0)},\widehat\theta^{(0)},\widehat\sigma^{(0)})$. At iteration $m$, solve
\begin{align*}
&(\widehat{\beta}^{(m)},\widehat\theta^{(m)})\in\argmin_{\substack{(b,t)\in\mathcal{B}}}\left(\frac{2\widehat\sigma^{(m-1)}}{c}\left(\left|\bold{D}_{\bold{X}}^{-1}b_{S_Q}\right|_1+\left|\bold{D}_{\bold{Z}}t_{S_\perp^c}\right|_1\right)+\mathcal{O}(b,t)^2\right),\\
&\widehat\sigma^{(m)}=\mathcal{O}\left(\widehat{\beta}^{(m)},\widehat{\theta}^{(m)}\right),
\end{align*}
then replace $m$ by $m+1$, and iterate until convergence.
\end{algm}
Step 1 can be computationally intensive, whereas Step 2 is trivial. To solve Step 1, we use FISTA with partial smoothing (\cite{BT09,BT12}). Both terms in the minimization problem are convex but nonsmooth; the first involves an $\ell^1$-norm, and the second a maximum. The smoothing is partial because, following \cite{BT12}, we smooth only the maximum, for which we use log-sum-exp smoothing, replacing it by
$$g_\mu(b,t)\triangleq \mu\log\left(\sum_{l\in[d_Z]}
\exp\left(\frac{\widehat{\sigma}_l^2(b,t)}{\mu}\right)+\exp\left(\frac{1}{\mu \underline{r}_n^2}\left(
\left(\bold{D}_{\bold{Z}}\right)_{l,l}\left(\mathbb{E}_n[Z_lU(b)]-t_l\right)
\right)^2\right)\right).$$
Based on Proposition 4.1 and Theorem 3.1 of \cite{BT12}, in practice we take $\mu=\epsilon/(2\log 2d_Z)$ and $\epsilon=0.1$. Smaller values of $\epsilon$ improve the approximation of the maximum but increase the computational burden. After smoothing we are left with the sum of an $\ell^1$-norm and a smooth function, to which we apply FISTA. 

\paragraph{Monte-Carlo for C-STIV.}
We modify the design in Section \ref{sec:sim} to allow for endogenous IVs and apply C-STIV. We take $n=2000$, $d_X=6$, $S_I^c=\{1,5\}$ and $d_Z=50$. We consider a problem with smaller dimensions than for the NV-STIV experiment for computational reasons. Though it is possible to use FISTA to compute C-STIV in applications, in our experiment we need to compute the estimator hundreds of times. The 45 IVs with indices $S_{\perp}^c=\{5,6,...,49\}$ are possibly endogenous. We modify the design by setting $Z_5=\sqrt{1-0.8^2}E+0.8U(\beta^*)$ with $E$ an independent standard Gaussian. This preserves the variance matrix of $Z$ as the identity but implies that $\E[ZU(\beta^*)]=\theta^*$ has one nonzero entry given by $\theta^*_5=0.8$. This is a challenging design since there are fewer IVs known to be exogenous than there are regressors. Table \ref{mc_des4_res} reports C-STIV confidence sets over 1000 replications. The C-STIV estimator detects the endogenous IV, though is downwards biased due to the shrinkage. The confidence sets using a sparsity certificate correctly detect the endogenous IV with frequency 0.74 for $\widetilde{s}=1$, which decreases as $\widetilde{s}$ increases. The confidence sets based on estimated support detect the endogenous IV in every replication.
\setlength{\tabcolsep}{4pt}
\begin{table}
\caption{0.95 C-STIV confidence sets for detection of endogenous IVs}\label{mc_des4_res}
\centering
{\footnotesize
\begin{tabular}{lcccccccccccccc}
\hline
\multicolumn{10}{c}{$d_Z=50,d_X=6,n=2000,\pi=0.8$}\\\hline
&\multicolumn{3}{c}{C-STIV}& SC 4,1 &SC 4,2&SC 4,3&SC 4,4&SC 4,5 & ES \\
\cmidrule(lr){2-4} \cmidrule(lr){5-10}	
				&p2.5 & p50 & p97.5&\multicolumn{6}{c}{Median width/2}  \\
\cmidrule(lr){2-4} \cmidrule(lr){5-10}	
$\theta^*_5(=0.8)$&0.66      &    0.71       &    0.77       &    0.69       &    0.71       &    0.73       &    0.76       &    0.78       &    0.25    
\\
$\theta^*_6(=0)$   &     0       &         0       &         0       &    0.69       &    0.71      &    0.74       &    0.76       &    0.78       &         0
\\
\hline
$S(\widehat\theta)\supseteq S(\theta^*)$&&1&Power&.74&.51&.26&.11&.04&1\\
$S(\widehat\theta)= S(\theta)$&&1       &             &(.71,.77)              &(.48,.54)        & (.23,.29)         & (.09,.13)    & (.03,.05)       & (.996,1)   \\
                    \hline
\multicolumn{10}{p{120mm}}{\scriptsize \textbf{Notes:} 1000 replications. `SC $s,\widetilde{s}$' use sparsity certificates $s,\widetilde{s}$. `ES' use estimated support. SC/ES use one grid point for $c$. $\underline{r}_n=0.07$. `C-STIV' uses $c=0.99$. `Power' is the frequency with which the confidence sets do not include $\theta_5=0$.
}
\end{tabular}
}
\vspace{-3mm}
\end{table}

\subsubsection{Confidence Bands under Conditional Homoskedasticity}\label{conhomban}
Let us consider confidence bands when we maintain (C4.\ref{Bi}). We present these for a structural model with approximation errors, as in the proof of Theorem \ref{tCI}. The confidence bands are $\widehat{C}_{\varPhi}\triangleq\left[\widehat{\underline{C}}_{\varPhi},\widehat{\overline{C}}_{\varPhi}\right]$ where
\begin{align}
&\widehat{\underline{C}}_{\varPhi}\triangleq\widehat{\varPhi\beta}-\widehat{q},\widehat{\overline{C}}_{\varPhi}\triangleq\widehat{\varPhi\beta}+\widehat{q}\label{bndOmega1},\quad\widehat{q}\triangleq\frac{q_{G_{\varPhi}|\bold{Z}\widehat\Lambda^\top}(1-\alpha)+3\zeta_n}{\sqrt{n}}\widehat{\sigma}(\widehat{\beta})\bold{D}_{\bold{Z}\widehat{\Lambda}^{\top}}^{-1}1,\notag
\end{align}
$G_{\varPhi}=\sqrt{n}|\bold{D}_{\bold{Z}\widehat{\Lambda}^{\top}}\widehat{\Lambda}\mathbb{E}_n[ZE]|_{\infty}$, and $\bold{E}\in\R^n$ is a standard Gaussian vector independent of $\bold{Z}\widehat{\Lambda}^{\top}$.  
For the analysis, we introduce the deterministic upper bound $v^{\sigma_{W(\beta)}}_n$ such that on $\mathcal{G}\cap\mathcal{G}_{A1}$ of probability at least  $1-\alpha_n^{S}$, 
$|\widehat{\sigma}(\widehat\beta)-\sigma_{W(\beta)}|\le v^{\sigma_{W(\beta)}}_n$, 
which is obtained from \eqref{esigma}.  
We also replace \eqref{d.g.p.4a}-\eqref{d.g.p.4b2} in Assumption \ref{assCI0} by (N.\ref{Ni}) holds for $Z$ and $M_{Z}$, and
\begin{itemize}
\item[(iv')] (N.\ref{Ni}) holds for $\Lambda Z$ and $M_{\Lambda Z}$; 
\item[(v')] $\hspace{-.2cm}\left|\hspace{-.1cm}\left(\hspace{-.1cm}\max\left(\hspace{-.1cm}\E\hspace{-.1cm}\left[\hspace{-.1cm}\left(\hspace{-.1cm}\left(D_{\Lambda Z}\Lambda\right)_{f,\cdot} \hspace{-.05cm}G(\beta)/\sigma_{W(\beta)}\hspace{-.1cm}\right)^{2+q_1}\hspace{-.1cm}\right],
\E\left[\hspace{-.1cm}\left(\hspace{-.05cm}\left(D_{\Lambda Z}\Lambda\right)_{f,\cdot} \hspace{-.05cm}G(\beta)E/\sigma_{W(\beta)}\hspace{-.1cm}\right)^{2+q_1}\hspace{-.1cm}\right]\hspace{-.05cm}
\right)\hspace{-.1cm}\right)_{f=1}^{d_{\varPhi}}\right|_{\infty}\hspace{-.4cm}
\le B_n^{q_1}$;
\item[(vi')] $\hspace{-.15cm}\max\left(\E\left[\left(\left|D_{\Lambda Z }\Lambda G(\beta)\right|_{\infty}\hspace{-.1cm}/\hspace{-.05cm}(B_n\sigma_{W(\beta)})\right)^{q_2}\hspace{-.05cm}\right],
\E\left[\left(\left|D_{\Lambda Z }\Lambda G(\beta)E\right|_{\infty}\hspace{-.1cm}/\hspace{-.05cm}(B_n\sigma_{W(\beta)}\right)^{q_2}\hspace{-.05cm}\right]\right)\hspace{-.1cm}\le\hspace{-.1cm} 2$.
\end{itemize}
The loss $|D_{\Lambda G(\beta)}\cdot|_{\infty,\infty}$ is replaced by 
$|D_{\Lambda Z}\cdot|_{\infty,\infty}$. 

The coverage is guaranteed, with coverage error 
$\alpha^{A2}_n\triangleq
2\zeta'_n+\zeta''_n+\varphi(d_{\varPhi},\tau_n)+\iota(d_{\varPhi},n)+\alpha_n(\Lambda Z)
+\alpha^{S}_n+\alpha^{BC}_n$, where
$(\zeta'_n)^2\triangleq 
3\alpha_n+\alpha_n(ZE)'+\alpha^{BC}_n
+\iota(d_{\varPhi},n)$ and $\zeta''_n\triangleq 
\alpha_n+
\alpha^{BC}_n+\alpha_n^{S}+\iota(d_{\varPhi},n)$, if we assume:  
\begin{ass}\label{assCI1b} 
For all $(\beta,\Lambda)\in\mathcal{I}_{\varPhi}$, we have
\begin{enumerate}[\textup{(}i\textup{)}]
\item $\max(v_n^{\Lambda,\beta}B_Z+\tau_n,v^{\sigma_{W(\beta)}}_n,2v_n^D)<1$;
\item\label{ezeta} $\zeta_n\ge\max\left(2\sqrt{n}v^{\Lambda,\beta}_nr_n^E,2\sqrt{n}v^{\Lambda,\beta}_nr_n/1_n,
4v_n^D\log\left(2d_\varPhi/\alpha_n\right)
/(1-2v_n^D),v_n^R\right)$; 
\end{enumerate}
where
\begin{align*}
v_n^{D}&\triangleq (v_n^{\Lambda,\beta}B_Z+\tau_n)(1-v^{\sigma_{W(\beta)}}_n)+v^{\sigma_{W(\beta)}}_n,\\
v^R_n&\triangleq  \sqrt{n}\left(\left|D_{\Lambda G(\beta)}\right|_{\infty}\left(\underline{r}'_n v^{\Sigma(\Lambda)}_n v^{\beta}_n+|\overline{V}(\beta)|_{\infty}\right)
+v_{d_X}\sqrt{1+\tau_n}\right)/(\sigma_{W(\beta)}(1-v_n^D)),\\
r_n^E&\triangleq \underline{r}_n^E
2\log\left(2n/\alpha_n\right)\sqrt{1+\tau_n},
\end{align*}
and $\underline{r}_n^E$ is obtained like $\underline{r}_n$ for Class 1 replacing $\alpha$ by $\alpha_n$. 
\end{ass}
The main arguments of the proof are detailed for Class 4 later in this appendix. The specific elements are the following. 
On $\mathcal{E}\triangleq\underline{\mathcal{G}}'\cap\mathcal{G}_{A1}
\cap\mathcal{E}_{T}'^{c}\cap\mathcal{E}_Z^c$
 of probability $1-\alpha^{BC}_n$, we have
\begin{align*}
&\left(D_{\Lambda Z}\right)_{f,f}\left|\mathbb{E}_n\left[\left(\widehat{\Lambda}_{f,\cdot}Z\right)^2\right]^{1/2}-\mathbb{E}\left[(\Lambda_{f,\cdot}Z)^2\right]^{1/2}\right|\\
&\le
\left(D_{\Lambda Z}\right)_{f,f}\left(\mathbb{E}_n\left[\left(\left(\widehat{\Lambda}_{f,\cdot}-\Lambda_{f,\cdot}\right)Z\right)^2\right]^{1/2}+
\left|\mathbb{E}_n\left[(\Lambda_{f,\cdot}Z)^2\right]^{1/2}-\mathbb{E}\left[(\Lambda_{f,\cdot}Z)^2\right]^{1/2}\right|\right)\\
&\le v_n^{\Lambda,\beta}B_Z+\tau_n\le v_n^D,
\end{align*}
so $\left|\bold{D}_{\bold{Z}\widehat{\Lambda}^{\top}}D_{\Lambda Z}^{-1}\right|_{\infty}
\le 1/(1-v_n^D)$ and, on $\mathcal{E}\cap\mathcal{G}$,  
$\left|\bold{D}_{\bold{Z}\widehat{\Lambda}^{\top}}D_{\Lambda Z}^{-1}\right|_{\infty}\sigma_{W(\beta)}/\widehat{\sigma}(\widehat{\beta})\le 1/(1-v_n^D)$. 
We now use
\begin{align*}
&\frac{\sqrt{n}}{\widehat{\sigma}(\widehat{\beta})}\bold{D}_{\bold{Z}\widehat{\Lambda}^{\top}}\left(\widehat{\varPhi\beta}-\varPhi\beta-\overline{V}(\beta)\right)=R+\frac{\sqrt{n}}{\widehat{\sigma}(\widehat{\beta})}\bold{D}_{\bold{Z}\widehat{\Lambda}^{\top}}\widehat{\Lambda}\mathbb{E}_n[G(\beta)],\\
&R\triangleq\frac{\sqrt{n}}{\widehat{\sigma}(\widehat{\beta})}\bold{D}_{\bold{Z}\widehat{\Lambda}^{\top}}\left(\varPhi-\widehat{\Lambda}\mathbb{E}_n[ZX^{\top}]\right)\bold{D}_{\bold{X}}\widehat\Delta-\frac{\sqrt{n}}{\widehat{\sigma}(\widehat{\beta})}\bold{D}_{\bold{Z}\widehat{\Lambda}^{\top}}\overline{V}(\beta)+\frac{\sqrt{n}}{\widehat{\sigma}(\widehat{\beta})}\bold{D}_{\bold{Z}\widehat{\Lambda}^{\top}}\widehat{\Lambda}\mathbb{E}_n[ZV(\beta)].
\end{align*} 
On $\mathcal{E}\cap\mathcal{G}$ of probability $1-(\alpha^{S}_n+\alpha^{BC}_n)$, we have 
$|R|_{\infty}\le v_n^R$. 
Define 
\begin{align*}
&T_{\varPhi}\triangleq\left|\frac{\sqrt{n}}{\widehat{\sigma}(\widehat\beta)}\bold{D}_{\bold{Z}\widehat{\Lambda}^{\top}}\widehat{\Lambda}\mathbb{E}_n[G(\beta)]\right|_{\infty}\hspace{-.3cm},\ T_{\varPhi1}=\left|\frac{\sqrt{n}}{
\sigma_{W(\beta)}
}D_{\Lambda Z}\widehat{\Lambda} \mathbb{E}_n[G(\beta)]\right|_{\infty}\hspace{-.3cm},\ 
T_{\varPhi0}\triangleq\left|\frac{\sqrt{n}}{\sigma_{W(\beta)}}D_{\Lambda Z}\Lambda \mathbb{E}_n[G(\beta)]\right|_{\infty}\hspace{-.3cm},\\ 
&G_{\varPhi1}\triangleq\left|\sqrt{n}D_{\Lambda Z}\widehat{\Lambda} \mathbb{E}_n[ZE]\right|_{\infty},\ 
G_{\varPhi0}\triangleq\left|\sqrt{n}D_{\Lambda Z}\Lambda \mathbb{E}_n[ZE]\right|_{\infty}.
\end{align*}
On $\mathcal{E}\cap\mathcal{G}$, $\left|T_{\varPhi}-T_{\varPhi1}\right|\le T_{\varPhi1}v_n^D/(1-v_n^D)$ and $\left|T_{\varPhi1}-T_{\varPhi0}\right|\le v^{\Lambda,\beta}_n\sqrt{n}r_n/1_n
$, so
$$\left|T_{\varPhi}-T_{\varPhi0}\right|\le (T_{\varPhi0}+v^{\Lambda,\beta}_n\sqrt{n}r_n/1_n)v_n^D/(1-v_n^D)+ v^{\Lambda,\beta}_n\sqrt{n}r_n/1_n.$$
Let $\mathcal{E}_{2,E}\triangleq\left\{\forall l\in[d_Z],|\mathbb{E}_n[Z_lE]|> \underline{r}_n^E\mathbb{E}_n[(Z_lE)^2]^{1/2}\right\}$. On $\mathcal{E}\cap\mathcal{E}_{2,E}\cap \mathcal{E}_{ZE}'^c\cap\{|\bold{E}|_{\infty}\le 2\log(2n/\alpha_n)\}$, 
$\left|G_{\varPhi}-G_{\varPhi1}\right|\le G_{\varPhi1}v_n^D/(1-v_n^D)$ and 
$\left|G_{\varPhi1}-G_{\varPhi0}\right|\le \sqrt{n}v^{\Lambda,\beta}_nr_n^E$,
so
$$\left|G_{\varPhi}-G_{\varPhi0}\right|\le
(G_{\varPhi0}+\sqrt{n}v^{\Lambda,\beta}_nr_n^E)v_n^D/(1-v_n^
{D,2})+\sqrt{n}v^{\Lambda,\beta}_nr_n^E.$$
By Assumption \ref{assCI1b} \eqref{ezeta},  we have
\begin{align*}
(\zeta_n-\sqrt{n}v^{\Lambda,\beta}_nr_n^E)(1-v_n^D)/v_n^D-\sqrt{n}v^{\Lambda,\beta}_nr_n^E&\ge \zeta_n(1-2v_n^D)/(2v_n^D),\\
(\zeta_n-\sqrt{n}v^{\Lambda,\beta}_nr_n/1_n)(1-v_n^D)/v_n^D-\sqrt{n}v^{\Lambda,\beta}_nr_n/1_n&\ge \zeta_n(1-2v_n^D)/(2v_n^D).
\end{align*}
By \eqref{ezeta} and $2\log\left(2d_\varPhi/\alpha_n\right)\ge q_{N_{\varPhi0}}(1-\alpha_n)$, where 
$N_{\varPhi0}\triangleq\left|D_{\Lambda Z}\Lambda E_{Z}\right|_{\infty}$ and $E_{Z}$
is a Gaussian vector of covariance $\mathbb{E}[ ZZ^{\top}]$, we get 
$$\mathbb{P}\left(|T_\varPhi-T_{\varPhi 0}|>\zeta_n \right)\leq\zeta''_n\quad \text{and}\quad 
\mathbb{P}\left(\mathbb{P}\left(|G_{\varPhi}-G_{\varPhi0}|>\zeta_n| \bold{Z}\widehat\Lambda^\top\right)>\zeta'_n\right)<\zeta'_n.$$ 
This yields, as in the proof of Class 4, $\mathbb{P}\left(T_{\varPhi}\ge q_{G_\varPhi|\bold{Z}\widehat{\Lambda}^\top}(1-\alpha)+2\zeta_n\right)<\alpha+2\zeta'_n+\zeta''_n+\varphi(d_{\varPhi},\tau_n)+\iota(d_{\varPhi},n)+\alpha_n(\Lambda Z)$, hence the result.  

\subsubsection{Results in Previous Versions of this Paper: \cite{GRT}}
The interested reader can find results for 
STIV confidence sets with a high-dimensional version of 2SLS and its failure in various situations. There are results for other approaches  than the one based on sparsity certificates for NV-STIV. 
 The C-STIV is a simple modification of the NV-STIV and was introduced to answer a referee's comment in 2011 on ways to avoid loosing $\left|\bold{D}_{\bold{Z}}\bold{Z}^{\top}\right|_{\infty}$ when setting $\widehat{r}=\underline{r}_n\left|\bold{D}_{\bold{Z}}\bold{Z}^{\top}\right|_{\infty}$. It was the STIV estimator in the revision of this paper between 2012 and 2014 (\href{https://arxiv.org/pdf/1812.11330.pdf}{first revision (2012)}). 
 In the previous versions we also propose an assumption which allows to obtain tighter bounds in the same spirit as the treatment of the regressors in $S_I$ for the C-STIV. 
We propose confidence bands with bias correction 
using sample splitting. We also study the combination of the confidence bands with an upper bound on the bias obtained from the identification robust confidence sets in case we suspect the "bias" of the debiased estimator might not be negligible. 
 An alternative identification robust confidence set, also not involving test inversion, called SNIV set, relies on semidefinite relaxations and can deliver tight sets when $d_Z<d_X$ and when there are endogenous IVs but $S_{\perp}$ is not available.

\subsection{Proofs of Results in the Appendix}\label{proofAp}
\noindent\noindent{\bf Proof of Proposition \ref{p4} and \ref{p4b}.} We prove the bounds for the sensitivities based on $\widehat{K}_{S}$, those for  the sensitivities based on $\widehat{\overline{K}}_{S}$ are obtained similarly. Parts \eqref{p410} and \eqref{p411} are easy.\\ 
The upper bound in the first display in \eqref{p42i} follows from $|\Delta_{S_0}|_q\ge|\Delta_{S_0}|_{\infty}$.  We obtain the lower bound as follows. Because
$|\Delta_{S_0}|_q\le|\Delta_{S_0}|_1^{1/q}|\Delta_{S_0}|_{\infty}^{1-1/q}$, we get that, for $\Delta\ne0$,
\begin{equation}\label{Gcone0}
\frac{\left|\widehat{\Psi}\Delta\right|_{\infty}}{|\Delta_{S_0}|_q}\ge
\frac{\left|\widehat{\Psi}\Delta\right|_{\infty}}{|\Delta_{S_0}|_{\infty}}
\left(\frac{|\Delta_{S_0}|_{\infty}}{|\Delta_{S_0}|_1}\right)^{1/q}
\end{equation}
and use $|\Delta_{S_0}|_1\le |S_0||\Delta_{S_0}|_{\infty}$. 
Furthermore, for $\Delta \in \widehat{K}_{S}$, by definition of the set, we have $|\Delta_{S^c\cap S_Q}|_1\le |\Delta_{S\cap S_Q}|_1+c\widehat{r} |\Delta|_1+ c(1-\widehat{r}) |\Delta_{S_I^c}|_1$,
which, by adding $|\Delta_{(S\cap S_Q)\cup S_Q^c}|_1$ on both sides, gives
\begin{equation}\label{eq:boundl11}
|\Delta|_1\le \gamma(c\widehat{r})\left(2|\Delta_{S\cap S_Q}|_1+|\Delta_{S_Q^c}|_1+
c(1-\widehat{r})|\Delta_{S_I^c}|_1\right).
\end{equation}
From \eqref{eq:boundl11} and the fact that $\Delta_{S^c\cap S(\widehat{\beta})^c}=0$, we deduce
\begin{equation}\label{eother}
|\Delta|_1
\le \gamma(c\widehat{r})\left(2|S\cap S_Q|+\left|S_Q^c\cap\left(S\cup S(\widehat{\beta})\right)\right|+c(1-\widehat{r})\left|S_I^c\cap\left(S\cup S(\widehat{\beta})\right)\right|\right)\left|\Delta_{\widehat{S}(S,S(\widehat{\beta}))}\right|_{\infty}.
\end{equation}
Let us obtain an alternative lower bound for the case $c\in(0,1)$. The condition $\Delta\in \widehat{K}_{S}$ can be written as $|\Delta_{S^c\cap S_Q}|_1\le |\Delta_{S\cap S_Q}|_1+c(\widehat{r}-1) |\Delta_{S_I}|_1+ c |\Delta|_1$
which implies $|\Delta_{S^c\cap S_Q}|_1\le |\Delta_{S\cap S_Q}|_1+ c |\Delta|_1$ and, by adding $|\Delta_{(S\cap S_Q)\cup S_Q^c}|_1$ on both sides, if $c\in(0,1)$, we have 
\begin{equation}\label{Gconeb3}
|\Delta|_1\le\gamma(c)\left(2 |\Delta_{S\cap S_Q}|_1+  |\Delta_{S_Q^c}|_1\right).
\end{equation}
Using $\Delta_{S^c\cap S(\widehat{\beta})^c}=0$, this yields
\begin{equation}\label{eq:boudl12}
|\Delta|_1\le \gamma(c)\left(2|S\cap S_Q|+\left|S_Q^c\cap\left(S\cup S(\widehat{\beta})\right)\right|\right)\left|\Delta_{\widehat{S}(S,S(\widehat{\beta}))}\right|_{\infty}.
\end{equation}
Combining \eqref{eother} and \eqref{eq:boudl12} yields 
\begin{equation}\label{galbound}
|\Delta|_1\le\widehat{c}_{\kappa}(S)\left|\Delta_{\widehat{S}(S)}\right|_{\infty}. 
\end{equation}
This yields the second display in \eqref{p42i}. The first lower bound in the first display in \eqref{p42i} uses
\begin{equation*}
\frac{\left|\widehat{\Psi}\Delta\right|_{\infty}}{\left|\Delta_{S_0}\right|_q}\ge
\frac{\left|\widehat{\Psi}\Delta\right|_{\infty}}{\left|\Delta_{S_0\cup\widehat{S}(S,S(\widehat{\beta}))}\right|_q}\ge
\frac{\left|\widehat{\Psi}\Delta\right|_{\infty}}{\left|\Delta_{S_0\cup\widehat{S}(S,S(\widehat{\beta}))}\right|_{\infty}}
\left(\frac{\left|\Delta_{S_0\cup\widehat{S}(S,S(\widehat{\beta}))}\right|_{\infty}}{\left|\Delta_{S_0\cup\widehat{S}(S,S(\widehat{\beta}))}\right|_1}\right)^{1/q}
\end{equation*}
and $|\Delta_{S_0\cup\widehat{S}(S)}|_1\le\widehat{c}_{\kappa}(S)|\Delta_{S_0\cup \widehat{S}(S)}|_{\infty}$ which can be deduced from \eqref{galbound}.\\ 
To prove \eqref{bk1} it suffices to note that, by definition of the set $\widehat{K}_{S}$,
\begin{align}
|\Delta|_1&\le \left(\frac{2}{\widehat\kappa_{\ell^1_{S\cap S_Q},S}} +\frac{1}{\widehat\kappa_{\ell^1_{S_Q^c},S}}+
\frac{c}{\widehat\kappa_{\widehat{g},S}}\right)\left|\widehat{\Psi}\Delta\right|_{\infty},
\end{align}
by \eqref{eq:boundl11}, $$|\Delta|_1\le \gamma(c\widehat{r})\left(
\frac{2}{\widehat\kappa_{\ell^1_{S\cap S_Q},S}} +\frac{1}{\widehat\kappa_{\ell^1_{S_Q^c},S}}+\frac{c(1-\widehat{r})}{\widehat\kappa_{\ell^1_{S_I^c}, S}}  
\right)\left|\widehat{\Psi}\Delta\right|_{\infty},$$
and, by \eqref{Gconeb3}, 
$$|\Delta|_1\le \gamma(c)\left(\frac{2}{\widehat\kappa_{\ell^1_{S\cap S_Q},S}}+ \frac{1}{\widehat\kappa_{\ell^1_{S_Q^c},S}}\right)\left|\widehat{\Psi}\Delta\right|_{\infty}.$$ 
The bound \eqref{p46i2} is obtained by rewriting $\Delta\in \widehat{K}_{S}$ as $(1-c\widehat{r})|\Delta_{S_I}|_1+(1-c)|\Delta_{S_I^c}|_1\le 2|\Delta_{S\cap S_Q}|_1+  |\Delta_{S_Q^c}|_1$, which yields
\begin{align}
\widehat{r}|\Delta_{S_I}|_1+|\Delta_{S_I^c}|_1&\le \widehat{r}\gamma(c\widehat{r})\left(2|\Delta_{S\cap S_Q}|_1+  |\Delta_{S_Q^c}|_1
+\frac{1-\widehat{r}}{\widehat{r}}|\Delta_{S_I^c}|_1\right)\label{exactks0}\\
&\le \left|\widehat{\Psi}\Delta\right|_{\infty}\widehat{r}\gamma(c\widehat{r})\left(
\frac{2}{\widehat\kappa_{\ell^1_{S\cap S_Q},S}}+
\frac{1}{\widehat\kappa_{\ell^1_{S_Q^c},S}}
+\frac{1-\widehat{r}}{\widehat{r}\widehat\kappa_{\ell^1_{S_I^c}, S}}\right)\notag.
\end{align}
The second upper bound follows from noticing that, if $\widehat\kappa_{\widehat{g},S}>0$, we have
$$\frac{1}{\widehat\kappa_{\widehat{g},S}}=\sup_{\Delta\in \widehat{K}_{S}:\ \left|\widehat{\Psi}\Delta \right|_{\infty}=1}\left(\widehat{r}\left|\Delta_{S_I}\right|_1+\left|\Delta_{S_I^c}\right|_1\right)\le 
\sup_{\Delta\in \widehat{S}_{S}:\ \left|\widehat{\Psi}\Delta \right|_{\infty}=1}\widehat{r}\left|\Delta_{S_I}\right|_1+
\sup_{\Delta\in \widehat{K}_{S}:\ \left|\widehat{\Psi}\Delta \right|_{\infty}=1}\left|\Delta_{S_I^c}\right|_1.$$
Let us now prove \eqref{p44i}. Since for all $k$ in $S_0$, $|\Delta_{S_0}|_{\infty}\ge|\Delta_k|$, for all $k$ in $S_0$, 
$$\widehat\kappa_{\ell^\infty_{S_0},S}=\min_{\Delta\in \widehat{K}_{S_0}}\frac{\left|\widehat{\Psi}\Delta\right|_{\infty}}{|\Delta_{S_0}|_{\infty}}\le \min_{\Delta\in \widehat{K}_{S}}\frac{\left|\widehat{\Psi}\Delta\right|_{\infty}}{|\Delta_k|}=\widehat\kappa_{\ell_{k},S}.$$
Thus $\widehat\kappa_{\ell^\infty_{S_0},S}\le \min_{k\in S_0}\widehat\kappa_{\ell_{k},S}$. 
But one also has 
$$
\widehat\kappa_{\ell^\infty_{S_0},S}=\min_{k\in S_0}\min_{\Delta\in \widehat{K}_{S}:\ |\Delta_k|=|\Delta_{S_0}|_{\infty}=1}\left|\widehat{\Psi}\Delta\right|_{\infty}\ge\min_{k\in S_0}\min_{\Delta\in \widehat{K}_{S}:\ |\Delta_k|=1}\left|\widehat{\Psi}\Delta\right|_{\infty}.\quad\quad\quad\quad\quad\quad\quad\square
$$

\noindent {\bf Proof of Lemma \ref{lconc}.} We prove the middle statement:
\begin{align*}
\mathbb{P}\left(\mathcal{E}_{A}'\right)&=\mathbb{P}\left(\left|\left(\mathbb{E}_n\left[\frac{A_{l}^2}{\E\left[A_l^2\right]}-1\right]\right)_{l\in[d_A]}\right|\ge \tau_n\right)\\
&\le \frac{1}{\tau_n^2}\E\left[\max_{l\in[d_A]}\left|\mathbb{E}_n\left[\frac{A_{l}^2}{\E\left[A_l^2\right]}-1\right]\right|^2\right]\quad \mbox{(by\ the Chebyshev inequality)}\\
&\le \frac{C_{{\rm N}}(d_A)}{n\tau_n^2}\E\left[\max_{l\in[d_A]}\left|\left(\frac{A_{l}^2}{\E\left[A_l^2\right]}-1\right)\right|^2\right]\quad \mbox{(by the Nemirovski inequality)}\\
&\le \alpha_n(A)'.
\end{align*}
The proof of the remaining statements is the same.\hfill $\square$\vspace{0.2cm}

\noindent{\bf Proof of lemmas \ref{thrm:DLBsensitivities}.} 
Clearly, on $\mathcal{E}_X'^c$, the following holds
\begin{align}
&\forall b\in\R^{d_X},\ \ell\in\mathcal{L},\ \sqrt{1-\tau_n}\ell\left( D_X^{-1}b\right)\le \ell\left(\bold{D}_{\bold{X}}^{-1}b\right)\le \sqrt{1+\tau_n}\ell\left( D_X^{-1}b\right),\label{thrm:DLBsensitivitiesl}.
\end{align}
 Assume now that we work on the event $\mathcal{G}_{A1}$. Let $S\subseteq[{d_X}],\ell\in\mathcal{L}$, and $\overline{\Delta}\triangleq  D_X^{-1}\bold{D}_{\bold{X}}\Delta$. 
Due to \eqref{thrm:DLBsensitivitiesl}, we have 
$\sqrt{1-\tau_n}\ell(\overline{\Delta})\le \ell(\Delta)\le \sqrt{1+\tau_n}\ell(\overline{\Delta})$.
This, the fact that $\widehat{r}\le r$, and manipulations on the $\ell^1$-norm of subvectors used previously, yield 
$\overline{\Delta}\in K_{S}$ if $\Delta\in\widehat{K}_{S}$ and $\overline{\Delta}\in \overline{K}_{S}$ if $\Delta\in\widehat{\overline{K}}_{S}$. Now, because $\mathcal{G}_{A1}\subseteq \mathcal{E}_Z'^c\cap \mathcal{E}_{ZX^{\top}}^c$, we obtain
\begin{align}
\left|\widehat{\Psi}\Delta\right|_{\infty}&\ge\min_{l\in[{d_Z}]}\left(\bold{D}_{\bold{Z}}D_Z^{-1}\right)_{l,l}
\left|D_Z\mathbb{E}_n\left[ZX^{\top}\right] D_X D_X^{-1}\bold{D}_{\bold{X}}\Delta\right|_{\infty}\notag\\
&\ge\left(
\left|\Psi\overline{\Delta}\right|_{\infty}
-\tau_n\left|\overline{\Delta}\right|_1
\right)/\sqrt{1+\tau_n}.\label{eNPIV0}
\end{align}
\eqref{thrm:DLBsensitivities1b} 
is obtained from the definition of  $\kappa_{\ell^1,S}$ and $\overline{\kappa}_{\ell^1,S}$ and  that, on $\mathcal{G}_{A1}$, 
$\ell(\Delta)\le \sqrt{1+\tau_n}\ell(\overline{\Delta}).\square$\vspace{0.2cm}

\noindent{\bf Proof of Proposition \ref{p5}.}  \eqref{p50i}, the first identity in \eqref{p5ii}, and \eqref{p501i} are obtained like the similar results in Proposition \ref{p4}. By H\"older's inequality, Condition IC, for all $\Delta\in K_S$ and $q\in[1,\infty]$,  
$$|\Delta|_1\le u_{\kappa}|\Delta_S|_1\le u_{\kappa}|S|^{1-1/q}|\Delta_S|_q,$$
from which we deduce \eqref{p5iib}. A similar result also holds for the sensitivities (see \cite{GRT}).
We obtain \eqref{galboundd} by similar arguments as those used for \eqref{galbound}. 
This yields the second identity in \eqref{p5ii}.
To prove \eqref{p5i}, let $\lambda\in\R^{d_Z}$ such that $|\lambda|_1\le1$, $k\in [d_X]$, and $\Delta\in K_S$. By the inverse triangle inequality
$$\left|\lambda^{\top}\Psi\Delta-\lambda^{\top}\Psi_{\cdot,k}\Delta_k\right|\le\left(\sum_{k'\ne k}\left|\Delta_{k'}\right|\right)\max_{k'\ne k}\left|\lambda^{\top}\Psi_{\cdot,k'}\right|,$$
which yields $\left|\lambda^{\top}\Psi_{\cdot,k}\right|\left|\Delta_k\right|\le\left(\sum_{k'\ne k}\left|\Delta_{k'}\right|\right)\max_{k'\ne k}\left|\lambda^{\top}\Psi_{\cdot,k'}\right|+\left|\lambda^{\top}\Psi\Delta\right|$, hence
\begin{align}
\left(\left|\lambda^{\top}\Psi_{\cdot,k}\right|+\max_{k'\ne k}\left|\lambda^{\top}\Psi_{\cdot,k'}\right|\right)\left|\Delta_k\right|&\le\left|\Delta\right|_1\max_{k'\ne k}\left|\lambda^{\top}\Psi_{\cdot,k'}\right|+\left|\Psi\Delta\right|_{\infty},\label{econtinue}\\
&\le c_{\kappa}(S)\max_{k'\ne k}\left|\lambda^{\top}\Psi_{\cdot,k'}\right|\left|\Delta_{\overline{S}\cup S_0}\right|_{\infty}+\left|\Psi\Delta\right|_{\infty}.\notag
\end{align}
The last display uses $|\Delta|_1\le c_{\kappa}(S)|\Delta_{\overline{S}}|_{\infty}$. 
For $k$ such that $|\Delta_k|=|\Delta_{\overline{S}\cup S_0}|_{\infty}$, 
\begin{equation}\label{eqcont}
\max_{\lambda\in\R^{d_Z}:\ |\lambda|_1\le1}\left(\left|\lambda^{\top}\Psi_{\cdot,k}\right|-(c_{\kappa}(S)-1)\max_{k'\ne k}\left|\lambda^{\top}\Psi_{\cdot,k'}\right|\right)\left|\Delta_{\overline{S}\cup S_0}\right|_{\infty}\le\left|\Psi\Delta\right|_{\infty}
\end{equation}
and we conclude by taking the minimum over $k\in \overline{S}\cup S_0$ and the definition of the $\ell^{\infty}_{\overline{S}\cup S_0}$ sensitivity. We can remove the first absolute value by changing $\lambda$ in $-\lambda$. 
To prove \eqref{p5iii} we start from \eqref{econtinue}. By definition of the $\ell^{1}$ population sensitivity we have 
$$\left(\left|\lambda^{\top}\Psi_{\cdot,k}\right|+\max_{k'\ne k}\left|\lambda^{\top}\Psi_{\cdot,k'}\right|\right)\left|\Delta_k\right|\le\left(\frac{\max_{k'\ne k}\left|\lambda^{\top}\Psi_{\cdot,k'}\right|}{\kappa_{\ell^1,S}}+1\right)\left|\Psi\Delta\right|_{\infty}.$$
We conclude by setting $\Delta_k=1$. The other items are proved similarly to Proposition \ref{p4}. 
\hfill$\square$\vspace{0.2cm}

\noindent{{\bf Proof of Proposition \ref{tCLambda}.}}
Take $(\beta,\Lambda)\in\mathcal{I}_{\varPhi}$.
Set $\widehat\Delta'\triangleq(\widehat{\Lambda}-\Lambda)\bold{D}_{\bold{Z}}^{-1}$, and 
$\overline{\widehat\Delta'}\triangleq\widehat\Delta'\bold{D}_{\bold{Z}}D_Z^{-1}$. 
Clearly, on $\underline{\mathcal{G}}'$, $\Lambda$ belongs to $\widehat{\mathcal{I}}_{\varPhi}\left(\underline{r}'_n,\widehat{\Sigma}(\Lambda\right)$. We now work on the event in the statement of the theorem. 
We start by proving \eqref{tCLambdai}. 
The arguments in the proof of Proposition \ref{t1} yield
\begin{align}
\left|\widehat\Delta'\widehat{\Psi}^{\top}\right|_{\infty} &\le
\underline{r}'_n\left(\widehat{\nu}+\widehat{\Sigma}\left(\Lambda\right)\right)\label{stop2CCB}\\
\left|\widehat\Delta'_{S(\Lambda)^c}\right|_{1}&\le
\left|\widehat\Delta'_{S(\Lambda)}\right|_{1}+\frac{\lambda}{\widehat{\rho}^{ZX}}\left(\widehat{\Sigma}\left(\Lambda\right)-\widehat{\Sigma}\left(\widehat\Lambda\right)\right)\nonumber
\end{align}
and, by those of the proof of Proposition \ref{prop:cstiv},
$\widehat{\Sigma}(\Lambda)- \widehat{\Sigma}(\widehat{\Lambda})\le \widehat{\rho}^{ZX}|\widehat\Delta'|_1$. As a result, $\overline{\widehat\Delta'}\in K_{S(\Lambda)}'\subseteq\widehat{K}_{S(\Lambda)}'$ and, using the definition of $\widehat\kappa_{\ell^1_{S(\Lambda)},S(\Lambda)}'$ and of the objective function in \eqref{Mhat} in the first display and \eqref{stop2CCB} in the second display,
\begin{align*}
&\widehat{\nu}\le\frac{\widehat{\rho}^{ZX}\left|\widehat\Delta'\widehat{\Psi}^{\top}\right|_{\infty}}{\lambda\widehat\kappa_{\ell^1_{S(\Lambda)},S(\Lambda)}'}
+\widehat{\Sigma}\left(\Lambda\right),\quad 
\widehat{\nu}+\widehat{\Sigma}\left(\Lambda\right)\le2\widehat{\Sigma}\left(\Lambda\right)\gamma\left(\frac{\underline{r}'_n\widehat{\rho}^{ZX}}{\lambda\widehat\kappa_{\ell^1_{S(\Lambda)},S(\Lambda)}'}\right),\\
&\left|\widehat\Delta'\widehat{\Psi}^{\top}\right|_{\infty}
\le 2\underline{r}'_n\widehat{\Sigma}\left(\Lambda\right)\gamma\left(\frac{\underline{r}'_n\widehat{\rho}^{ZX}}{\lambda\widehat\kappa_{\ell^1_{S(\Lambda)},S(\Lambda)}'}\right).
\end{align*}
Let us now show the results of item \eqref{tCLambdaii}. Take $S\subseteq[d_{\varPhi}]\times[{d_Z}]$. We have
\begin{align*}
|\widehat\Delta_{S^c}'|_1&\le\left|\widehat\Delta_S'\right|_1+2\left|\Lambda_{S^c}\bold{D}_{\bold{Z}}^{-1}\right|_1+\lambda\left|\widehat\Delta'\right|_1
\end{align*}
and distinguish the two cases\\
Case 1: $2|\Lambda_{S^c}\bold{D}_{\bold{Z}}^{-1}|_1\le|\widehat\Delta_S'|_1$ for which the rest is usual,\\
Case 2: $2|\Lambda_{S^c}\bold{D}_{\bold{Z}}^{-1}|_1>|\widehat\Delta_S'|_1$ for which we have
$$ \left|\widehat\Delta'\right|_1=\left|\widehat\Delta_{S^c}'\right|_1+\left|\widehat\Delta_S'\right|_1\le 
2\frac{3+\lambda}{1-\lambda} \left|\Lambda_{S^c}\bold{D}_{\bold{Z}}^{-1}\right|_1,$$
hence
$$ \left|\overline{\widehat\Delta'}\right|_1\le 
\frac{2}{1_{n}}\frac{3+\lambda}{1-\lambda}  \left|\Lambda_{S^c} D_{Z}^{-1}\right|_1.\quad\quad\quad\quad\quad\square$$

\noindent{\bf Proof of theorems \ref{th:nonvalid} and \ref{th:nonvalid1}.} We work on $\mathcal{G}_{\perp}\cap\underline{\mathcal{G}}_{\not\perp}$ and denote by $\widetilde{\Delta}\triangleq\bold{D}_{\bold{Z}}(\widehat{\theta} -\theta)$.\\ 
First, we show that $\theta\in\widehat{\mathcal I}_{\not\perp}(\underline{r}^{\not\perp}_n,\widehat{\Sigma}_{\not\perp}(\beta,\theta))$ by the following computations
\begin{eqnarray*}\label{eq:proof:nv}
\left|\bold{D}_{\bold{Z}}\left(\mathbb{E}_n[ZU(\widehat{\beta})]-\theta\right)_{ S_{\perp}^c}\right|_\infty &\le& \left|\bold{D}_{\bold{Z}}\left(\mathbb{E}_n[ZU(\beta)]-\theta\right)_{ S_{\perp}^c}\right|_\infty +\left|\left(\widehat{\Psi}\bold{D}_{\bold{X}}^{-1}\left(\widehat{\beta}-\beta\right)\right)_{ S_{\perp}^c}\right|_\infty\\
&\le&\underline{r}^{\not\perp}_n \widehat{\Sigma}_{\not\perp}\left(\beta,\theta\right) +\widehat{\delta}.
\end{eqnarray*}
The second constraint in the definition of $\widehat{\mathcal I}_{\not\perp}(\underline{r}^{\not\perp}_n,\widehat{\Sigma}_{\not\perp}(\beta,\theta))$ is satisfied because, by convexity, 
$\widehat{\Sigma}_{\not\perp}(\widehat{\beta},\theta) \le\widehat{\Sigma}_{\not\perp}\left(\beta,\theta\right)+ \widehat{\delta}^{\Sigma}$.
Now, because $\theta\in\widehat{\mathcal I}_{\not\perp}(\underline{r}^{\not\perp}_n,\widehat{\Sigma}_{\not\perp}(\beta,\theta))$
and $(\widehat{\theta},\widehat{\widetilde{\sigma}})$ minimizes \eqref{def:STIV_est_nonvalid}, 
\begin{equation}\label{eq:proof:nv6}
\left|\widetilde{\Delta}_{S\left(\theta\right)^c}\right|_1 \le
\left|\widetilde{\Delta}_{S\left(\theta\right)}\right|_1+
\widetilde{c}\left(\widehat{\Sigma}_{\not\perp}\left(\beta,\theta\right)-\widehat{\widetilde{\sigma}}\right).
\end{equation}
Similar to the proof of \eqref{dubovicC}, using the second constraint in the definition of $\widehat{\mathcal I}_{\not\perp}(\underline{r}^{\not\perp}_n,\widehat{\widetilde{\sigma}})$, 
\begin{equation}\label{eq:proof:nv6a}
\widehat{\Sigma}_{\not\perp}\left(\beta,\theta\right)-\widehat{\widetilde{\sigma}}\le\widehat{\Sigma}_{\not\perp}\left(\beta,\theta\right)-\widehat{\Sigma}_{\not\perp}\left(\widehat\beta,\widehat\theta\right)+\widehat{\delta}^{\Sigma}\le  2\widehat{\delta}^{\Sigma}+\underline{r}^{\not\perp}_n\left|\widetilde{\Delta}_{S_{\perp}^c}\right|_1.
\end{equation}
This and (\ref{eq:proof:nv6}) yield
$$
\left|\widetilde{\Delta}_{S\left(\theta\right)^c}\right|_1 \le
\left|\widetilde{\Delta}_{S\left(\theta\right)}\right|_1+
\widetilde{c}\underline{r}^{\not\perp}_n\left|\widetilde{\Delta}_{S_{\perp}^c}\right|_1 + 2\widetilde{c}\widehat{\delta}^{\Sigma}
$$
and, equivalently,
\begin{equation}\label{eq:proof:nv7}
\left|\widetilde{\Delta}_{S\left(\theta\right)^c}\right|_1 \le \frac{1+\widetilde{c}\underline{r}^{\not\perp}_n}{1-\widetilde{c}\underline{r}^{\not\perp}_n}\left|\widetilde{\Delta}_{S\left(\theta\right)}\right|_1 +
\frac{2\widetilde{c}}{1-\widetilde{c}\underline{r}^{\not\perp}_n}\widehat{\delta}^{\Sigma}.
\end{equation}
Next, using the first constraint in the definition of $\widehat{\mathcal I}_{\not\perp}(\underline{r}^{\not\perp}_n,\widehat{\widetilde{\sigma}})$ and $\widehat{\mathcal I}_{\not\perp}(\underline{r}^{\not\perp}_n,\widehat{\Sigma}_{\not\perp}(\beta,\theta))$, 
we find
\begin{align}
\left|\bold{D}_{\bold{Z}}\left(\widehat{\theta}-\theta\right)_{S_{\perp}^c}\right|_\infty &\le \left|\bold{D}_{Z}\left(\mathbb{E}_n[ZU(\widehat{\beta})]-\widehat{\theta}\right)_{S_{\perp}^c}\right|_\infty+ \left|\bold{D}_{\bold{Z}}\left(\mathbb{E}_n[ZU(\widehat\beta)]-\theta\right)_{S_{\perp}^c}\right|_\infty\notag\\ 
&\le \underline{r}^{\not\perp}_n\left(\widehat{\widetilde{\sigma}} +\widehat{\Sigma}_{\not\perp}\left(\beta,\theta\right)\right)
+2\widehat{\delta}\label{ethetasup}.
\end{align}
This and \eqref{eq:proof:nv6a} yield
\begin{equation}\label{eq:proof:nv8}
\left|\widetilde\Delta_{S_{\perp}^c}\right|_\infty \le
\underline{r}^{\not\perp}_n\left(2\overline{\widetilde{\sigma}} + \underline{r}^{\not\perp}_n\left|\widetilde\Delta_{S_{\perp}^c}\right|_1+\widehat{\delta}^{\Sigma}\right)
+2\widehat{\delta}.
\end{equation}
On the other hand, (\ref{eq:proof:nv7}) implies 
\begin{align}
\left|\widetilde\Delta_{S_{\perp}^c}\right|_1 
&\le
\frac{2\left|S\left(\theta\right)\right|}{1-\widetilde{c}\underline{r}^{\not\perp}_n}\left|\widetilde{\Delta}_{S_{\perp}^c}\right|_\infty +
\frac{2\widetilde{c}\widehat{\delta}^{\Sigma}}{1-\widetilde{c}\underline{r}^{\not\perp}_n}.\label{eq:proof:nv9}
\end{align}
\eqref{eq:th:nonvalidCSinfty} 
follows by simple manipulations of (\ref{eq:proof:nv8})-(\ref{eq:proof:nv9}). As before, we obtain 
\begin{equation}\label{eq:proof:nv10}
{\widehat{\widetilde{\sigma}}}\le|\widetilde{\Delta}_{S\left(\theta\right)}|_1/\widetilde{c} + 
\widehat{\Sigma}_{\not\perp}\left(\beta,\theta\right)
\le \left|S\left(\theta\right)\right||\widetilde{\Delta}_{S_{\perp}^c}|_{\infty}/\widetilde{c}+\widehat{\Sigma}_{\not\perp}\left(\beta,\theta\right),
\end{equation}
which, together with \eqref{ethetasup}, 
yield the following bound used to obtain Theorem \ref{th:nonvalid1}
\begin{align*}
\left|\bold{D}_{\bold{Z}}\left(\widehat{\theta}-\theta\right)_{S_{\perp}^c}\right|_\infty
&\le 2\gamma\left(\underline{r}^{\not\perp}_n\left|S\left(\theta\right)\right|/\widetilde{c}\right)
\left(\underline{r}^{\not\perp}_n\widehat{\Sigma}_{\not\perp}(\beta,\theta)+\widehat{\delta}\right).
\end{align*}
The rest is as before.\hfill$\square$\vspace{0.2cm}

\noindent{\bf Proof of propositions \ref{ThoracleSE} and \ref{ThoracleSEO}.} Take $\beta$ in $\mathcal{I}$, set $\widehat{\Delta}\triangleq\bold{D}_{\bold{X}}^{-1}(\widehat{\beta}-\beta)$,  and work on $\mathcal{G}\cap\mathcal{G}_{A1}$. We have, for $g\in[{d_G}]$, using the triangle inequality in the second and fourth display, and the definition of $\mathcal{G}$ and the Cauchy-Schwartz inequality in the third,
\begin{align*}
\left|\bold{D}_{\bold{Z}}\mathbb{E}_n[ZU_g(\beta)]\right|_{\infty}
&\le\left|\bold{D}_{\bold{Z}}\mathbb{E}_n[ZW_g(\beta)]\right|_{\infty}+\left|\bold{D}_{\bold{Z}}\mathbb{E}_n[ZV_g(\beta)]\right|_{\infty}\\
&\le r_n\E_n[W_{g}(\beta)^2]^{1/2}+\E_n[V_{g}(\beta)^2]^{1/2}\\
&\le r_n\widehat{\sigma}_g(\beta)+(r_n+1)\E_n[V_g(\beta)^2]^{1/2}\le r_n\widehat{\sigma}_g(\beta)+(r_n+1)\widehat{v}_g.
\end{align*}
Hence, $\beta\in\widehat{\mathcal{I}}_{E}\left(r_n,\widehat{\sigma}(\beta)\right)$ and
\begin{equation}\label{PsiE}
\left|\widehat{\Psi}\widehat{\Delta}_{\cdot,g}\right|_{\infty}\le r_n\left(\widehat{\sigma}_g+\widehat{\sigma}_g(\beta)\right)
+2(r_n+1)\sqrt{1+\tau_n}v_{g,d_X}.
\end{equation}
Moreover, by the inverse triangle inequality, we have 
\begin{align*}
\widehat{\sigma}_g(\beta)&\ge\E_n[W_g(\beta)^2]^{1/2}-\E_n[V_g(\beta)^2]^{1/2}\ge \sqrt{1-\tau_n}\sigma_{W_g(\beta)}-\sqrt{1+\tau_n}v_{g,d_X}.
\end{align*}
Hence, by convexity, we have
\begin{align}
\widehat{\sigma}_g(\beta)
-\widehat\sigma_g\left(\widehat\beta\right)
&\le\min\left(
r_n+(r_n+1)\sqrt{1+\tau_n}v_{g,d_X}\gamma\left(\sqrt{1-\tau_n}\sigma_{W_g(\beta)}-\sqrt{1+\tau_n}v_{g,d_X}\right),1\right)|\widehat{\Delta}_{S_I}|_1\notag\\
&\quad+\left|\widehat{\Delta}_{S_I^c}\right|_1\notag\\
&\le\min\left(r_n+(r_n+1)\max\left(0,1_{n}\frac{\sigma_{ W_g(\beta)}}{v_{g,d_X}}-1\right)^{-1},1\right)\left|\left(\widehat{\Delta}_{S_I}\right)_{\cdot,g}\right|_1+\left|\left(\widehat{\Delta}_{S_I^c}\right)_{\cdot,g}\right|_1\notag\\
&\le r_n(\beta)\left|\left(\widehat{\Delta}_{S_I}\right)_{\cdot,g}\right|_1+\left|\left(\widehat{\Delta}_{S_I^c}\right)_{\cdot,g}\right|_1\label{pourfin}.
\end{align} 
Hence we obtain the first inequality in \eqref{ThoracleSEi}. Denoting by $\widehat{\Psi}_{X}\triangleq D_X\mathbb{E}_n[XX^{\top}]D_X$, the second inequality comes from 
$
\mathbb{E}_n[(X^{\top}(\widehat{\beta}_{\cdot,g}-\beta_{\cdot,g}))^2]\le|\widehat{\Psi}_X\Delta_g|_{\infty}\left|\Delta_g\right|_{1}.
$\\ 
The third inequality comes from 
$
|\widehat{\sigma}_g(\widehat\beta)-\sigma_{W_g(\beta)}|\le\mathbb{E}_n[(X^{\top}(\widehat{\beta}_{\cdot,g}-\beta_{\cdot,g}))^2]^{1/2}+|\sigma_{W_g(\beta)}-\widehat{\sigma}_g(\beta)|
$
and $\max\left(\sqrt{1+\tau_n}-1,1-\sqrt{1-\tau_n}\right)\le \tau_n$.
By definition of the estimator, we have
\begin{align*}
\hspace{-.5cm}\left|\widehat{\Delta}_{S^c\cap S_Q}\right|_1
\le &\left|\widehat{\Delta}_{S\cap S_Q}\right|_1+2\left|\bold{D}_{\bold{X}}^{-1}\beta_{S^c\cap S_Q}\right|_1\\
&\hspace{-1.6cm}+c\sum_{g\in[{d_G}]}\left(\min\left(r_n+(r_n+1)\max\left(0,\frac{1_{n}\sigma_{W_g(\beta)}}{v_{g,d_X}}-1\right)^{-1},1\right)\left|\left(\widehat{\Delta}_{S_I}\right)_{\cdot,g}\right|_1+\left|\left(\widehat{\Delta}_{S_I^c}\right)_{\cdot,g}\right|_1\right)\\
\le&\left|\widehat{\Delta}_{S\cap S_Q}\right|_1+2\left|\bold{D}_{\bold{X}}^{-1}\beta_{S^c\cap S_Q}\right|_1+c\left(r_n(\beta)
\left|\widehat{\Delta}_{S_I}\right|_1+\left|\widehat{\Delta}_{S_I^c}\right|_1\right)\\
\sum_{g=1}^{d_G}\left|\widehat{\Psi}\widehat{\Delta}_{\cdot,g}\right|_{\infty}\le& r_n\sum_{g=1}^{d_G} \left(\widehat{\sigma}_g+\widehat{\sigma}_g(\beta)\right)+2(r_n+1)\sqrt{1+\tau_n}\sum_{g=1}^{d_G} v_{g,d_X}\\
\le& \frac{r_n}{c}\left(\left|\bold{D}_{\bold{X}}^{-1}\beta_{S_Q}\right|_1-\left|\bold{D}_{\bold{X}}^{-1}\widehat{\beta}_{S_Q}\right|_1\right)+2 r_n\sum_{g=1}^{d_G}\widehat{\sigma}_g(\beta)+2(r_n+1)\sqrt{1+\tau_n}\sum_{g=1}^{d_G} v_{g,d_X}.
\end{align*}
The second inequality from \eqref{ThoracleSEii} is obtained in a similar manner as in the proof of Proposition \ref{tCLambda}. 
The last statement is obtained using that, by similar arguments as those leading to 
\eqref{eNPIV0}, 
\begin{equation*}
\mathbb{E}_n\left[\left(X^{\top}\left(\widehat{\beta}_{\cdot,g}-\beta_{\cdot,g}\right)\right)^2\right]\le 
\left|D_X^{-1}\left(\widehat{\beta}_{\cdot,g}-\beta_{\cdot,g}\right)\right|_{1}^2|\widehat{\Psi}_X|_{\infty}.\quad\quad\quad\quad\quad\quad\square
\end{equation*}

\end{document}